\documentclass[leqno,11pt,oneside]{amsart}

 \usepackage[english]{babel}

\usepackage{xcolor}
\usepackage{graphicx}
\usepackage{caption}
\usepackage{relsize}
\usepackage{amsthm}
\usepackage{tikz}
\usetikzlibrary{decorations.pathreplacing}

\usepackage{enumitem}
\usepackage{amsmath,amsfonts,amssymb,amsthm}
\usepackage{soul}

\usepackage{bm}

\newtheorem*{theorem*}{Theorem}
\newtheorem{teo}{Theorem}[section]
\newtheorem{prop}[teo]{Proposition}
\newtheorem{ddef}[teo]{Definition}
\newtheorem{example}[teo]{Example}
\newtheorem{cor}[teo]{Corollary}

\newtheorem{lem}[teo]{Lemma}
\newtheorem*{cor*}{Corollary}

\newtheorem{assertion}{Assertion}

\newtheorem*{lem*}{Lemma}

\newtheorem*{teorA'}{Theorem A'}

\newtheorem{teoRM}{Theorem}

\newtheorem*{fact*}{Fact}

\newtheorem*{assertion*}{Assertion}

\theoremstyle{definition}

\newcommand{\val}{{\rm Val}}
\newcommand{\cont}{\mathrm{Cont}}

\newcommand{\ord}{{\rm ord}}

\newcommand{\C}{\mathbb{C}}
 
\newcommand{\R}{\mathbb{R}}

\newcommand{\N}{\mathbb{N}}
\newcommand{\Z}{\mathbb{Z}}
\newcommand{\Q}{\mathbb{Q}}

\newcommand{\PP}{\mathrm{P}}

\newcommand{\ma}{{\rm m}}
\newcommand{\Top}{{\rm Top}}
\newcommand{\Bot}{{\rm Bot}}

\newcommand{\Ne}{\mathcal{N}}

\newcommand{\dd}{\textup{d}}

\newcommand{\A}{\mathcal{A}}
\newcommand{\B}{\mathcal{B}}

\newcommand{\D}{\mathcal{D}}
\newcommand{\RR}{\mathcal{R}}

\newcommand{\F}{\mathcal{F}}

\newcommand{\E}{\mathcal{E}}

\newcommand{\supp}{\mathcal{S}}

\newcommand{\cl}[1]{\mbox{$\mathcal{#1}$}}

\newcommand{\sep}{{\rm Sep}}
\newcommand{\Sep}{{\rm Sep}}

\newcommand{\oo}{{\rm o}}
\newcommand{\In}{{\rm in}}

\newcommand{\hash}{\scriptscriptstyle \#}

\newcommand*\xbar[1]{ %
   \hbox{ %
     \vbox{%
       \hrule height 0.3pt 
       \kern0.35ex
       \hbox{%
         \kern-0.1em
         \ensuremath{#1}%
         \kern-0.1em
       }%
     }%
   }%
}

\newcommand*\xxbar[1]{%
   \hbox{%
     \vbox{%
       \hrule height 0.3pt 
       \kern0.4ex
       \hbox{%
         \kern-0.1em
         \ensuremath{#1}%
         \kern-0.1em
       }%
     }%
   }%
}

\makeindex

\usepackage{lineno}

\begin{document}



\thanks{First author  supported by a grant   CNPq-Brazil PQ-306011/2023-9.
Second author supported by the Spanish grant PID2023-149508NB-I00 funded by MICIU/AEI /10.13039/501100011033 and by FEDER, UE.
Third author   supported by   a grant FAPEMIG-CNPq APQ-06573-24. He also thanks Vicerrectorado de Investigación y Transferencia and Agencia Canaria de Investigación, Innovación y Sociedad de la Información del Gobierno de Canarias for supporting his visit to Tenerife in 2025.
First and third  authors  supported by CNPq Projeto Universal 408687/2023-1 and by   CAPES/PRINT-UFMG PRPG08/2024.}

\setcounter{section}{0}
\setcounter{teo}{0}
\setcounter{exe}{0}

\title{Newton non-degenerate holomorphic foliations and their separatrices}

\author{Arturo Fern\'andez-P\'erez}
\address{Departamento de Matem\'atica ---  Universidade Federal de Minas Gerais}
\curraddr{Av. Ant\^onio Carlos 6627 --- 31270-901 --- Belo Horizonte, BRAZIL.}
\email{fernandez@ufmg.br}

\author{Evelia R. Garc\'{i}a Barroso}
\address{Dpto. Matem\'{a}ticas, Estad\'{\i}stica e Investigaci\'on Operativa\\
Instituto Universitario de Matem\'aticas y Aplicaciones (IMAULL)\\
Universidad de La Laguna.} 
\curraddr{ Apartado de Correos 456 --- 38200 --- La Laguna, Tenerife, SPAIN.}
\email{ergarcia@ull.es}

\author{Rog\'erio Mol}
\address{ Departamento de Matem\'atica ---  Universidade Federal de Minas Gerais}
\curraddr{Av. Ant\^onio Carlos 6627 --- 31270-901 --- Belo Horizonte, BRAZIL.}
\email{rmol@ufmg.br}

\maketitle

\noindent\hfill
\begin{minipage}{0.65\textwidth}
\begin{flushright}
    \textit{Dedicated to the memory of Professor Arkadiusz P\l oski, whose work was a profound inspiration to this research.}
\end{flushright}
\end{minipage}
\vspace{2em}

\begin{abstract}
In this paper we study germs of holomorphic foliations in the complex plane that are Newton non-degenerate.
This is an algebraic condition attached to the Newton polygon of the foliation,
analogous to  the one defined by Kouchnirenko for analytic curves
and, more generally, hypersurfaces. As in this case, the condition of Newton non-degeneracy
for foliations can be identified through the comparison of two numerical invariants: the Newton number
--- a computation involving simple geometric data of the Newton polygon --- and the Milnor number.
The Newton number is never greater than the Milnor number, and the foliation is Newton
non-degenerate if and only if they are equal.
The assumption that a foliation is Newton non-degenerate imposes conditions on its set of separatrices.
We prove that, under this hypothesis, the foliation has no weak separatrices and, hence, all its separatrices are analytic. 
We also give a description of the set of separatrices
of a Newton non-degenerate foliation, which is  sufficient to characterize the foliation as
Newton non-degenerate when it is of generalized curve type.
This  allows us to prove that the property of being Newton non-degenerate is
topologically invariant for foliations of generalized curve type.
\end{abstract}

 \medskip \medskip

  \section{Introduction}
\label{sec-introduction}

This article presents a geometric study of germs of \emph{singular holomorphic foliations} in the complex plane --- referred to hereafter simply as
\emph{holomorphic foliations} or \emph{foliations} ---
using their \emph{Newton polygons} as a main tool. For a  foliation $\F$, defined at $(\C^{2},0)$ by a germ of holomorphic $1$-form $\omega$, our    main purpose
is to obtain information on its  invariant irreducible formal curves --- its \emph{separatrices} --- or,
equivalently, on the formal solutions of the  analytic ordinary differential equation corresponding to   $\omega$. We have as a fundamental  assumption
an algebraic condition of generic nature called \emph{Newton non-degeneracy}. Some basic material on foliations is presented in Section
\ref{sec-holomorphic-foliations} and the main definitions concerning  the notions of  Newton polygon and   Newton non-degeneracy
are given in Section  \ref{sec-kouchnirenko-theorem}. In particular, most terms in italics  in this Introduction are  explained
in these two sections.

The study of algebraic curves and their singularities has long been a central theme in algebraic geometry. Classical investigations by Newton already revealed that the local behavior of plane curves can be analyzed through expansions of algebraic equations, a method that was later formalized through the Newton-Puiseux theorem. These ideas provided early tools for understanding the structure of singular points of algebraic curves.
The concept of {\it non-degeneracy} with respect to Newton polyhedra was studied by Kouchnirenko in the 1970s (see \cite{kouchnirenko1976}). In this framework, one associates to a polynomial, to a power series or, more generally, to a Laurent series, its Newton polyhedron. A function is said to be non-degenerate if, for every compact face of its Newton polyhedron, the corresponding truncated function defines a hypersurface without singularities in the algebraic torus. One of the key results of Kouchnirenko shows that, under this non-degeneracy condition, several important invariants of singularities can be computed directly from the combinatorial geometry of the Newton polyhedron. In particular, the Milnor number of an isolated hypersurface singularity can be expressed in terms of the
volumes of the region bounded by the Newton polyhedron
and the coordinate hyperplanes, and of its sections by coordinate subspaces.
These ideas were further developed by Varchenko and Khovanskii, who established deep relations between Newton polyhedra, monodromy of singularities, and toric geometry \cite{Varchenko1976, Khovanskii1977}. In particular, Varchenko related the zeta function of monodromy to the Newton polyhedra, while Khovanskii developed systematic methods for studying algebraic varieties using Newton polyhedra and toroidal techniques.

The  use of the Newton polygon and the Newton-Puiseux method in the investigation of solutions of analytic ordinary differential equations
was first proposed by  H. Fine  in 1889 \cite{fine1889} and   later presented in E. Ince's book \cite{ince1944}.
 At the end of the 20th century,
J. Cano, in his Ph.D. thesis (see the corresponding article \cite{jcano1993}) developed and refined this method, providing an algorithmic proof for the existence of an analytic separatrix for a holomorphic foliation (the so-called \emph{Separatrix theorem}, whose original proof was given in \cite{camacho1982}).
These fundamental works support the idea that the Newton polygon and the  algorithmic processes founded on it
encode substantial geometric  information of the foliation itself. Some works   on local holomorphic foliations,
published in the subsequent  years, illustrate this perception.
For example, the works  \cite{rouille1999, corral2003, nancy2018, FS-GB-SM2022, mariana2025}   employ the Newton polygon to characterize   families of foliations
(of \emph{generalized curve type},
of \emph{second type}) and to  study   their polar curves.
Another point of view appears in \cite{jcano2024}, where the information collected in the Newton polygon is used to bound from  above
topological invariants of an analytic solution.

The  concept  of Newton non-degeneracy for  real analytic planar vector fields  was proposed   by M. Brunella and M. Miari in
\cite{brunella1990}. Having it as a hypothesis, they proved  that a local vector field which is not a center-focus is topologically equivalent to its principal part, that is, the polynomial vector field obtained by collecting the terms lying on the Newton polygon.
A. \`Esterov, in \cite{esterov2006}, used Newton polyhedra and the notion of Newton non-degeneracy in order to calculate indices of holomorphic  $1$-forms.
More recently, in \cite{molinasamper2022},  B. Molina-Samper studied questions involving the reduction of singularities of Newton non-degenerate foliations.
In dimension two, for instance, she proves that a generalized curve type foliation is Newton non-degenerate
if and only if it has a toric pre-reduction of singularities, that is, if and only if the foliation
can be transformed
into one having at most \emph{pre-simple singularities} by a sequence of quadratic blow-ups locally expressed as
a monomial map. Our results here go beyond this point: saddle-node singularities
are admitted (see, for instance, Theorem \ref{th-weak} below) and, as far as a full reduction of singularities is considered, that is, one leading to \emph{simple singularities}, it can
be of non-toric type (see Theorems \ref{teo-newton-curve1} and \ref{teo-newton-curve2} and also Example \ref{ex-dicfamily-g2}).

The concept of  Newton  non-degeneracy is linked to the Newton polygon, inheriting thus the dependency
on the choice of analytic coordinates.
These are assumed to be \emph{convenient}, which means that the Newton polygon touches both the vertical and the horizontal axes.
For a foliation, Newton  non-degeneracy is
defined similarly to   analytic curves, expressed by means of  polynomial equations involving  monomial terms of the defining
$1$-form $\omega$ associated with  points with integer coordinates lying on  the Newton polygon (see the formulas in \eqref{eq-degeneracy} and Definition \ref{def-Newton-nondeg}).

Just as it happens in Kouchnirenko's theorem for analytic curves (Corollary \ref{cor-mu-nu-curves}), for a foliation $\F$, Newton non-degeneracy can be identified through the comparison of two numerical invariants.
One of them is the celebrated \emph{Milnor number} $\mu_{0}(\F)$, a topological invariant  of the foliation (\cite{camacho1984}, see also
formula \eqref{def-milnor-number} below). The other is the
 \emph{Newton number}, defined as  $\nu_{0}(\F) = 2 \A - \ell - h +1$, where $\A$ is the area of the region delimited
by the Newton polygon and the coordinate axes, while $\ell$ and $h$ are the lengths of its   projections
onto the horizontal and vertical axes (formulas
 \eqref{def-newton-number}). In convenient coordinates, we always have $\nu_{0}(\F) \leq \mu_{0}(\F)$, with
 equality holding  if and only if the foliation is Newton non-degenerate. This is Kouchnirenko's theorem for foliations, our  Theorem \ref{teo-fol-nondeg2}.
Its proof is carried out in Section \ref{sec-kouchnirenko-theorem}
and relies on the study of Newton non-degeneracy of pairs of
functions carried out by A. P\l{}oski in \cite{ploski1990}, which is
presented in  \S\ref{subsec-Newton-pairs}.

The property of Newton  non-degeneracy imposes conditions on the set of separatrices of a foliation.
This idea is the essence of the remaining material in this article. 
Theorem \ref{th-weak}, in Section \ref{sec-well-oriented-saddle-node},   a first result
in this direction, asserts that  the existence of a \emph{weak separatrix} --- that is, an invariant formal branch corresponding,
in the \emph{reduction of singularities} of the foliation, to a separatrix associated with the zero eigenvalue of
a saddle-node singularity --- implies that $\F$ is Newton degenerate.
We should emphasize that this    phenomenon actually concerns the existence of weak separatrices and not a more general
condition involving the
 existence of saddle-node singularities in the reduction of singularities (Examples \ref{ex-newton-deg} and \ref{ex-newton-ndeg}).
Considering that a  purely formal separatrix is a weak separatrix, we have as a consequence that,
if $\nu_{0}(\F) = \mu_{0}(\F)$ for some system of convenient coordinates, then
all separatrices of $\F$ are analytic. This interesting  numerical criterion for analyticity,
 which is sufficient but not necessary
(see Example \ref{ex-newton-deg}),
 is a condition
that can be read only on the  initial terms of the defining $1$-form $\omega$.

A Newton non-degenerate analytic curve   is characterized in \cite{garciabarroso2007} 
(see also
\cite{garciabarroso2009}) as an object  called there
 \emph{Newton curve}  (Definition \ref{def-Newtoncurve}).
This concept is invariant under   topological equivalences and,   as a consequence, for analytic curves, the property of Newton non-degeneracy is   topologically invariant.
For foliations, results of a similar nature hold,  this time involving the set of separatrices.
For this purpose, in Section \ref{sec-newton-germs} we establish the concept of  \emph{Newton set of separatrices}
(Definition
\ref{def-Newtontype-separatrices} and the pair of definitions preceding it), which, in a simplistic way, is a Newton curve, corresponding to the \emph{isolated separatrices}
of $\F$, combined with some special families of curves formed by its
\emph{dicritical separatrices},  studied in
Section \ref{sec-dicritical-family}.
Two meaningful aspects  are present in both Newton curves and Newton sets of separatrices:
the topological types of their branches are simple and there is a minimality condition on the contact of different branches.
We finish Section \ref{sec-newton-germs} by  proving Theorem \ref{teo-newton-curve1}, which asserts that a Newton non-degenerate
foliation has a Newton set of separatrices.
In Section \ref{sec-generalizecurve-nondeg} we prove  Theorem \ref{teo-newton-curve2}, a partial converse result: a foliation  of generalized curve type having a Newton set of separatrices is Newton non-degenerate.
The proof relies on the notion of effective foliation
\cite{mariana2025}, for which the
\emph{balanced divisor of separatrices} play the role of a curve
of reference.
We remark that the hypothesis of generalized curve type cannot be dropped
(Example \ref{ex-newton-set-not-enough}).
Finally, using the fact that the concept of Newton set of separatrices is invariant under topological equivalences,
we prove in  Theorem \ref{teo-topological} that, for foliations of generalized curve type, the property of being
Newton non-degenerate (in suitable coordinates) is   topologically invariant.

The strategy for proving Theorems \ref{th-weak} and
\ref{teo-newton-curve1}
consists in following a fixed separatrix through three types of
transformations and examining their effect on the Newton polygon:
\emph{blow-up maps}, which act towards the reduction of singularities of the foliation,
\emph{horizontalization procedures}, which bring oblique tangent lines to the horizontal position, and
\emph{ramification maps}, which   reduce all arguments to the case of smooth separatrices.
They are explained in the three subsections of Section \ref{sec-basic-transformations}. The final goal
is to reach elementary models called \emph{primary Newton polygons}, which correspond either to a regular foliation  or
to a singular foliation defined by a vector field  with non-nilpotent linear part --- this means, one having a \emph{pre-simple singularity},
which is   \emph{simple} if a non-resonance condition
 on the eigenvalues is satisfied. In this context, resonances, which may  correspond to
 dicritical blow-ups in the reduction of singularities,  bring about one peculiar  aspect of the study
 of  Newton polygons of   foliations: in addition to their   edges, their  vertices matter.
This contrasts with  what happens in the case of analytic curves, but adds a special flavor to this theory.

We conclude this introduction by expressing the hope that the results presented here shed further light on the role of the Newton polygon and contribute to its application in the study of singular holomorphic foliations.

 \section{Some concepts on holomorphic foliations}
\label{sec-holomorphic-foliations}
This section presents some basic concepts of the theory of singular holomorphic foliations.  As a supporting text for
the material here, we recommend the book \cite{cano-cerveau-deserti2013}.
A germ of \emph{singular holomorphic foliation} at $(\C^2,0)$ is the object $\F$ defined by a germ of holomorphic $1$-form $\omega$
with isolated  singularity at the origin, up to  multiplication by a unit in $\C\{x,y\}$. As mentioned in the Introduction,
throughout this text, we will refer to this object only as \emph{holomorphic foliation} or simply as \emph{foliation}.
In a small punctured neighborhood of $0 \in \C^{2}$, the distribution of directions determined by  $\omega$ defines a regular holomorphic foliation, which engenders the decomposition of this set into a disjoint
union of immersed analytic curves called \emph{leaves}.
If $\F$ is defined by  $\omega = A(x,y) \dd x + B(x,y) \dd y$, where $A, B \in \C\{x,y\}$ are relatively prime,
then the germ of holomorphic vector field
${\bf v} = -B(x,y)  \partial / \partial x +  A(x,y)  \partial / \partial y$ induces the same
distribution of directions   and hence can also be used to define the foliation.
 A germ of irreducible formal curve (a formal \emph{branch}) $S$ at $(\C^2,0)$,
 defined by an irreducible $f \in \C[[x,y]]$, is a \emph{separatrix}
of $\F$
if $ \omega \wedge \dd f = f \eta$  for some formal $2$-form $\eta$.
If $S$ is analytic, that is, defined by $f \in \C\{x,y\}$, then $\eta$ is holomorphic and
$S \setminus \{0\}$ is a leaf of $\F$. In this case, we say that $S$ is an \emph{analytic separatrix}.
Analytic separatrices always exist by \cite{camacho1982}.
We will denote the family of all separatrices of $\F$ by $\Sep(\F)$.

The     \emph{order} or \emph{algebraic multiplicity}  $\ma_{0}(f)$  of a series $f \in \C[[x,y]]$ (or of a fractional power series, in any number of variables)  is the smallest degree of its non-zero terms.
The   \emph{order} or \emph{algebraic multiplicity} of the foliation $\F$ is defined
as $\ma_{0}(\F)= \min \{ \ma_{0}(A), \ma_{0}(B)\}$.

A  foliation $\F$ as above has a \emph{reduction of singularities}   \cite{seidenberg1968,camacho1982}.
This means that there exists  a finite composition of quadratic blow-ups
$\pi: (\tilde{M},\D) \to (\C^{2},0)$, where $\D = \pi^{-1}(0)$
is the exceptional divisor, formed by  a finite union of projective lines with normal crossings and $\tilde{M}$ is a germ of surface
around $\D$ which, by taking successive strict transform foliations of individual blow-ups,
 transforms $\F$ into a foliation
$\tilde{\F} = \pi^{*} \F$,  defined   on $\tilde{M}$, having only
\emph{simple singularities} --- a concept to be explained in the next paragraph.
Besides, it is asked that     $\pi$  desingularizes
 the set of separatrices of $\F$, meaning that all branches in $\Sep(\F)$ are taken by $\pi$ into disjoint smooth curves, all of them transverse
 to the exceptional divisor $\D$. A reduction of singularities  is not unique, but we can fix a minimal
reduction of singularities, which is unique up to blow-up isomorphism. Along the text, unless otherwise stated,  reductions of singularities will be assumed to be minimal.

As mentioned above, one of the goals of the process of reduction of singularities is
to obtain reduced models called simple singularities.
A germ of holomorphic foliation has a \emph{pre-simple} singularity  if
it is defined by a germ of holomorphic vector field whose linear part is non-nilpotent.
Denoting by $\lambda_{1}, \lambda_{2}$ its eigenvalues, we say that the singularity is   \emph{simple}
if   one of them is zero or, when both are non-zero,
 they satisfy $\lambda = \lambda_{1} / \lambda_{2} \not\in \Q_{+}$.
The former is a \emph{saddle-node} singularity, the latter is a \emph{non-degenerate} one.
For a saddle-node singularity,  the foliation is defined in formal coordinates $(x,y)$  by Dulac's normal form:
\[ \omega =    y ^{k+1}  \dd x  -  x( 1+  \zeta y^{k} ) \dd y    , \]
where $k \geq 1$ and $\zeta \in \C$.
A non-degenerate singularity is given, in analytic coordinates, by a germ of holomorphic $1$-form of the type
\[ \omega =    y  A^{*}(x,y)  \dd x  -  x  B^{*}(x,y) \dd y    , \]
where $A^{*}, B^{*} \in \C\{x,y\}$.
A simple singularity
has exactly two separatrices, both smooth and transverse to each other,  the coordinate axes in the above normal forms, corresponding to the eigenspaces of the two eigenvalues.
In the non-degenerate case, both separatrices are analytic. In the saddle-node case, the one corresponding
to the non-zero eigenvalue is analytic --- the \emph{strong separatrix} --- whereas the one associated with
the zero eigenvalue can be a purely formal curve --- the \emph{weak separatrix}. The
weak separatrix is analytic if and only if Dulac's normal form is valid in analytic coordinates.

An irreducible component   $D \subset \D$ is classified as \emph{dicritical}, if it is not invariant by $\tilde{\F}$, or \emph{non-dicritical},
otherwise.
As explained in \S\ref{sec-blow-ups} below, this is equivalent to the fact that the blow-up map originating $D$ is dicritical or non-dicritical.
Each dicritical component   $D \subset \D$ is  fibered by leaves of $\tilde{\F}$.
These  leaves, except for those contained in $\D$, go down by $\pi$ into an infinite family of equisingular analytic branches of  $\Sep(\F)$, the separatrices \emph{associated with} $D$. Denote this family by $\sep(D)$.
All simple singularities of  $\tilde{\F}$ lie over non-dicritical components of $\D$.
They are   finite in number. Each simple singularity  positioned at a regular point of $\D$ (which we call \emph{trace point})
has exactly one separatrix outside $\D$. Such a separatrix is taken by $\pi$ into a branch of $\Sep(\F)$. Separatrices obtained in this
way are  called \emph{isolated separatrices}.
Note that a singularity of $\tilde{\F}$ lying at a singular point of $\D$
(which we call \emph{corner}) is in the intersection of two non-dicritical component, which correspond to its two separatrices.
We say that $\F$ is a \emph{dicritical} foliation if   the exceptional divisor
of its reduction of singularities contains some dicritical component. Otherwise, we say that $\F$ is \emph{non-dicritical}.
Equivalently, $\F$ is dicritical if and only if it has infinitely many separatrices.

If   a simple singularity of $\tilde{\F}$ is
a saddle-node, we say that it is a \emph{tangent saddle-node} if its weak separatrix is contained in $\D$. Otherwise, it is
said to be a \emph{well-oriented saddle-node}. In this case, the singularity necessarily  is a trace point of $\D$ and its weak
separatrix goes down to a branch of $\Sep(\F)$, which  can be purely formal and inherits the name \emph{weak separatrix}.
We note that all purely formal separatrices of  $\Sep(\F)$ are obtained in this way.

Two special families of foliations, characterized by their simple models in the reduction of singularities,
 will be considered in Section \ref{sec-generalizecurve-nondeg}. A foliation $\F$ is
of \emph{generalized curve  type} if it has no saddle-nodes in the exceptional divisor of its reduction of singularities and
it is of \emph{second type} if it contains no tangent saddle-nodes. Evidently, the family of generalized curve
foliations is contained in the family of second type foliations. Also, in
Section \ref{sec-generalizecurve-nondeg}, we will need the notion of
\emph{balanced divisor of separatrices} of the foliation $\F$,
an object that carries relevant combinatorial information
of the dicritical components of the reduction of singularities.
It is defined as a formal sum
\[ \B = \sum_{S \in \sep(\F)} a_{S} \cdot S , \]
where:
\begin{itemize}
  \item $a_{S} \in \{-1,0,1\}$;
  \item $a_{S} \neq 0$ for only finitely many $S \in \sep(\F)$;
  \item $a_{S} = 1$ whenever $S$ is an isolated separatrix;
  \item if $D \subset \D$ is a dicritical component, then $\sum_{S \in \sep(D)} a_{S} = 2 - \val(D)$,
\end{itemize}
with $\val(D)$ denoting the \emph{valence} of $D$, defined as the number of points of intersection of $D$ and the
components of $\D$ other than $D$ itself. In other words,  $\val(D)$ is equal to the number of corner points
of $\D$ lying on $D$.
Defining $\ma_{0}(\B) =  \sum_{S \in \sep(\F)} a_{S}   \ma_{0}(S)$, the following fact is proved in \cite{genzmer2007}:
$ \ma_{0}(\B) -1 \leq  \ma_{0}(\F)$,
with equality holding if and only if  $\F$ is a second type foliation. Actually, in \cite{genzmer2007}, the excess of $\ma_{0}(\F)$ over $\ma_{0}(\B)$ is quantified in terms of indices associated with tangent saddle-nodes, but we do not need
this information in the sequel.

\section{Kouchnirenko's theorem for foliations}
\label{sec-kouchnirenko-theorem}

\subsection{Newton non-degeneracy for pairs of functions}
\label{subsec-Newton-pairs}
The results in this subsection were published   by A. P\l{o}ski \cite{ploski1990}.
We have chosen to present them here in a form that aligns with our notation and suits our needs.
They will be used below to prove Theorem \ref{teo-fol-nondeg2}, concerning the inequality involving
Milnor and Newton numbers for a foliation.
We should mention that, although we establish definitions and state results for convergent power series, the material developed here
also works for formal power series.

Let us consider    $\C^{2}$ with a system of analytic coordinates $(x,y)$.
Given a set $\supp \subset \N \times \N  \subset \R^{2}$, where $\N = \{0,1,2, \ldots\}$ is the set of natural numbers, define $\Ne(\supp)$ to be the convex hull of the Minkowski sum $\supp + (\R_{\geq 0})^{2}$.
If $f = \sum_{i,j} f_{ij}x^{i}y^{j} \in \C\{x,y\}$ is a convergent power series, with $f(0,0) =0$,   then the \emph{support} of $f$ is the set
$\supp_f = \{(i,j) \in \N \times \N; f_{ij} \neq 0 \}$.
We then define the \emph{Newton diagram} of $f$ as the set $\Ne_{f} := \Ne(\supp_{f})$.
Evidently, this object depends on the analytic coordinates $(x,y)$.
If both $f(0,y)$ and $f(x,0)$ are non-zero, then
we say that  $(x,y)$ are  \emph{convenient} coordinates  or that $f$ itself is \emph{convenient} for
the coordinates $(x,y)$. This condition  means that $\Ne_{f}$ intersects both coordinate axes.
Thus, the boundary of $\Ne_{f}$ is formed by a compact polygonal line  and by
a pair of half-lines contained in the coordinate axes.  More generally, if $f$ is not convenient, the boundary of $\Ne_f$ still
consists of a compact polygonal line together with a pair of half-lines, each
contained either in a coordinate axis or in a line parallel to one of the
coordinate axes. In either case, we call the compact polygonal line, denoted
$\Delta_f$, the \emph{Newton polygon} of $f$.

When dealing with an edge $E$  of $\Delta_f$ we will preferably consider its  \emph{co-slope} or \emph{inclination}  $\kappa \in \Q_{>0}$, defined to
be the negative of the inverse of its slope.
Observe that, if ordered from the left to the right, the edges of $\Delta_f$ have growing co-slopes.
 If $E$ is an edge of $\Delta_f$ of co-slope $\kappa$, then its topmost   vertex will be denoted
 by $\Top(E)$, while its bottommost vertex by $\Bot(E)$.
If $\Top(E) = (i,j)$ and $\Bot(E) = (i',j')$,
then, for $E$, we define  its \emph{height} as $h(E) = j - j'$, its  \emph{length} as $\ell(E) = i' - i$
and its \emph{elevation} as  $e(E) = j'$.
Its  co-slope is
  $\kappa =  \frac{\ell(E) }{h(E)}$.


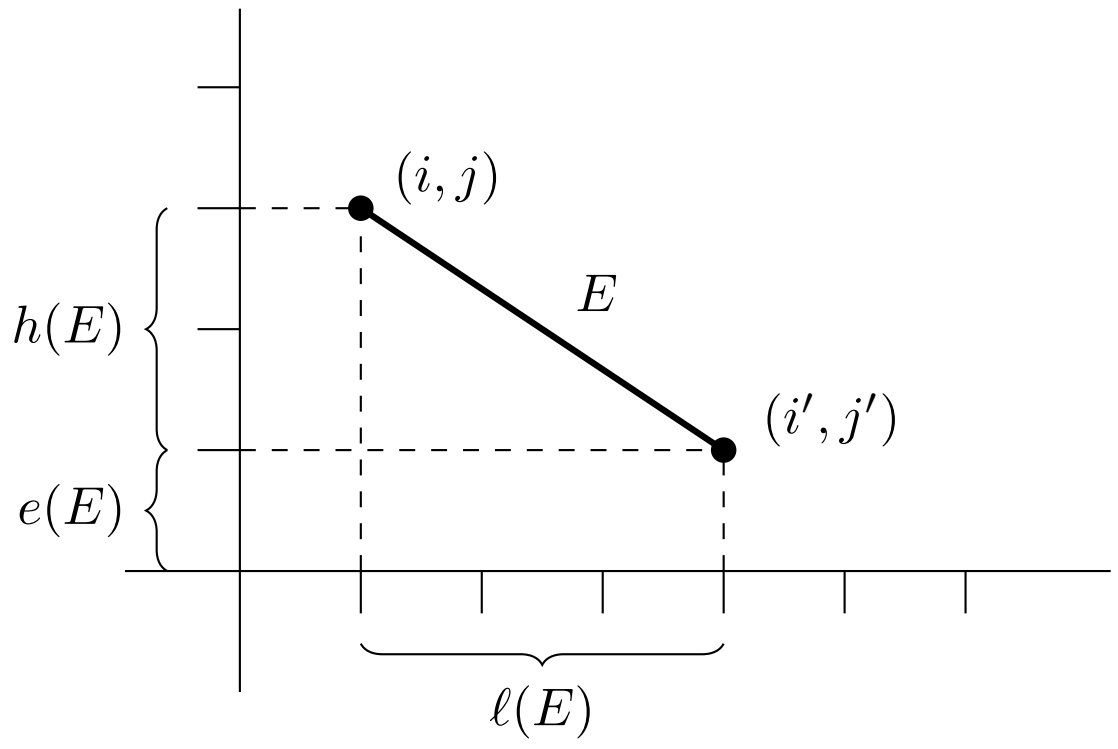
\begin{figure}[h]
\centering
\begin{tikzpicture}[scale=0.8]
  \draw (-0.95,0) -- (7.2,0);
  \draw (0,-1) -- (0,4.65);
  \foreach \x in {1,...,6} \draw (\x,0) -- (\x,-0.35);
  \foreach \y in {1,...,4} \draw (0,\y) -- (-0.35,\y);
  \draw[dashed] (0,3) -- (1,3);
  \draw[dashed] (0,1) -- (4,1);
  \draw[dashed] (1,0) -- (1,3);
  \draw[dashed] (4,0) -- (4,1);
  \draw[very thick] (1,3) -- (4,1);
  \fill (1,3) circle (3pt);
  \fill (4,1) circle (3pt);
  \node[right] at (1.1,3.25) {\small $(i,j)$};
  \node[right] at (4.15,1.25) {\small $(i',j')$};
  \node at (2.95,2.3) {\small $E$};
  \draw[decorate,decoration={brace,amplitude=4pt,mirror}]
    (-0.6,3) -- (-0.6,1) node[midway,left=4pt] {\small $h(E)$};
  \draw[decorate,decoration={brace,amplitude=4pt,mirror}]
    (-0.6,1) -- (-0.6,0) node[midway,left=4pt] {\small $e(E)$};
  \draw[decorate,decoration={brace,amplitude=4pt,mirror}]
    (1,-0.6) -- (4,-0.6) node[midway,below=4pt] {\small $\ell(E)$};
\end{tikzpicture}
\end{figure}

We also define
\[  h(\Ne_{f}) := \sum_{E \subset \Delta_{f}} h(E)
 \qquad \text{and} \qquad
\ell(\Ne_{f}) := \sum_{E \subset \Delta_{f}} \ell(E),
\]
where the sum runs over all edges of  $\Delta_{f}$.
If the coordinates $(x,y)$ are convenient, then
$h(\Ne_{f}) = (f,x)_{0}$ and  $\ell(\Ne_{f}) = (f,y)_{0}$,
where $(\, \cdot \, , \, \cdot \, )_{0}$ denotes the intersection number at $0 \in \C^{2}$.

For a pair of edges $E_1, E_2 \subset \Delta_{f}$, define their \emph{bracket} by
\[ [E_1, E_2] = \min \{ \ell(E_{1}) h(E_{2}), h(E_{1}) \ell(E_{2}) \} .\]
Observe that
\begin{equation}
\label{eq-bracket-coslope} \kappa(E_1) \leq \kappa(E_2) \quad \Leftrightarrow \quad \frac{\ell(E_1) }{h(E_1)} \leq \frac{\ell(E_2) }{h(E_2)}
\quad \Leftrightarrow \quad [E_1, E_2] = \ell(E_{1}) h(E_{2}).
\end{equation}
The \emph{principal part} of $f$ with respect to an edge $E \subset \Delta_{f}$ is the polynomial
\[ f_{E}(x,y)  = \sum_{  (i,j) \in E \cap \supp_f } f_{ij}x^{i}y^{j}. \]
The following definition was introduced in \cite{kouchnirenko1976}:

\begin{ddef}
\label{def-function-nondeg}
{\rm A convenient power series $f\in\mathbb{C}\{x,y\}$ is said to be
\emph{Newton non-degenerate} if, for every edge
$E\subset\Delta_f$, the principal parts $(f_x)_E$ and
$(f_y)_E$ of $f_x$ and $f_y$ with respect to
$E$  have no common zeros in $(\mathbb{C}^*)^2 = (\mathbb{C} \setminus \{0\})^2$ .}
\end{ddef}

In other words, $f$ is Newton non-degenerate if and only if, for every edge $E \subset \Delta_{f}$, the polynomials  $(f_x)_{E}$ and $(f_y)_{E}$
in  $\C[x,y]$ have no common factor in $\C[x,y]$
other than a power of $x$ or a power of $y$.
Definition \ref{def-function-nondeg}
 can be extended to pairs of power series as follows:

\begin{ddef}
\label{def-pair-nondeg}
{\rm
Let $f,g \in \C \{x,y\}$ be convenient power series. We say that the map germ $(f,g): (\C^{2},0) \to (\C^{2},0)$ is
\emph{Newton non-degenerate} if one of the following conditions holds:
\begin{enumerate}
  \item the Newton polygons $\Delta_{f}$ and $\Delta_{g}$ have no pairs of parallel edges;
 
  \item   whenever  $E_{1} \subset \Delta_{f}$ and  $E_{2} \subset \Delta_{g}$ are  parallel edges, the principal partes
  $f_{E_{1}}(x,y)$ and $g_{E_{2}}(x,y)$ have no common zeros in $(\mathbb{C}^*)^2$.
\end{enumerate}
}\end{ddef}

For a pair  of convenient power series $f,g \in \C \{x,y\}$, we define their bracket by
\[ [\Ne_{f}, \Ne_{g}] = \sum_{E_{1} \subset \Delta_{f}} \sum_{E_{2} \subset \Delta_{g}}[E_1,E_2] .\]

Since, for $g=g_{1}g_{2}$, the Newton polygon $\Delta_{g}$ is obtained from $\Delta_{g_1}$ and $\Delta_{g_2}$ by merging the edges of equal co-slope, adding their lengths and heights, while edges of distinct co-slopes pass unchanged into $\Delta_g$, formula \eqref{eq-bracket-coslope} shows the bracket $[\,\cdot\,,\,\cdot\,]$ is additive in each argument. Directly from the definition,

\[ [\Ne_{f}, \Ne_{f}] = \sum_{E_{1} \subset \Delta_{f}}[E_1,E_1]
+  \sum_{E_{1},E_{2} \subset \Delta_{f}, E_{1} \neq E_{2}}[E_1,E_2] .\]

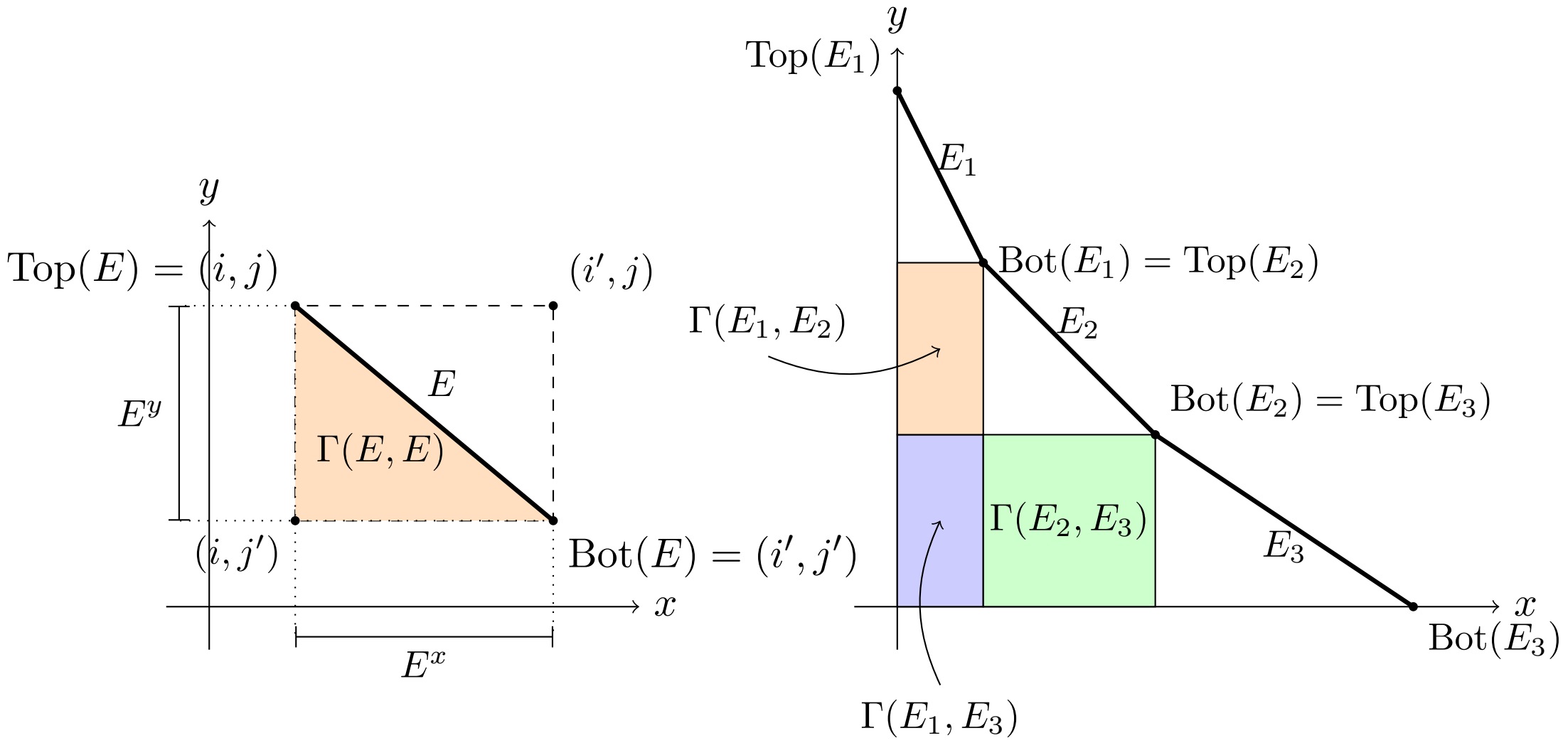
\begin{figure}[h]
\begin{tikzpicture}[scale=0.8]

\begin{scope}
  \draw[->] (-0.5,0) -- (5,0) node[right] {$x$};
  \draw[->] (0,-0.5) -- (0,4.5) node[above] {$y$};

  \draw[dashed] (1,1) -- (1,3.5) -- (4,3.5) -- (4,1) -- cycle;
  \fill[orange!25] (1,1) -- (1,3.5) -- (4,1) -- cycle;
  \draw[very thick] (1,3.5) -- (4,1);

  \fill (1,3.5) circle (1.5pt) node[above left] {$\Top(E)=(i,j)$};
  \fill (4,1) circle (1.5pt) node[below right] {$\Bot(E)=(i',j')$};
  \fill (1,1) circle (1.5pt) node[below left] {\small $(i,j')$};
  \fill (4,3.5) circle (1.5pt) node[above right] {\small $(i',j)$};

  \node at (2,1.8) {\small $\Gamma(E,E)$};
  \node at (2.7,2.6) {\small $E$};

  \draw[|-|] (1,-0.35) -- (4,-0.35) node[midway,below] {\small $E^x$};
  \draw[|-|] (-0.35,1) -- (-0.35,3.5) node[midway,left] {\small $E^y$};
  \draw[dotted] (1,3.5) -- (1,-0.35);
  \draw[dotted] (4,1) -- (4,-0.35);
  \draw[dotted] (1,3.5) -- (-0.35,3.5);
  \draw[dotted] (4,1) -- (-0.35,1);
\end{scope}

\begin{scope}[xshift=8cm]
  \draw[->] (-0.5,0) -- (7,0) node[right] {$x$};
  \draw[->] (0,-0.5) -- (0,6.5) node[above] {$y$};

  \draw[very thick] (0,6) -- (1,4);
  \draw[very thick] (1,4) -- (3,2);
  \draw[very thick] (3,2) -- (6,0);

  \fill[blue!20] (0,0) -- (1,0) -- (1,2) -- (0,2) -- cycle;
  \draw (0,0) -- (1,0) -- (1,2) -- (0,2) -- cycle;

  \fill[green!20] (1,0) -- (3,0) -- (3,2) -- (1,2) -- cycle;
  \draw (1,0) -- (3,0) -- (3,2) -- (1,2) -- cycle;

  \fill[orange!25] (0,2) -- (1,2) -- (1,4) -- (0,4) -- cycle;
  \draw (0,2) -- (1,2) -- (1,4) -- (0,4) -- cycle;

  \fill (0,6) circle (1.5pt) node[above left] {\small $\Top(E_1)$};
  \fill (1,4) circle (1.5pt) node[right] {\small $\Bot(E_1)=\Top(E_2)$};
  \fill (3,2) circle (1.5pt) node[above right] {\small $\Bot(E_2)=\Top(E_3)$};
  \fill (6,0) circle (1.5pt) node[below right] {\small $\Bot(E_3)$};

  \node at (0.7,5.2) {\small $E_1$};
  \node at (2.1,3.3) {\small $E_2$};
  \node at (4.5,0.7) {\small $E_3$};

\node (labE1E2) at (-1.5,3.3) {\small $\Gamma(E_1,E_2)$};
\draw[->, bend right=25] (labE1E2.south) to (0.5,3);
\node (labE1E3) at (0.5,-1.3) {\small $\Gamma(E_1,E_3)$};
\draw[->, bend left=25] (labE1E3.north) to (0.5,1);
  \node at (2,1) {\small $\Gamma(E_2,E_3)$};
\end{scope}
\end{tikzpicture}
\caption{Left: the rectangle $E^x\times E^y$ split by the diagonal $E$ into two
triangles; the shaded lower-left triangle is $\Gamma(E,E)$. Right: three
consecutive edges $E_1,E_2,E_3\subset\Delta_f$ with $\kappa(E_1)<\kappa(E_2)<\kappa(E_3)$,
and the rectangles $\Gamma(E_1,E_2)$, $\Gamma(E_1,E_3)$, and $\Gamma(E_2,E_3)$.}
\label{pict:Gamma}
\end{figure}

For an edge $E\subset\Delta_f$, let $E^x$ and $E^y$ denote its projections
onto the coordinate axes. Then $E^x\times E^y$ is a rectangle, divided by
the diagonal $E$ into two right triangles. Writing $\Top(E)=(i,j)$ and
$\Bot(E)=(i',j')$, the vertices of $E^x\times E^y$ are $(i,j')$, $(i',j')$,
$(i,j)$, $(i',j)$, and the diagonal $E$ joins $(i,j)$ to $(i',j')$; the two
triangles are then the one containing $(i,j')$ and the one containing
$(i',j)$. We denote by $\Gamma(E,E)$ the triangle containing the vertex
$(i,j')$, which we call the \emph{lower-left triangle} associated with $E$
(see the left side of Figure \ref{pict:Gamma}).
Observe that
$[E,E] = \ell(E) h(E) = 2 A(\Gamma(E,E))$, where $A$ denotes the area.
Now let $E_1,E_2\subset\Delta_f$ be two distinct edges. We define
$\Gamma(E_1,E_2)$ to be the rectangle $E^x_{1}\times E^y_{2}$ if
$\kappa(E_1)<\kappa(E_2)$, and $E^x_{2}\times E^y_{1}$ otherwise.
In either case, \eqref{eq-bracket-coslope} yields $[E_{1},E_{2}] = A(\Gamma(E_{1},E_{2}))$ (see the right side of Figure \ref{pict:Gamma}).
The rectangles and triangles defined above will be referred to as
\emph{basic polygons}.

Given a convenient power series  $f \in \C \{x,y\}$, denote by $\A(\Ne_{f})$ the area
of the polygonal region $\RR_{f}$  delimited by $\Delta_{f}$ and the coordinate axes.  Note that, although $\Ne_f$ itself is unbounded, $\A(\Ne_f)$ denotes the
(finite) area of the bounded region $\RR_f$, consistently with the notation
$h(\Ne_f)$ and $\ell(\Ne_f)$ introduced above.
Observe that
\begin{equation}
\label{eq-Rf}
 \RR_{f} = \left(\bigcup_{E \subset \Delta_{f}} \Gamma(E,E) \right) \cup \left( \bigcup_{E_{1},E_{2} \subset \Delta_{f}, \kappa(E_{1}) < \kappa(E_{2})}  \Gamma(E_{1},E_{2}) \right),
\end{equation}
with the basic polygons  appearing in this decomposition overlapping, at most, along  common edges.
We have:

\begin{prop}
\label{prop-bracket}
 Let $f,g \in \C \{x,y\}$ be convenient power series. Then:
\begin{enumerate}[label=(\alph*)]
  \item $[\Ne_{f}, \Ne_{f}] = 2 \A(\Ne_{f})$;
  \smallskip \smallskip
  \item  $[\Ne_{f}, \Ne_{g}] =  \A(\Ne_{f} + \Ne_{g}) - \A(\Ne_{f}) - \A(\Ne_{g})$, where $\Ne_f+\Ne_g$ denotes the Minkowski sum.
\end{enumerate}
\end{prop}
\begin{proof}
(a)
We have, by \eqref{eq-Rf},
\begin{eqnarray*}
 \A(\Ne_{f}) = A(\RR_{f}) &=& \sum_{E \subset \Delta_{f}} A(\Gamma(E,E))
 + \sum_{E_{1},E_{2} \subset \Delta_{f}, \kappa(E_{1}) < \kappa(E_{2})} A(\Gamma(E_{1},E_{2})) \\
    &=&  \sum_{E \subset \Delta_{f}} \frac{1}{2}[E,E] +  \sum_{E_{1},E_{2} \subset \Delta_{f}, \kappa(E_{1}) < \kappa(E_{2})} [E_{1},E_{2}]\\
  &=& \frac{1}{2} \left( \sum_{E \subset \Delta_{f}} [E,E] +  \sum_{E_{1},E_{2} \subset \Delta_{f}, \kappa(E_{1}) \neq \kappa(E_{2})} [E_{1},E_{2}] \right) = \frac{1}{2} [\Ne_{f}, \Ne_{f}].
\end{eqnarray*}
(b) Let $\Ne = \Ne_{f} + \Ne_{g}$ and denote by $\Delta$ the compact part of its polygonal border and by
$\RR$ the polygonal region  delimited by $\Delta$ and the coordinate axes.
We have that $\Delta$  is the polygonal line of height $h(\Delta) =  h(\Delta_{f}) + h(\Delta_{g})$ and length
$\ell(\Delta) =  \ell(\Delta_{f}) + \ell(\Delta_{g})$ built through the juxtaposition of copies $\tilde{E}_{1}$ of edges $E_{1}$ of
$\Delta_{f}$  and  copies $\tilde{E}_{2}$ of edges $E_{2}$ of
$\Delta_{g}$, ordered by growing co-slope. Note that edges $E_{1} \subset \Delta_{f}$
and $E_{2} \subset \Delta_{g}$ of the same co-slope $\kappa$ produce a single edge of $\Delta$ of co-slope $\kappa$,
which can be considered as the juxtaposition of $\tilde{E}_{1}$ and $\tilde{E}_{2}$.
We can thus decompose $\RR$ into the union of  three types of basic polygons:
\begin{enumerate}
  \item Rectangles $\Gamma(\tilde{E}_{1}, \tilde{E}_{2})$, where $\tilde{E}_{1}, \tilde{E}_{2} \subset \Delta$ are edges corresponding
to $E_{1} \subset \Delta_{f}$ and $E_{2} \subset \Delta_{g}$; the sum of their areas is $[\Ne_{f}, \Ne_{g}]$.
  \item Basic polygons of the form $\Gamma(\tilde{E}_{1}, \tilde{E}_{1}')$, where $\tilde{E}_{1}, \tilde{E}_{1}' \subset \Delta$
  correspond to $E_{1},E_{1}' \subset \Delta_{f}$ (possibly equal); the sum of their areas is $A(\RR_{f}) = \A(\Ne_{f})$.
   \item Basic polygons of the form $\Gamma(\tilde{E}_{2}, \tilde{E}_{2}')$, where $\tilde{E}_{2}, \tilde{E}_{2}' \subset \Delta$
  correspond to $E_{2},E_{2}' \subset \Delta_{g}$ (possibly equal); the sum of their areas is $A(\RR_{g}) =\A(\Ne_{g})$.
\end{enumerate}
The result then follows by observing that the sum of the areas of all basic polygons listed above is $A(\RR) = \A(\Ne_{f} + \Ne_{g})$.
\end{proof}

The following theorem characterizes Newton non-degeneracy for a pair of power series.

\begin{teo}
\label{teo-nondeg1} Let $f,g \in \C\{x,y\}$ be convenient power series. Then
\[ (f,g)_{0} \geq [\Ne_f,\Ne_g], \]
with equality holding if and only if the pair $(f,g)$ is Newton non-degenerate.
\end{teo}
\begin{proof}
By Newton-Puiseux theorem (see, for instance, \cite{casas2000}), given an edge $E \subset \Delta_{f}$, there are
 $h(E)$ $y$-roots of order $\kappa = \kappa(E) =  \frac{\ell(E) }{h(E)}$, i.e. fractional power series $y_{E}(x) \in \C\{x\}^{\circledast} := \bigcup_{n \in \N} \C\{x^{\frac{1}{n}}\}$, each of them of the form $y_{E}(x) =  a x^{\kappa} + \cdots$, where $a \in \C^{*}$, such that  $f(x,y_{E}(x)) = 0$. In particular, the initial term
 $\In(y_{E}(x)) = a x^{\kappa}$ satisfies $f_{E}(x,a x^{\kappa}) = 0$.
Similarly, switching the roles of $x$ and $y$, there are $\ell(E)$ $x$-roots of the form
$x_{E}(y) =  b y^{\rho} + \cdots$, where $b \in \C^{*}$, of order $\rho = \frac{1}{\kappa} =  \frac{h(E) }{\ell(E) }$,
the initial term $\In(x_{E}(y)) = b y^{\rho}$ satisfying $f_{E}(b y^{\rho},y) = 0$.
The proof of the theorem relies on the following lemma:

\begin{lem}
\label{lem-nodeg-pair}
Let $f,g \in \C\{x,y\}$ be convenient power series. The following assertions are equivalent:
\begin{enumerate}
  \item the map germ $(f,g): (\C^{2},0) \to (\C^{2},0)$ is Newton non-degenerate;
  \item for every pair of fractional power series $y_{f}(x), y_{g}(x)   \in \C\{x\}^{\circledast}$, which are, respectively,
  $y$-solutions of $f(x,y) = 0$ and  $g(x,y) = 0$,
   we have
\[ \ord_{x} (y_{f}(x) - y_{g}(x)) = \min \{\ord_{x}  y_{f}(x),\ord_{x}  y_{g}(x)\} .\]
\end{enumerate}
\end{lem}
\begin{proof}[Proof {\rm (of the lemma)}]
We first observe that the minimality condition of item (2) is satisfied if and only if
$\In(y_{f}(x)) \neq \In(y_{g}(x))$ for every pair of $y$-roots $y_{f}(x)$ and $y_{g}(x)$ of $f(x,y)$ and  $g(x,y)$.
Now, recalling the Newton-Puiseux method, if $E_{1} \subset \Delta_{f}$ is an edge, then
the principal part    $f_{E_{1}}$  is a product of a monomial in the variables $x,y$  and a polynomial
$f_{E_{1}}^{*}$ given by the product of $h(E_{1})$ terms of the form $ y - a x^{\alpha}$, where $\In(y_{f}(x)) = a x^{\alpha}$
and $y_{f}(x)$ runs over
the $y$-roots of $f(x,y)$ of order $\alpha$. A similar statement is true for each
edge $E_{2} \subset \Delta_{g}$. If $\Delta_f$ and $\Delta_g$ have no pair of parallel edges, every $y$-root of $f$
has an order different from that of every $y$-root of $g$, so the minimality
condition holds trivially for every pair. It remains to consider the case where
$\Delta_f$ and $\Delta_g$ have parallel edges.
Now, if we have parallel edges
$E_{1} \subset \Delta_{f}$ and  $E_{2} \subset \Delta_{g}$, say of co-slope $\kappa = \alpha = \beta$,
then $f_{E_{1}}^{*}$ and $g_{E_{2}}^{*}$ have
no common factors if and only if all terms   $y - a x^{\kappa}$ of $f_{E_{1}}^{*}$ and $y - b x^{\kappa}$ of $g_{E_{2}}^{*}$
are different or,
equivalently, the initial terms  $a x^{\kappa}$ of the $y$-roots of $f$ and $b x^{\kappa}$ of the $y$-roots of $g$ never coincide.
This gives the lemma's proof.
\end{proof}

The proof of the theorem also requires the following result:
\medskip

\par \noindent
{\bf Zeuthen's rule} \cite[Proposition 2.3]{ploski1990}{\bf .}
\emph{Let $f,g \in \C\{x,y\}$ be convenient power series. Then
\[ (f,g)_{0} = \sum_{i,j}  \ord_{x}(y_{f,i}(x) - y_{g,j}(x)) ,\]
 where $y_{f,i}(x)$ and $y_{g,j}(x)$ are,  respectively, the $y$-roots of $f(x,y)$ and of $g(x,y)$.
}

\begin{proof}[Proof {\rm (of Zeuthen's rule)}]  Applying  Weierstrass preparation theorem,   write $f = u f^{*}$ and $g = v g^{*}$,
where $u,v \in  \C\{x,y\}$ are units and $f^{*}, g^{*} \in \C\{x\}[y]$ are such that
\[ f^{*}(x,y) =  \prod_{i} (y - y_{f,i}(x)) \qquad \text{and} \qquad   g^{*}(x,y) =  \prod_{j} (y - y_{g,j}(x)) .\]
Since $u,v$ are units, $(f,g)_{0} = (f^{*},g^{*})_{0}$. By \cite[Th. 4.17]{hefez2003}, the
intersection number of two Weierstrass polynomials equals the $x$-order of their
resultant in $y$, so
\[ (f^{*},g^{*})_{0} = \ord_{x} {\rm Res}_{y}(f^{*},g^{*}) .\]
Since $f^{*}$ and $g^{*}$ are monic in $y$, their resultant equals the product of
the differences of their roots, so
\[ {\rm Res}_{y}(f^{*},g^{*}) = \prod_{i,j} (y_{f,i}(x) - y_{g,j}(x)) .\]
Finally, since $\ord_{x}$ is additive over products,
\[ \ord_{x} \prod_{i,j} (y_{f,i}(x) - y_{g,j}(x)) = \sum_{i,j} \ord_{x}(y_{f,i}(x) - y_{g,j}(x)) .\]
The combination of the four displayed equalities gives the rule.
\end{proof}

Proceeding with the proof  of Theorem \ref{teo-nondeg1}, we have, by Zeuthen's rule,
\begin{eqnarray*}
  (f,g)_{0} &=& \sum_{i,j}  \ord_{x}(y_{f,i}(x) - y_{g,j}(x))\\
    &\geq& \sum_{i,j}  \min \{\ord_{x}(y_{f,i}(x)), \ord_{x}(y_{g,j}(x)) \},
\end{eqnarray*}
with equality if and only if the pair $(f,g)$ is Newton non-degenerate, by Lemma \ref{lem-nodeg-pair}.
However,
\begin{eqnarray*}
\sum_{i,j}  \min \{\ord_{x}(y_{f,i}(x)), \ord_{x}(y_{g,j}(x)) \}
   &=& \sum_{E_{1} \subset \Delta_{f}} \sum_{E_{2} \subset \Delta_{g}}  h(E_{1})   h(E_{2})
   \min \left\{ \frac{\ell(E_{1})}{h(E_{1})},  \frac{\ell(E_{2})}{h(E_{2})}     \right\}  \\
   &=& \sum_{E_{1} \subset \Delta_{f}} \sum_{E_{2} \subset \Delta_{g}} [E_{1},E_{2}] = [\Ne_{f},\Ne_{g}]
\end{eqnarray*}
\end{proof}

\begin{example}
\label{ex-pair-nondeg}
{\rm
Let $f=y^{2}-x^{3}$ and $g=y-x^{2}$. The Newton polygons
$\Delta_{f}$ and $\Delta_{g}$ each consist of a single edge, of
co-slope $\kappa(\Delta_{f})=3/2$ and $\kappa(\Delta_{g})=2$
respectively. Since these co-slopes differ, $\Delta_{f}$ and
$\Delta_{g}$ have no pair of parallel edges: the pair $(f,g)$ is
Newton non-degenerate by Definition \ref{def-pair-nondeg}. The
bracket is
\[ [\Ne_{f},\Ne_{g}] = \min\{\ell(E_{f})h(E_{g}),\,h(E_{f})\ell(E_{g})\}
= \min\{3\cdot1,\,2\cdot2\} = 3, \]
matching the intersection number $(f,g)_0=3$, computed via the
resultant $\mathrm{Res}_y(f,g)=x^{3}(x-1)$. Theorem
\ref{teo-nondeg1} is verified: $(f,g)_0=[\Ne_{f},\Ne_{g}]$.
}\end{example}

\subsection{Newton non-degenerate foliations}
\label{subsec-Newton-and-foliations}

Let $\F$ be a   foliation at $(\C^{2},0)$ defined by the germ of holomorphic $1$-form
$\omega = A(x,y) \dd x + B(x,y) \dd y$.
Let us write
 \[ \omega = \sum_{i,j} \omega_{ij} ,\]
where
\[ \omega_{ij} = a_{ij} x^{i-1}y^{j} \dd x +  b_{ij} x^{i}y^{j-1} \dd y\]
is called \emph{differential  binomial} of $\omega$ of \emph{bidegree} $(i,j)$.
We remark that, under this decomposition, $a_{0j} = 0$ and $b_{i0}=0$ for all $i,j \in \N$ and,
in particular, $\omega_{00} = 0$.
The \emph{support} of $\omega$ is
\[ \supp_\omega = \supp_{xA} \cup \supp_{yB} = \{(i,j) \in \N \times \N; \omega_{ij} \neq 0 \}.\]
Its \emph{Newton diagram} $\Ne_\omega$ is $\Ne(\supp_\omega)$ and its \emph{Newton polygon} $\Delta_\omega$
is the union of compact edges of  $\Ne_\omega$. Observe that, if $\varphi \in  \C\{x,y\}$ is a unit and
$\omega' = \varphi \omega$,  then $\Delta_{\omega'} = \Delta_{\omega}$ and, hence, once a system of analytic
coordinates $(x,y)$ is fixed,
 the Newton polygon is in fact an object associated with  the foliation $\F$ itself.
  However, we will avoid the notation $\Ne_{\F}$ or $\Delta_{\F}$ and use
 $\Ne_\omega$ or $\Delta_{\omega}$ instead, with  reference only to
  the $1$-form $\omega$ that defines $\F$. This is more suitable for the calculations to be performed  and also
stresses the dependency on the choice of the coordinates.
 Note also that
  $\Delta_\omega \cap \supp_\omega = \Delta_{\omega'} \cap \supp_{\omega'}$ and, for every $(i,j)$ in this set,   $\omega_{ij}' = \alpha \omega_{ij}$,
 where $\alpha = \varphi(0,0) \in \C^{*}$.
Let us call the points in $\Delta_\omega \cap \supp_\omega$ by
 \emph{particles}. Note that, in particular, all vertices of  $\Delta_\omega$ are particles.
   We say that a particle is a \emph{contribution} of $A$ or of $B$ if it belongs, respectively,
 to $\supp_{xA}$ or to $\supp_{yB}$. If $(i,j)$ is a particle, then we define its \emph{multiplier} as
 $\lambda_{ij} = - a_{ij}/ b_{ij}$ when $(i,j)$ is a contribution of $B$, setting  $\lambda_{ij} = \infty$, otherwise.
 For a given edge $E \subset \Delta_\omega$ we denote the set of its particles by $\PP(E)$.

Given an edge $E$ of  $\Delta_{\omega}$, its \emph{degeneracy equations}   are
\begin{equation}
\label{eq-degeneracy}
  (xA)_{E} = \sum_{(i,j) \in \PP(E)} a_{ij} x^{i}y^{j}  \qquad \text{and}  \qquad
(yB)_{E} = \sum_{(i,j) \in \PP(E)} b_{ij} x^{i}y^{j} .
\end{equation}
We say that $\omega$ is \emph{Newton $E$-degenerate} if  $(xA)_{E}=0$ and  $(yB)_{E}=0$ have common solutions outside the coordinate
axes or, equivalently, if  $(xA)_{E}$ and  $(yB)_{E}$ have a non-trivial (i.e. neither $x$ nor $y$) common factor.
Otherwise,   $\omega$ is said to be \emph{Newton $E$-non-degenerate}.
We say that $\omega$ is \emph{Newton non-degenerate} if it is Newton $E$-non-degenerate
with respect to all its edges.
Observe that the replacement of $\omega$ by another $1$-form defining $\F$ entails  the  multiplication of the degeneracy equations by a non-zero scalar.
Thus, for a fixed system of analytic coordinates, the notion of $E$-degeneracy is a property of the foliation $\F$ itself, independent of the defining $1$-form.

The coordinates $(x,y)$ are \emph{convenient} for $\omega$ if $\Delta_{\omega}$ intersects both coordinate axes.
This means that $\supp_\omega$ contains points of the form $(i_0,0)$ and $(0,j_{0})$ or, in other words,  that
$A(x,y)$ has a non-zero term $a_{i_0 0} x^{i_0-1}$ and $B(x,y)$ has a non-zero term $b_{0 j_0} y^{j_0-1}$.
In turn, this is equivalent to the fact that the coordinate axes are not invariant by $\F$.
  With the concepts above, we have the following definition:
\begin{ddef}
\label{def-Newton-nondeg}
{\rm
The foliation $\F$, with defining $1$-form $\omega$, is \emph{Newton non-degenerate} if there are convenient coordinates in which $\omega$ is  Newton non-degenerate.
}
\end{ddef}

If our coordinates are non-convenient then either the topmost point of $\Delta_{\omega}$ is of the form $(1,j_{0})$
or its bottommost point is of the kind $(i_0,1)$.
Indeed, if $\Delta_{\omega}$
does not touch the vertical axis, then $x=0$ is $\F$-invariant    and, in this case, $\F$ is defined by a $1$-form  $\omega = A(x,y) \dd x + B(x,y) \dd y$, where $x$ must divide  $B$ but cannot divide $A$. Thus $A$ has a non-zero term of type $a_{1 j_{0}} y^{j_{0}}$ which, if $j_{0}$ is minimal, yields the vertex
$(1,j_{0})$ of $\Delta_{\omega}$. A similar argument applies when $\Delta_{\omega}$ does not intersect the horizontal axis.
Note that an analytic change of coordinates that is  highly tangent to the identity ---
that is, one given by a map  of the form $\Psi(x,y) = (x + \cdots, y + \cdots)$, the dots standing for
terms of sufficiently high order
 --- preserves the edges of
$\Delta_{\omega}$ and their  degeneracy equations.
Let us explain this fact. Let  $(\alpha,0)$ and $(0,\beta)$ be points lying in the coordinate axes
belonging to the supporting line of, respectively,  the last and  the first edge of $\Delta_{\omega}$, and
set $\gamma = \max\{i+j; (i,j) \in \Delta_{\omega}\}$.
Define $M = \max\{\alpha, \beta, \gamma\}$ and choose $\Psi$ such that its terms represented by  dots
have order at least  $M$. Since $A(0,0) = B(0,0) = 0$  (so every particle has $i+j\geq2$),
these terms  generate  terms in $\Psi^{*} \omega$ of order at least $M$, that in turn  produce points
$(i,j) \in \supp_{\Psi^{*} \omega}$ with $i+j > M$.
Hence,   the terms corresponding to points in $\Delta_{\omega}$ of both $\omega$ and $\Psi^{*} \omega$ are the same.
It is thus possible to obtain new convenient coordinates preserving the data of $\Delta_{\omega}$ and introducing the one or two missing edges that should touch the coordinate axes.
We also introduce the following definition:  the coordinates $(x,y)$ are \emph{semi-convenient} if  $\Delta_{\omega}$ intersects only the horizontal axis or, equivalently, if the $y$-axis   is $\F$-invariant  but not the $x$-axis.
Along this article, we  always assume that our initial coordinates are convenient, although,
as it will be apparent  below,
  semi-convenient coordinates appear in the process of reduction of singularities.

We can define the \emph{element} of $\Delta_{\omega}$ of co-slope $\kappa \in \Q_{>0}$.
This is an edge or a vertex of $\Delta_{\omega}$ chosen in the following way.   Take, among all lines of equation $i + \kappa j = \delta$ that
intersect the Newton diagram $\Ne_{\omega}$, the one with minimal $\delta$. It intersects $\Ne_{\omega}$ over    $\Delta_{\omega}$, either on
an   edge of co-slope $\kappa$  or a vertex. This edge or vertex is   the element of $\Delta_{\omega}$ of co-slope $\kappa$,
represented, in both cases, by $E_{\kappa}$.
The context will make it clear whether it refers to an edge or to a vertex.
The corresponding line is the $\emph{supporting line}$ of co-slope $\kappa$.
Elements that are not edges will be needed in the sequel, in order to handle
the resonant case and dicritical singularities. Observe that, if we take $\kappa = 1$, then each particle $(i,j) \in E_{1}$ has   minimum order $i + j -1$. Hence, if $i + j = \delta$ is the equation of the supporting line of co-slope $\kappa = 1$,
we have that $\ma_{0}(\F) = \delta -1$.

 If $E_{\kappa}$ is the element of  $\Delta_{\omega}$ of co-slope $\kappa$, then the \emph{initial polynomial} of $\omega$
 with respect to $E_{\kappa}$ is
\[ \Phi_{E_{\kappa}}(z) =
(xA)_{E_{\kappa}}(1,z) +
 \kappa  (y B)_{E_{\kappa}}(1,z) =
 \sum_{(i,j) \in P(E_{\kappa})} ( a_{ij} + \kappa  b_{ij}) z^{j} . \]
If $\Phi_{E_{\kappa}}(z) \equiv 0$, we say that
the element  $E_{\kappa}$ is \emph{dicritical} or that $\Delta_{\omega}$ is \emph{$\kappa$-dicritical}.
This is evidently equivalent to asking that
$(xA)_{E_{\kappa}}(x,y) + \kappa  (y B)_{E_{\kappa}}(x,y) \equiv 0$ and
occurs if and only if $\lambda_{ij} = \kappa$   for every particle
$(i,j) \in E_{\kappa}$,   where $\lambda_{ij}$ is the multiplier of $(i,j)$.
As we will see below, in the case $\kappa  \in \Z_{>0}$, this is equivalent to the fact that,
in a sequence of blow-ups with centers following the $x$-direction,  the $\kappa$-th blow-up is dicritical.
If $E_{\kappa}$ is \emph{non-dicritical}, i.e.,  if $\Phi_{E_{\kappa}}(z) \not\equiv 0$,
 then we have
a factorization  $\Phi_{E_{\kappa}}(z) = z^{\ell}\Phi_{E_{\kappa}}^{*}(z)$,
where $\ell \geq e(E_{\kappa})$
and $\Phi_{E_{\kappa}}^{*}$ is a polynomial of degree at most $h(E_{\kappa})$ satisfying $\Phi_{E_{\kappa}}^{*}(0) \neq 0$.
We have that $\Phi_{E_{\kappa}}^{*}$ is  non-constant if and only if there are at least two particles on
$E_{\kappa}$ with multipliers different from $\kappa$.

Writing
\[ \alpha_{\kappa}(z) = (xA)_{E_{\kappa}}(1,z) \qquad \text{and}  \qquad
  z \beta_{\kappa}(z)= (y B)_{E_{\kappa}}(1,z)   ,\]
we have (recall that
$b_{i,0}=0$, so every particle contributing to $B$ has $j\geq1$;
hence every term of $(yB)_{E_\kappa}(1,z)$ has $z$-exponent
$j\geq1$, and $(yB)_{E_\kappa}(1,z)$ has no constant term)
\begin{equation}
\label{eq-principal-polynomial}
 \Phi_{E_{\kappa}}(z) =  \alpha_{\kappa}(z) +  \kappa z \beta_{\kappa}(z) =
  \sum_{(i,j) \in P(E_{\kappa})} a_{ij} z^{j} + \kappa z \left( \sum_{(i,j) \in P(E_{\kappa})} b_{ij} z^{j-1} \right) .
\end{equation}

For $\kappa \in \Q_{>0}$ and $c \in \C^{*}$,  consider a  fractional
power series     $y = c x^{\kappa} + \cdots$,  where the dots represent terms of higher rational order,
and set $\Psi(x) = (x, c x^{\kappa} + \cdots)$. Suppose that the supporting line of co-slope $\kappa$ has equation
$i + \kappa j = \delta$. Given a particle $(i,j)$ of  $E_{\kappa}$, write  $\omega_{ij} = a_{ij} x^{i-1}y^{j} \dd x + b_{ij} x^{i}y^{j-1} \dd y$,
the differential binomial of bidegree $(i,j)$. Then
\begin{eqnarray*}
  \Psi^{*} \omega_{ij} &=& ( a_{ij}  x^{i-1}( c x^{\kappa})^{j} +  b_{ij}  x^{i}( c x^{\kappa})^{j-1} \kappa x^{\kappa-1} + \cdots ) \dd x \medskip \\
    &=& ( ( a_{ij}  +  \kappa b_{ij} ) c^{j}     x^{\delta -1}  + \cdots ) \dd x,
\end{eqnarray*}
where the dots refer to terms of order at least $\delta$ in the variable $x$.
The term of least order of $\Psi^{*} \omega$ is generated in this way  by considering all particles over $E_{\kappa}$ and
has the form  $\Phi_{E_{\kappa}}(c)  x^{\delta -1}  \dd x$.
Hence, if $S$ is a separatrix of  $\F$, different from the $y$ axis, and $y = c x^{\kappa} + \cdots$ is one of its
Puiseux series, we have $\Psi^{*}\omega \equiv 0$ and, therefore, $\Phi_{E_{\kappa}}(c) = 0$. This is the cornerstone of J. Cano's
algorithmic proof of the Separatrix theorem \cite[Th. 1]{jcano1993}. We will say that $E_{\kappa}$ is the element \emph{associated with} $S$
or that $S$ is an $E_{\kappa}$-separatrix.
If $E_{\kappa}$ is dicritical then, by definition, $\Phi_{E_{\kappa}}(c) = 0$ is satisfied for every $c \in \C^{*}$.
This must happen, for instance, when the element
 $E_{\kappa}$ associated with $S$ is a single vertex, since, in this case, $\Phi_{E_{\kappa}}$ would be a monomial if non-vanishing.
Therefore, when  $E_{\kappa}$ is non-dicritical,   it is necessarily an edge of $\Delta_{\omega}$. In this case, in the  factorization
$\Phi_{E_{\kappa}}(z) = z^{\ell}\Phi_{E_{\kappa}}^{*}(z)$, the polynomial $\Phi_{E_{\kappa}}^{*}(z)$ is non-constant.

The following simple fact will be valuable to understand a Newton non-degenerate foliation:
\begin{prop}
\label{prop-dicritical-edge}
If $E \subset \Delta_{\omega}$ is a dicritical edge of the Newton polygon of the foliation $\F$,
then $\omega$ is Newton $E$-degenerate.
\end{prop}
\begin{proof} Denote by $\omega = A(x,y) \dd x  + B (x,y)\dd y$ the $1$-form defining $\F$ and by $\kappa \in \Q_{>0}$ the co-slope of $E$.
The edge $E$ being dicritical is equivalent to the equation
$(xA)_{E} = - \kappa  (y B)_{E}$, which in turn is equivalent to
$a_{ij}= -\kappa b_{ij}$ for every  $(i,j)\in E$.  So,
every particle of $E$ is a contributions of both $A$ and $B$, and this fact
holds, in particular, for $\Top(E)$ and $\Bot(E)$.
As a consequence,   $(xA)_E$ and $(yB)_E$ are non-monomial.
Now, by the above equation, this pair of polynomials have the
same irreducible factors.
Being non-monomial,
they have at least a non-trivial common factor.
This implies that $\omega$ is Newton $E$-degenerate, as wished.
\end{proof}

For a germ of holomorphic foliation $\F$, fixing a system of convenient local analytic coordinates $(x,y)$, we define its \emph{Newton number}
in the same way we did above for power series:
\begin{equation}
\label{def-newton-number}
 \nu_{0}(\F) = 2 \cl{A}(\Ne_{\omega}) - \ell(\Ne_{\omega}) - h(\Ne_{\omega}) + 1 ,
\end{equation}
where $\omega = A(x,y) \dd x + B(x,y) \dd y$ is a germ of analytic $1$-form defining $\F$.
Observe that $\nu_0(\F)$ is well defined, since $\Ne_\omega$ (and hence $\Delta_\omega$) does not
depend on the choice of the $1$-form $\omega$ defining $\F$, as remarked in \S\ref{subsec-Newton-and-foliations}.
The \emph{Milnor number} of $\F$ at $0 \in \C^{2}$ (see \cite{camacho1982}) is defined as
\begin{equation}
\label{def-milnor-number}
 \mu_{0}(\F) = (A, B)_{0}  = \dim_{\C} \frac{\C \{x,y \}}{(A,B)},
\end{equation}
where $(A,B)$ denotes the ideal of $\C \{x,y \}$ generated by $A$ and $B$.
These numbers can be compared, providing a version of Kouchnirenko's theorem for foliations:

\begin{teoRM}
\label{teo-fol-nondeg2}
Let $\F$ be a germ of holomorphic foliation at $(\C^2,0)$ and let $(x,y)$ be convenient analytic coordinates.  Then
\[  \nu_{0}(\F) \leq \mu_{0}(\F) \]
with equality holding if and only if $\F$ is Newton non-degenerate.
\end{teoRM}

\begin{proof}  Suppose, as above, that  $\F$ is defined by the $1$-form $\omega = A \dd x + B \dd y$, with
$A,B \in \C\{x,y\}$.
Define  the functions
\begin{equation}
\label{eq-geradores}
\begin{cases}
  f_{1} \  = \ \alpha_{1} x A + \beta_{1} y B   \\
  f_{2} \ = \  \alpha_{2} x A + \beta_{2} y B ,
\end{cases}
\end{equation}
where  the pairs $(\alpha_{1},\beta_{1}), (\alpha_{2},\beta_{2}) \in \C^{2}$ are generic,
meaning that $ \Ne_{\omega} = \Ne_{f_{1}} = \Ne_{f_{2}}$ and also that $\alpha_{1} \beta_{2} -  \beta_{1} \alpha_{2} \neq 0$.
Such pairs exist because, for $k=1,2$, $\supp_{f_k}\subset\supp_\omega$ always, and a vertex of
$\Delta_\omega$ remains a vertex  of $\Delta_{f_k}$ as long as the corresponding coefficient
$\alpha_k a_{ij}+\beta_k b_{ij}$ does not vanish. This excludes only
finitely many ratios $\alpha_k/\beta_k$ (one per vertex of  $\Ne_{\omega}$), together with the single ratio making
$\alpha_1\beta_2-\beta_1\alpha_2=0$, so a generic choice avoids all of them.
Since $\Delta_{\omega} =  \Delta_{f_{1}} = \Delta_{f_{2}}$ it follows directly from
Definition 
\ref{def-pair-nondeg}, applied to the pairs $(xA,yB)$ and $(f_{1},f_{2})$,
 that $\F$ is Newton non-degenerate if and only if  the pair $(f_{1},f_{2})$ is Newton non-degenerate.
 We have
\begin{eqnarray*}
  (x A,  y B)_{0} &=& (x,  y)_{0} +  (A,  y)_{0} +  (x, B)_{0} +  (A,B)_{0} \\
   &=& 1 +  \ord_{0}(A(x,0)) +  \ord_{0}(B(0,y))    + \mu_{0}(\F)  \\
   &=& 1 + (\ell(\Ne_{\omega})-1) + (h(\Ne_{\omega})-1) + \mu_{0}(\F) \\
    &=& -1 +  \ell(\Ne_{\omega})  +  h(\Ne_{\omega})  + \mu_{0}(\F).
\end{eqnarray*}
Hence
\[ \mu_{0}(\F) =   (x A,  y B)_{0} - \ell(\Ne_{\omega}) - h(\Ne_{\omega}) + 1 . \]
On the other hand, from equations \eqref{eq-geradores}, the ideals $(f_{1},f_{2})$ and $(x A,y B)$ of $\C\{x,y\}$ are equal and, hence
\[ (x A, y B )_{0} =  (f_{1},f_{2})_{0} .\]
Using this fact, Theorem \ref{teo-nondeg1} and Proposition \ref{prop-bracket} (take into account
that $\Ne_{\omega}$ is $\Ne_{f}$ for some $f \in \C\{x,y\}$), we have
\begin{eqnarray*}
 \mu_{0}(\F) &=&  (f_{1},f_{2})_{0}  - \ell(\Ne_{\omega}) - h(\Ne_{\omega}) + 1   \\
  & \geq & [\Ne_{f_{1}},\Ne_{f_{2}}]    - \ell(\Ne_{\omega}) - h(\Ne_{\omega}) + 1 \\
   & = & [\Ne_{\omega},\Ne_{\omega}]    - \ell(\Ne_{\omega}) - h(\Ne_{\omega}) + 1 \\
   & = & 2 \A(\Ne_{\omega})    - \ell(\Ne_{\omega}) - h(\Ne_{\omega}) + 1  = \nu_{0}(\F).
\end{eqnarray*}
The equality in the second line is equivalent, by Theorem \ref{teo-nondeg1}, to the pair
$(f_1,f_2)$ being Newton non-degenerate, which, as observed above, is in turn equivalent
to $\F$ being Newton non-degenerate, which concludes the proof.
\end{proof}

\begin{example}
\label{ex-fol-nondeg}
{\rm
Let $\F$ be defined by
\[\omega = (-12x^2y+2x^2-12xy^2+14xy-3y^3+11y^2)\,\dd x +
(12x^3+12x^2y-5x^2+3xy^2+xy+4y^2)\,\dd y .\]
 Since $A(x,0)=2x^2$
and $B(0,y)=4y^2$ are non-zero, the coordinates $(x,y)$ are convenient. The
Newton polygon $\Delta_\omega$ has vertices $(0,3)$ and $(3,0)$,
giving $\nu_0(\F)=2 \frac{9}{2} -3-3+1=4$. Computing
$\mu_0(\F)=(A,B)_0$ via the resultant, whose lowest-order term is
$x^{4}$, gives $\mu_0(\F)=4$. Theorem \ref{teo-fol-nondeg2} is
verified: $\nu_0(\F)=\mu_0(\F)=4$.
}\end{example}

For a convenient power series $f  \in \C\{x,y\}$, the \emph{Newton number} is defined as
\[ \nu_{0}(f) = 2 \cl{A}(\Ne_{f}) - \ell(\Ne_{f}) - h(\Ne_{f}) + 1 ,\]
whereas the \emph{Milnor number} is
\[ \mu_{0}(f) = (f_{x}, f_{y})_{0}  = \dim_{\C} \frac{\C \{x,y \}}{(f_{x},f_{y})}.\]
As a corollary of the proof of Theorem \ref{teo-fol-nondeg2}, we obtain the two-dimensional case of Kouchnirenko's result \cite[Th. I]{kouchnirenko1976}:
\begin{cor}
\label{cor-mu-nu-curves}
Let $f  \in \C\{x,y\}$ be a convenient power series. Then
\[  \nu_{0}(f) \leq \mu_{0}(f), \]
with equality holding if and only if   $f$ is Newton non-degenerate.
\end{cor}
\begin{proof}  The proof follows the same steps as that of Theorem
\ref{teo-fol-nondeg2}, taking $A = f_{x}$ and $B = f_{y}$.
\end{proof}


\section{Three basic transformations}
\label{sec-basic-transformations}

\subsection{Blow-up maps}
\label{sec-blow-ups}

Fix analytic coordinates   $(x,y)$   at $0 \in \C^{2}$.
Let $\sigma:(M,D) \to (\C^{2},0)$ be a quadratic blow-up at $0 \in \C^{2}$.
The blown-up space $M$ is covered by two coordinate charts, with coordinates $(\tilde{x}, \tilde{y})$ and $(\tilde{x}', \tilde{y}')$, such that
$\sigma(\tilde{x}, \tilde{y}) = (\tilde{x}, \tilde{x} \tilde{y})$ and
$\sigma(\tilde{x}', \tilde{y}') = (\tilde{x}' \tilde{y}', \tilde{y}')$, in such a way that, in each coordinate chart, the exceptional divisor
$D$ given by $\tilde{x}=0$ and $\tilde{y}' = 0$. We will focus on the coordinates $(\tilde{x}, \tilde{y})$ and say that  they are the  coordinates \emph{in
the first chart}.

Let  $\F$ be a foliation at the origin of $\C^{2}$, defined by the holomorphic $1$-form
$\omega = A(x,y) \dd x + B(x,y) \dd y$.
The foliation $\F$ is lifted by $\sigma$ to a germ of foliation in $(M,D)$,  the \emph{strict transform}
foliation, denoted as $\tilde{\F} = \sigma^{*} \F$. We say that $\sigma$ is a  \emph{dicritical} blow-up if
$xA_{m}(x,y) + y B_{m}(x,y) \equiv 0$, where $A_{m}$ and $B_{m}$ are the homogeneous parts of $A$ and $B$ of order $m = \ma_{0}(\F)$.
In the language of this paper, this is equivalent to asking that   $\Phi_{E_{1}}  \equiv 0$ or, in other words,
 that  the Newton polygon $\Delta_{\omega}$ is $1$-dicritical. In this situation, $D$ is not invariant by $\tilde{\F}$.
Otherwise, if $xA_{m}(x,y) + y B_{m}(x,y) \not\equiv 0$, or equivalently, if $\Phi_{E_{1}}  \not\equiv 0$, then $\sigma$ is a \emph{non-dicritical} blow-up. In this case, $D$ is invariant by $\tilde{\F}$.

Let us define, for the foliation $\F$,
\[  \bar{\ma}  =
\begin{cases}
\ma_{0}(\F), & \mbox{if $\sigma$ is non-dicritical } \\
\ma_{0}(\F) + 1, & \mbox{if $\sigma$ is dicritical}.
\end{cases}
\]
The pull-back $1$-form $\sigma^{*} \omega$ has $\tilde{x}$ as a factor of multiplicity  $\bar{\ma}$.
Canceling this factor, we obtain the
 $1$-form $\tilde{\omega} = \sigma^{*} \omega / \tilde{x}^{\bar{\ma}}$, the
\emph{strict transform} of $\omega$, which is a holomorphic $1$-form with isolated singularities, inducing $\tilde{\F}$ in the coordinates $(\tilde{x},\tilde{y})$. When the reference to the blow-up map is needed, we will also denote
$\tilde{\omega} = \sigma^{\hash} \omega$.
We warn the reader that $\sigma^{*}$ has different
meanings depending on the object it is applied to: for the foliation
$\F$, $\sigma^{*}\F$ already denotes the strict transform $\tilde{\F}$,
whereas for the defining $1$-form $\omega$, $\sigma^{*}\omega$ denotes
only the raw pull-back.

 We have
\[\tilde{\omega} = \sigma^{\hash} \omega     = \tilde{A}(\tilde{x},\tilde{y}) \dd \tilde{x} + \tilde{B}(\tilde{x},\tilde{y}) \dd \tilde{y} = \left(\frac{A(\tilde{x},\tilde{x}\tilde{y}) + \tilde{y} B(\tilde{x},\tilde{x}\tilde{y}) }{\tilde{x}^{\bar{\ma}}} \right) \dd \tilde{x}
+ \frac{  \tilde{x} B(\tilde{x},\tilde{x}\tilde{y})}{\tilde{x}^{\bar{\ma}}} \dd \tilde{y} .\]
Let us denote  the differential binomial  of bidegree $(\tilde{i},\tilde{j})$ of  $\tilde{\omega}$ by
$\tilde{\omega}_{\, \tilde{i} \, \tilde{j}}$.
The differential  binomial of $\omega$ of bidegree $(i,j)$
\begin{equation}
\label{eq-binomial0} \omega_{i j} = a_{ij} x^{i-1}y^{j} \dd x +  b_{ij} x^{i}y^{j-1} \dd y
\end{equation}
is transformed into the differential binomial
\begin{equation}
\label{eq-binomial}
 \tilde{\omega}_{i+j-\bar{\ma},j} = (a_{ij} + b_{ij}) \tilde{x}^{i+j- \bar{\ma} -1}\tilde{y}^{j} \dd \tilde{x} +  b_{ij} \tilde{x}^{i+j- \bar{\ma}}\tilde{y}^{j-1} \dd \tilde{y}
\end{equation}
of bidegree $(\tilde{i},\tilde{j}) = (i+j-\bar{\ma},j)$. Observe that $\omega_{i,j} \neq 0$ if and only if $\tilde{\omega}_{i+j-\bar{\ma},j} \neq 0$. Note also that the multipliers are related in the following way: $\lambda_{\, \tilde{i} \, \tilde{j}} = \lambda_{ij} -1$ (with the convention $\infty - 1 = \infty$).

Our analysis is then focused on the germ of $\tilde{\omega}$ at the origin of the coordinates $(\tilde{x}, \tilde{y})$.
If the coordinates
$(x,y)$ are  convenient or semi-convenient, then these new coordinates are either convenient,
if $\sigma$ is dicritical, or semi-convenient, if $\sigma$ is non-dicritical, since only in this latter case the exceptional divisor,
of equation $\tilde{x}=0$, is invariant by $\tilde{\F}$.
Note that,  since transversality
(non-invariance)   is preserved under strict transform,
the
axis  $\{y=0\}$ being non-invariant implies that its strict transform $\{\tilde y=0\}$ is also non-invariant.
We thus consider the Newton polygon $\Delta_{\tilde{\omega}}$ in these coordinates.

Let  $E_{\kappa}$ be an edge of $\Delta_\omega$, with $\kappa > 1$, having   a supporting line of  equation $ i + \kappa j = \delta$.
If  $(i,j) \in E_{\kappa}$ is a particle, then,
setting
$\tilde{i} = i+j- \bar{\ma}$ and  $\tilde{j} = j$, we find the relation $\tilde{i} + (\kappa-1) \tilde{j} = \delta - \bar{\ma}$.
Also, if $(i,j) \in \supp_{\omega} \setminus E_{\kappa}$, then $ i + \kappa j > \delta$ and also
$\tilde{i} + (\kappa-1) \tilde{j} > \delta - \bar{\ma}$.
This means that  $E_{\kappa}$  is transformed  by  $\sigma^{\hash}$ into an edge $\tilde{E}_{\kappa-1}$ of the
Newton polygon $\Delta_{\tilde{\omega}}$, of co-slope $\kappa -1$
(an edge, since $\kappa>1$ gives
$\ell(\tilde{E}_{\kappa-1})=\ell(E_\kappa)-h(E_\kappa)>0$) whose supporting line has   equation
$\tilde{i} + (\kappa-1) \tilde{j} = \delta - \bar{\ma}$. Since the $j$-coordinate is preserved in this process, $E_{\kappa}$ and $\tilde{E}_{\kappa-1}$
have both the same height and the same elevation. Denote $\tilde{E}_{\kappa-1} = \sigma^{\hash} E_{\kappa}$.

\begin{lem}
\label{lem-Ndeg} Let  $E_{\kappa}$ be an edge of $\Delta_\omega$ of co-slope $\kappa > 1$. Then,
$\omega$ is Newton $E_{\kappa}$-degenerate if and only if $\tilde{\omega}$ is Newton $\tilde{E}_{\kappa-1}$-degenerate.
\end{lem}
\begin{proof} To simplify notation, set $E =  E_{\kappa}$ and $\tilde{E} = \tilde{E}_{\kappa-1}$.
Let us define
\[ \sigma^{\hash} (Ax)_{E} =   \frac{\sigma^{*} (A x)_{E}}{ \tilde{x}^{\bar{\ma}}} =
  \sum_{(i,j) \in P(E)} a_{ij} \tilde{x}^{i+j- \bar{\ma}} \tilde{y}^{j}\]
and
\[ \sigma^{\hash} (B y)_{E} =   \frac{\sigma^{*} (B y)_{E}}{ \tilde{x}^{\bar{\ma}}} =   \sum_{(i,j) \in P(E)} b_{ij} \tilde{x}^{i+j- \bar{\ma}} \tilde{y}^{j}.\]
Hence, by \eqref{eq-binomial}, we have
\[  (\tilde{A}\tilde{x})_{\tilde{E}}  = \sigma^{\hash} (A x)_{E} + \sigma^{\hash} (B y)_{E} \qquad \text{and} \qquad
(\tilde{B}\tilde{y})_{\tilde{E}}  = \sigma^{\hash} (B y)_{E}.
\]
Considering this, if $\phi(x,y) \in \C[x,y]$ is a non-trivial common factor for $(A x)_{E}$ and $(B y)_{E}$,
then  $\tilde{\phi}(\tilde{x},\tilde{y}) = \phi(\tilde{x},\tilde{x}\tilde{y}) / \tilde{x}^{r}$, where $r$ is the order of $\phi$,
is a common factor for $\sigma^{\hash}(A x)_{E}$ and $\sigma^{\hash} (B y)_{E}$ and, hence, also
for  $(\tilde{A}\tilde{x})_{\tilde{E}}$ and  $(\tilde{B}\tilde{y})_{\tilde{E}}$.
Conversely, if a non-trivial   $\tilde{\phi}(\tilde{x},\tilde{y}) \in \C[\tilde{x},\tilde{y}]$ is a common factor for
$(\tilde{A}\tilde{x})_{\tilde{E}}$ and  $(\tilde{B}\tilde{y})_{\tilde{E}}$, it is also a non-trivial
common factor for $\sigma^{\hash}(A x)_{E}$ and $\sigma^{\hash} (B y)_{E}$ and we obtain a common factor
for $(A x)_{E}$ and $(B y)_{E}$
by setting  $\phi(x,y) = x^{s} \tilde{\phi}(x, y/x)$, where $s$ is the degree  of $\tilde{\phi}$ as a polynomial in $\tilde{y}$.
\end{proof}

Let $(x,y)$ be a local system of analytic coordinates at $0 \in \C^{2}$.
We   say that
a sequence of blow-ups is \emph{horizontal} with respect to $(x,y)$ if
their centers lie at the base points of the successive strict transforms of  the $x$-axis.
By \emph{base point} of the strict transform of a formal curve  by a blow-up map we mean
its intersection with the exceptional divisor.
In this way, if $\pi$ is a horizontal sequence of $\ell$ blow-ups, by taking, in each single blow-up,
first chart coordinates, $\pi$ is expressed in coordinates $(\tilde{x},\tilde{y})$ such that
$(x,y) = \pi(\tilde{x},\tilde{y}) = (\tilde{x}, \tilde{x}^{\ell} \, \tilde{y})$.
Denote by $\pi^{\hash}\omega$
the iterated strict transform of $\omega$ along $\pi$, obtained by
successively applying the strict transform $\sigma^{\hash}$ of each individual
blow-up $\sigma$ in the composition.
Let $E_{\kappa} \subset \Delta_{\omega}$ be an edge of co-slope $\kappa > \ell$.
Applying the discussion above inductively, one
blow-up at a time, we have that, at each step, $E_{\kappa}$ loses one
unit of co-slope,   its height and its elevation remain unchanged, whereas  the
multiplier of each of its particles is modified by a
shift $\lambda \mapsto \lambda-1$. In particular, if  $E_\kappa$ is dicritical, equivalently, if all
its particles   have multiplier equal to $\kappa$, then its successive transforms   will
be dicritical,
with the common value of the multiplier decreasing by $1$ at each blow-up.
Hence,  $E_\kappa$
 is transformed by $\pi^{\hash}$ into an edge $\tilde{E}_{\tilde{\kappa}} \subset \Delta_{\tilde{\omega}}$,
where  $\tilde{\omega} = \pi^{\hash}\omega$,
of co-slope $\tilde{\kappa} = \kappa - \ell$, having the same height and elevation.
In particular, if $\kappa = \ell + 1$ and $\Delta_{\omega}$ is $\kappa$-dicritical, then
$\Delta_{\tilde{\omega}}$ is $1$-dicritical. In this case, after applying the blow-ups of $\pi$,
 the next horizontal blow-up, centered at the origin of the coordinates $(\tilde{x},\tilde{y})$, is dicritical.

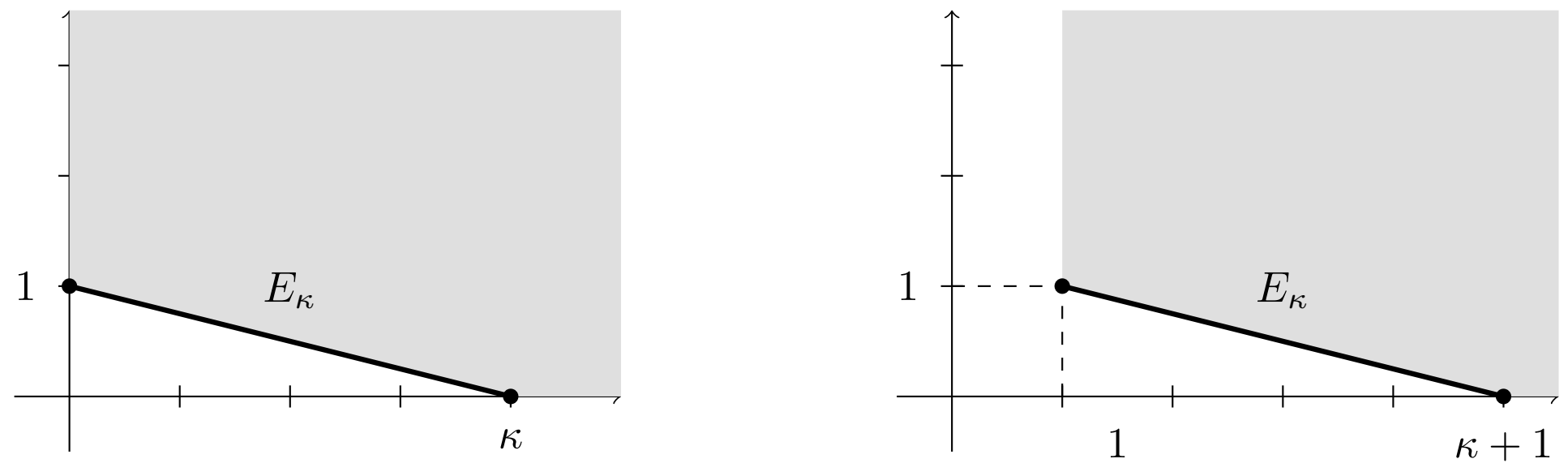
\begin{figure}[h]
\centering
\begin{tikzpicture}[scale=0.9]

\begin{scope}
  \draw[->] (-0.5,0) -- (5,0);
  \draw[->] (0,-0.5) -- (0,3.5);
  \foreach \x in {1,2,3,4} \draw (\x,-0.1) -- (\x,0.1);
  \foreach \y in {1,2,3} \draw (-0.1,\y) -- (0.1,\y);

  \fill[gray!25] (0,1) -- (0,3.5) -- (5,3.5) -- (5,0) -- (4,0) -- (0,1) -- cycle;

  \draw[very thick] (0,1) -- (4,0);

  \fill (0,1) circle (2pt);
  \fill (4,0) circle (2pt);

  \node[left] at (-0.15,1) {\small $1$};
  \node[below] at (4,-0.15) {\small $\kappa$};
  \node[above] at (2,0.65) {\small $E_{\kappa}$};

  \node[below] at (2,-1) {\small (a) Regular};
\end{scope}

\begin{scope}[xshift=8cm]
  \draw[->] (-0.5,0) -- (5.5,0);
  \draw[->] (0,-0.5) -- (0,3.5);
  \foreach \x in {1,2,3,4,5} \draw (\x,-0.1) -- (\x,0.1);
  \foreach \y in {1,2,3} \draw (-0.1,\y) -- (0.1,\y);

  \fill[gray!25] (1,1) -- (1,3.5) -- (5.5,3.5) -- (5.5,0) -- (5,0) -- cycle;

  \draw[dashed] (0,1) -- (1,1);
  \draw[dashed] (1,0) -- (1,1);

  \draw[very thick] (1,1) -- (5,0);

  \fill (1,1) circle (2pt);
  \fill (5,0) circle (2pt);

  \node[left] at (-0.15,1) {\small $1$};
  \node[below] at (1.5,-0.15) {\small $1$};
  \node[below] at (5,-0.15) {\small $\kappa+1$};
  \node[above] at (3,0.65) {\small $E_{\kappa}$};

  \node[below] at (3,-1) {\small (b) Singular};
\end{scope}

\end{tikzpicture}
\caption{Primary Newton polygons: regular and singular}
\label{fig: reg-sing polygons}
\end{figure}

Let   $\omega$ be the defining $1$-form of a foliation $\F$ and  $(x,y)$ be a system of
analytic coordinates.
We say that the Newton polygon $\Delta_{\omega}$   is \emph{primary} if it
is formed by a unique edge of height one and elevation zero. If $E_{\kappa}$ is this edge, of co-slope $\kappa$,
then there are two possibilities for its two vertices:  either they are $(0,1)$ and $(\kappa,0)$, and the coordinates are convenient, or they
are $(1,1)$ and $(\kappa+1,0)$, and the coordinates are semi-convenient. In the first case,
$\Delta_{\omega}$   is said to be \emph{regular} and, in  the second case, it is  called \emph{singular}
(see Figure \ref{fig: reg-sing polygons}).
For the regular case
\begin{equation}
\label{eq-regular-elementary}
(\omega)_{E_{\kappa}} = \omega_{0,1} + \omega_{\kappa,0}
=  b_{0,1}\dd y   +  a_{\kappa,0}x^{\kappa - 1}   \dd x ,
\end{equation}
where $b_{0,1}, a_{\kappa,0} \neq 0$ and,
thus, $\F$ is regular at $(x,y)=(0,0)$ and transverse to the $y$-axis.

Conversely, in convenient coordinates $(x,y)$, if $\F$ is regular
at $0 \in \C^{2}$ and transverse to the $y$-axis, then
$B(0,0)\neq0$ forces the vertex of $\Delta_\omega$
on the $y$-axis to be exactly $(0,1)$  \and,
since $(x,y)$ are convenient, $\Delta_\omega$
also reaches the horizontal axis, so its total height is $h =1$, forcing a single edge. Thus, $\Delta_\omega$
is  primary regular.
In the singular   case,   the $y$ axis is invariant since the coordinates are semi-convenient. In this case, we have
\begin{equation}
\label{eq-singular-elementary}
(\omega)_{E_{\kappa}} = \omega_{1,1} + \omega_{\kappa+1,0}
=  ( a_{1,1}y\dd x +  b_{1,1} x \dd y )  +  a_{\kappa+1,0}x^{\kappa}   \dd x,
\end{equation}
where $a_{\kappa+1,0} \neq 0$ and at least one out of $a_{1,1}$ and  $b_{1,1}$ is non-zero,
and $\F$ has a singularity at the origin.
Denote $\lambda = \lambda_{1,1}$  the \emph{multiplier} of $\Delta_{\omega}$.
If either $\lambda = 0$ or $\lambda = \infty$,
corresponding, respectively, to  $a_{1,1} =0$ or $b_{1,1}=0$,
then $\F$ has a saddle-node singularity at $0 \in \C^{2}$.
In the first case, the $y$-axis corresponds to weak direction and, in the second case, to the strong direction.
When both $a_{1,1}$ and $b_{1,1}$ are non-zero, we have $\lambda = \lambda_{1,1}= -a_{1,1}/b_{1,1}$.
In this case, if $\lambda \not\in \Q_{>0}$, then $\F$ has a simple non-degenerate
singularity.
If $\lambda  \in \Q_{>0}$, we   say that  $\Delta_{\omega}$    is \emph{resonant}.
This situation may be related to dicritical blow-ups and the existence of  infinitely many separatrices.

However, we remark that the foliation $\F$ defined by the Poincar\'e-Dulac normal form $\omega = (ny + \zeta x^{n}) \dd x - x \dd y$, where $\lambda = n \geq 2$ and
$\zeta \in \C^{*}$, is such that  $\Delta_{\omega}$
is primary resonant. The special feature here is that $\F$   has the $y$-axis as its unique separatrix and, therefore,  no separatrix is associated with the edge
$E_{n}$ of $\Delta_{\omega}$ of co-slope $\kappa = n$.

If $\Delta_{\omega}$ is regular primary, then it is clearly Newton non-degenerate.
In the   singular primary case, we have  a characterization of Newton degeneracy:
\begin{prop}
\label{prop-primary}
If $\Delta_{\omega}$ is singular primary, then $\omega$ is Newton degenerate if and only if the origin is a saddle-node
singularity of $\F$
whose strong separatrix is contained in the $y$-axis.
\end{prop}
\begin{proof}
The $1$-form  $\omega$, over   the
single edge $E_{\kappa}$ of $\Delta_{\omega}$, is given by
\eqref{eq-singular-elementary}.
The corresponding degeneracy equations  are, hence,
\[ (xA)_{E_{\kappa}} = x(a_{1,1}   y  +a_{\kappa+1,0}x^{\kappa} )
\qquad \text{and} \qquad (yB)_{E_{\kappa}} =  y  ( b_{1,1} x) .\]
Thus, $\omega$ is   Newton $E_{\kappa}$-degenerate
if and only if $b_{1, 1} = 0$ and, in this case, $a_{1, 1}, a_{\kappa+1,0} \neq 0$, and the non-trivial common factor
of $(xA)_{E_{\kappa}}$ and  $(yB)_{E_{\kappa}}$
is $a_{1, 1}y + a_{\kappa+ 1,0}x^{\kappa}$.
Recalling the discussion above, $b_{1,1}=0$ is
precisely the condition for $0 \in \C^{2}$ to be a saddle-node singularity
with the $y$-axis as its strong separatrix, which proves the
proposition.
\end{proof}

\begin{figure}[h]
\centering
\begin{tikzpicture}[scale=0.9]
  \draw[->] (-0.5,0) -- (7.5,0);
  \draw[->] (0,-0.5) -- (0,4.5);
  \foreach \x in {1,2,3,4,5,6,7} \draw (\x,-0.1) -- (\x,0.1);
  \foreach \y in {1,2,3,4} \draw (-0.1,\y) -- (0.1,\y);

  \fill[gray!25] (1,3) -- (1,4.5) -- (7.5,4.5) -- (7.5,0) -- (7,0) -- (3,1) -- (1,3) -- cycle;

  \draw[very thick] (1,3) -- (3,1);
  \draw[very thick] (3,1) -- (7,0);

  \draw[dashed] (0,1) -- (3,1);
  \draw[dashed] (3,0) -- (3,1);

  \fill (1,3) circle (2pt);
  \fill (3,1) circle (2pt);
  \fill (7,0) circle (2pt);

  \node[left] at (-0.15,1) {\small $1$};
  \node[below] at (3,-0.15) {\small $i_0$};
  \node[below] at (6.5,-0.15) {\small $i_0+\kappa$};
  \node[above] at (1.9,2.5) {\small $E_1$};
  \node[above] at (4.7,0.75) {\small $E_{\kappa}$};

\end{tikzpicture}
\caption{Pre-primary Newton polygon}
\label{fig: pre-primary}
\end{figure}
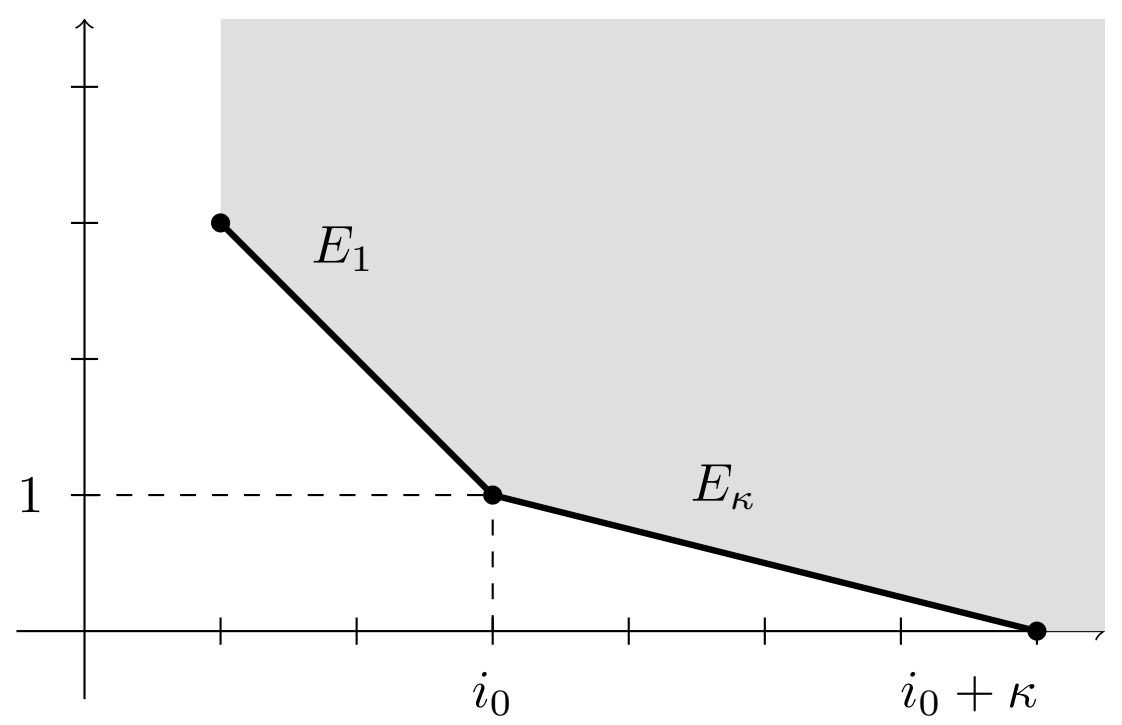

The Newton polygon  $\Delta_{\omega}$ is said to be \emph{pre-primary}
if it has exactly one edge $E_{\kappa}$ of co-slope $\kappa>1$, and
this edge satisfies $h(E_{\kappa}) = 1$ and $e(E_{\kappa}) = 0$ (see Figure \ref{fig: pre-primary}).
In this case, we have that $\kappa \in \N$. We remark that   $\Delta_{\omega}$ may have other edges of co-slope one or less.
 The edge $E_{\kappa}$ will be called \emph{distinguished}.
The only particles of $E_{\kappa}$ are its vertices,  $\Top(E_{\kappa}) = (i_0,1)$ and
$\Bot(E_{\kappa}) = (i_0 + \kappa,0)$, for some $i_{0} \in \N$.   Note that $\Top(E_{\kappa}) \in E_{1}$, the
element of co-slope one, and,
hence, $i_0 = \ma_{0}(\F)$.
Keeping the terminology above, the \emph{multiplier} of $\Delta_{\omega}$ is $\lambda = \lambda_{i_0,1}$, i.e. the multiplier of $\Top(E_{\kappa})$, and
 $\Delta_{\omega}$ is said to be   \emph{resonant} if $\lambda \in \Q_{>0}$.

\begin{prop}
\label{prop-preprimary-blowup}
Suppose that $\Delta_{\omega}$ is pre-primary. Let $\sigma $ be a blow-up as above,
and $\tilde{\omega} = \sigma^{\hash} \omega$ be the strict transform of $\omega$. Then
$\Delta_{\tilde{\omega}}$ is primary, being regular primary if and only if $\sigma$ is   dicritical, i.e. if
and only if $E_{1}$ is  dicritical.
\end{prop}
\begin{proof} Considering equations \eqref{eq-binomial0} and \eqref{eq-binomial}, the   differential
binomials $\omega_{i_0,1}$ and $\omega_{i_0 + \kappa,0}$ correspond, by $\sigma$, to
$\tilde{\omega}_{i_0+1 - \bar{\ma},1}$ and  $\tilde{\omega}_{i_0+\kappa- \bar{\ma},0}$, where
$\bar{\ma} = \ma_{0}(\F)$ if $\sigma$ is non-dicritical or $\bar{\ma} = \ma_{0}(\F) + 1$ if $\sigma$ is  dicritical.
Hence, $\Delta_{\tilde{\omega}}$ is formed by a single edge $\tilde{E}_{\kappa -1}$,
of vertices
$\Top(\tilde{E}_{\kappa -1}) = (i_0+1 - \bar{\ma},1)$ and $\Bot(\tilde{E}_{\kappa -1})= (i_0+\kappa- \bar{\ma},0 )$. This edge has   co-slope $\kappa -1$ (positive, since $\kappa>1$ by the pre-primary
hypothesis),    height $h(\tilde{E}_{\kappa -1}) = 1$  and elevation $e(\tilde{E}_{\kappa -1}) = 0$.
The  conclusion of the proposition then follows from the fact that $i_0 = \ma_{0}(\F)$, resulting in   $\Delta_{\omega}$ singular primary
in the non-dicritical case and regular primary
in the  dicritical case.
\end{proof}

Regarding the dicritical case of Proposition \ref{prop-preprimary-blowup},
if $D$ is the exceptional divisor of the blow-up $\sigma$, then $\tilde{\F}$ is regular and transverse to $D$
at all its points visible in the coordinates
$(\tilde{x},\tilde{y})$, that is, away from the point at infinity of
this chart (corresponding to the vertical direction, the origin in the second chart $(\tilde{x}',\tilde{y}')$).
 Indeed, it suffices to take into account Equation \eqref{eq-regular-elementary}
with coordinates $(\tilde{x}, \tilde{y})$ in place of $(x,y)$, knowing that $\tilde{y} \in \C$ parameterizes the points of $D$ considered.

In view of the Poincar\'e-Dulac normal form mentioned earlier,
we say that the  pre-primary Newton polygon $\Delta_{\omega}$, with  distinguished edge   $E_{\kappa}$, is \emph{Poncar\'e-Dulac} if its multiplier is $\kappa$.
This object has the following property:
\begin{lem}
\label{lem-poincare-dulac}
If the Newton polygon is Poincar\'e-Dulac, then the foliation $\F$ has no separatrices with horizontal tangent.
\end{lem}
\begin{proof}
Suppose, by contradiction,  that $\F$ admits a separatrix $S$ with horizontal tangent. Write one of its  Puiseux series in the form $y = c x^{k} + \cdots$,
where $k \in \Q$, $k >1$ and $c \in \C^{*}$.
We claim that $k = \kappa$, where $\kappa \in \N$ is the co-slope of the distinguished edge $E_{\kappa} \subset \Delta_{\omega}$.
Indeed, we cannot have $k < \kappa$ since, in this case,   $\Delta_{\omega}$ would be $k$-dicritical, which is not the case,
considering that the multiplier of $\Top(E_{\kappa})$ is $\kappa$.
Likewise, we cannot have $k>\kappa$ because, in this
case, $E_k$ would be the vertex $\Bot(E_\kappa)$, lying on the
horizontal axis, whose multiplier is $\infty\neq k$ for every
finite $k$ and, so, $E_k$ cannot be dicritical.
We note that the blow-up at $0 \in \C^{2}$ is non-dicritical. In fact,  $\Delta_{\omega}$ is pre-primary, so  $\Bot(E_{1}) = \Top(E_{\kappa})$ and this vertex has multiplier $\kappa >1$, ensuring that  $E_{1}$ is non-dicritical.
For this same reason, if $\pi$ denotes a sequence  of $\kappa -1$ horizontal blow-ups, all its individual blow-ups   will be non-dicritical.
Let   $(\tilde{x},\tilde{y})$ be the coordinates introduced by $\pi$ and denote by $\tilde{\F} = \pi^{*} \F$ the strict transform foliation,
defined by the $1$-form $\tilde{\omega} = \pi^{\hash} \omega$. The curve $\tilde{S} = \pi^{*}S$   is a separatrix of $\tilde{\F}$
 whose tangent line has
slope $c$. Since also $\tilde{x} =0$ is $\tilde{\F}$-invariant, the origin in the coordinates $(\tilde{x},\tilde{y})$
 is a singularity of $\tilde{\F}$.
Therefore, the Newton polygon $\Delta_{\tilde{\omega}}$ is singular primary, with edge $\tilde{E}_{1}$ of co-slope $\tilde{\kappa} = 1$
and  multiplier $\tilde{\lambda} = 1$.
We write $\tilde{\omega}= (b \tilde{x} + a \tilde{y} )\dd \tilde{x} - a \tilde{x} \dd \tilde{y} + \cdots$, where $a,b \neq 0$
and the dots stand for terms of higher degrees. The   initial polynomial of $\tilde{\omega}$ with respect  $\tilde{E}_{1}$ is
$\Phi_{\tilde{E}_{1}}  = b$. It has no roots and, hence,
$\tilde{\F}$ does not have  $\tilde{E}_{1}$-separatrices, which is an absurd given  the existence of $\tilde{S}$.
This contradiction implies that  $\F$  has no $E_{\kappa}$-separatrices.
\end{proof}
Recall that a pre-primary Newton polygon  may contain an edge of co-slope one.
Hence, in  Lemma \ref{lem-poincare-dulac}, $\F$ could in principle have separatrices with non-horizontal tangent lines associated
with this edge.

\subsection{Horizontalization procedure}
\label{sec-horizontalization}

As before, $\F$ is the foliation of equation $\omega = 0$, of which we select a  separatrix   $S$.
Our objective is to decide whether $\Delta_{\omega}$ is Newton degenerate, relying on information gathered along the base points
of the strict transforms of $S$ along the reduction of singularities of $\F$. In view of the computations in the preceding subsection, in order to have control over the evolution of the Newton polygon of $\F$ along these points, we will seek coordinates  with respect to which this sequence of blow-ups is horizontal.
In this process, we will need the \emph{horizontalization procedure}, that we describe
next.

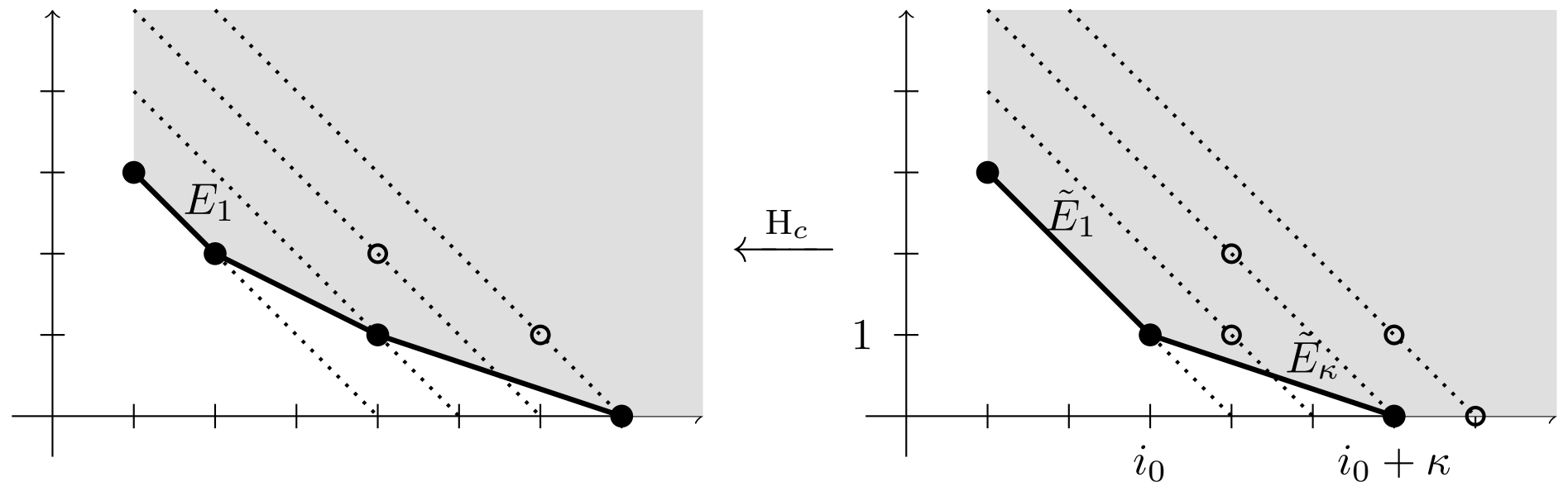
\begin{figure}[h]
\centering
\begin{tikzpicture}[scale=0.65]

\begin{scope}
  \draw[->] (-0.5,0) -- (8,0);
  \draw[->] (0,-0.5) -- (0,5);
  \foreach \x in {1,...,7} \draw (\x,-0.15) -- (\x,0.15);
  \foreach \y in {1,...,4} \draw (-0.15,\y) -- (0.15,\y);

  \fill[gray!25] (1,3) -- (1,5) -- (8,5) -- (8,0) -- (7,0)
    -- (4,1) -- (2,2) -- (1,3)--cycle;


  \draw[dotted, thick] (2,2) -- (4,0);
  \draw[dotted, thick] (1,4) -- (5,0);
  \draw[dotted, thick] (1,5) -- (6,0);
  \draw[dotted, thick] (2,5) -- (7,0);

  \draw[very thick] (1,3)--(2,2) -- (4,1) -- (7,0);

  \fill (1,3) circle (4pt);
  \fill (2,2) circle (4pt);
  \fill (4,1) circle (4pt);
  \fill (7,0) circle (4pt);

  \draw[thick] (4,2) circle (3pt);
  \draw[thick] (6,1) circle (3pt);

  \node[above right] at (1.4,2.2) {\small $E_1$};
\end{scope}

\node at (9,2.25) { $\xleftarrow{\ {\rm H}_c\ }$};

\begin{scope}[xshift=10.5cm]
  \draw[->] (-0.5,0) -- (8,0);
  \draw[->] (0,-0.5) -- (0,5);
  \foreach \x in {1,...,7} \draw (\x,-0.15) -- (\x,0.15);
  \foreach \y in {1,...,4} \draw (-0.15,\y) -- (0.15,\y);
  \fill[gray!25] (1,3) -- (1,5) -- (8,5) -- (8,0) -- (6,0)
    -- (3,1) -- (1,3) -- cycle;
  \draw[dotted, thick] (3,1) -- (4,0);
  \draw[dotted, thick] (1,4) -- (5,0);
  \draw[dotted, thick] (1,5) -- (6,0);
  \draw[dotted, thick] (2,5) -- (7,0);
  \draw[very thick] (1,3) -- (3,1) -- (6,0);
  \fill (1,3) circle (4pt);
  \fill (3,1) circle (4pt);
  \fill (6,0) circle (4pt);
  \draw[thick] (4,1) circle (3pt);
  \draw[thick] (4,2) circle (3pt);
  \draw[thick] (6,1) circle (3pt);
  \draw[thick] (7,0) circle (3pt);
  \node[left] at (-0.2,1) {\small $1$};
  \node[below] at (3,-0.15) {\small $i_0$};
  \node[below] at (6,-0.15) {\small $i_0+\kappa$};
  \node[above right] at (1.5,2) {\small $\tilde{E}_1$};
  \node[above] at (5,0.25) {\small $\tilde{E}_{\kappa}$};
\end{scope}

\end{tikzpicture}
\caption{Horizontalization procedure}
\label{fig:horizontalization}
\end{figure}

Suppose that  the separatrix $S$   has the oblique line $y - cx=0$ as a tangent direction, i.e., its Puiseux series are of
the form $y = cx + \cdots$, where $c \in \C^{*}$.
Let us do the coordinate change $(x,y) = {\rm H}_{c}(\bar{x},\bar{y}) = (\bar{x}, \bar{y} + c\bar{x})$, that takes the oblique
direction  $y - cx=0$ to the horizontal direction $\bar{y}= 0$ and, accordingly, $S$ to   the curve $\bar{S}= {\rm H}_{c}^{*} S$  with horizontal tangent line.
 If, by chance,  the line $y - cx=0$ is   invariant by $\F$, then we replace $ {\rm H}_{c}$ by its  composition with a diffeomorphism highly
tangent to the identity in order to have, at least, semi-convenient coordinates. For instance, we can take ${\rm H}_{c}(\bar{x},\bar{y}) = (\bar{x}, \bar{y} + c\bar{x} + \bar{x}^{\ell})$,
where $\ell \in \N$ is sufficiently large.
Anyhow, abusing notation, we will keep considering ${\rm H}_{c}(\bar{x},\bar{y}) = (\bar{x}, \bar{y} + c\bar{x})$ in our calculations, since they involve only  terms of low order, and these are unaffected by the extra term $\bar{x}^{\ell}$ if $\ell$ is sufficiently large.
Writing  $\omega = A(x,y) \dd x + B(x,y) \dd y$ and
  $\bar{\omega} = {\rm H}_{c}^{*} \omega$, we have
\begin{equation}
\label{eq-omegabar}
\bar{\omega} =    (A(\bar{x},\bar{y} + c \bar{x} ) + c B(\bar{x},\bar{y} + c \bar{x} ) ) \dd \bar{x}
+ B(\bar{x},\bar{y} + c \bar{x} ) \dd \bar{y} = \bar{A}(\bar{x},\bar{y}) \dd \bar{x} + \bar{B}(\bar{x},\bar{y}) \dd \bar{y}.
\end{equation}
Consider the element $E_{1} \subset \Delta_{\omega}$ and let $i + j = \delta$ be the equation of the supporting line of co-slope $\kappa = 1$.
A differential binomial
$ \omega_{ij} = a_{ij} x^{i-1}y^{j} \dd x +  b_{ij} x^{i}y^{j-1} \dd y $,
corresponding to a particle $(i,j)$ of $E_{1}$, produces, under this change of coordinates,
\begin{equation}
\label{eq-omegabar1} ( a_{ij} \bar{x}^{i-1}(\bar{y} + c \bar{x}) ^{j} +  c b_{ij} \bar{x}^{i}(\bar{y} + c \bar{x}) ^{j-1} ) \dd \bar{x} +
b_{ij} \bar{x}^{i}(\bar{y} + c \bar{x})^{j-1} \dd \bar{y}.
\end{equation}
Developing the binomial terms, we find a sum of differential binomials of bidegrees $(\bar{i},\bar{j})$
satisfying   $\bar{i} + \bar{j} = \delta$, where $\bar{j} \leq j$. Varying $(i,j)$ over $E_{1}$ and assembling differential
binomials of same bidegrees, we get differential binomials $\bar{\omega}_{\bar{i} \, \bar{j} }$ of $\bar{\omega}$, where $\bar{i} + \bar{j} = \delta$.
Thus, they generate the element of co-slope $\kappa = 1$ of $\Delta_{\bar{\omega}}$, say $\bar{E}_{1}$ (see Figure~\ref{fig:horizontalization}).
In \eqref{eq-omegabar1}, we see that $\Top(\bar{E}_{1}) = \Top(E_{1})$. Moreover, $\Top(\bar{E}_{1})$ has the same
multiplier $\lambda = -a_{ij}/b_{ij}$ as $\Top(E_{1})$.  However, some particles may be created or suppressed in this process.
For instance, since $c \in \C^{*}$ is the coefficient of the linear term of the Puiseux series of the invariant curve $S$, we have $\Phi_{E_{1}}(c) = 0$ and, then,
\[ \bar{\omega}_{\delta 0 } =  \left( \sum_{(i,j) \in P(E_{1})}  (a_{ij}  +   b_{ij}) c^{j} \right) \bar{x}^{\delta -1} \dd \bar{x}   =   \Phi_{E_{1}}(c) \bar{x}^{\delta -1} \dd \bar{x} = 0 .\]
Also, taking into account \eqref{eq-principal-polynomial},
\begin{equation}
\label{eq-elevation1} \bar{\omega}_{\delta -1, 1} =
\left(\sum_{(i,j) \in  P(E_{1})} j  a_{ij} c^{j-1} + c (j-1) b_{ij} c^{j-2}\right) \bar{x}^{\delta -2} \bar{y}  \dd \bar{x}
+  \left(\sum_{(i,j) \in  P(E_{1})}  b_{ij} c^{j-1} \right) \bar{x}^{\delta -1} \dd \bar{y}
\end{equation}
\[  =
\left( \sum_{(i,j) \in  P(E_{1})} j ( a_{ij} +  b_{ij}) c^{j-1} - \sum_{(i,j) \in  P(E_{1})} b_{ij} c^{j-1}\right) \bar{x}^{\delta -2} \bar{y}  \dd \bar{x}
+  \left(\sum_{(i,j) \in  P(E_{1})}  b_{ij} c^{j-1} \right) \bar{x}^{\delta -1} \dd \bar{y}   \]
\[ = ( \Phi_{E_{1}}'(c)  - \beta_{1}(c))  \bar{x}^{\delta -2} \bar{y}  \dd \bar{x}
+   \beta_{1}(c) \bar{x}^{\delta -1} \dd \bar{y} .  \]
Therefore,    $\bar{E}_{1}$ has elevation greater than one if and only if
the pair of equations
\begin{equation}
\label{eq-elevation1*}
 \Phi_{E_{1}}'(c)  = 0
\qquad \text{and} \qquad  \beta_{1}(c) = 0
\end{equation}
is satisfied.

The outcome of this  horizontalization procedure can give us information on Newton degeneracy:
\begin{prop}
\label{prop-fact1}
If $E_{1}$ is an edge of $\Delta_{\omega}$ and $e(\bar{E}_{1}) > 1$, then $\omega$ is  Newton $E_{1}$-degenerate.
\end{prop}
\begin{proof}
We start by introducing the following factorization:
\begin{equation}
\label{eq-ABhash}
(xA(x,y))_{E_{1}} = x A^{*}_{E_{1}}(x,y)
 \qquad \text{and} \qquad
 (yB(x,y))_{E_{1}} = y B^{*}_{E_{1}}(x,y) ,
\end{equation}
where
\[  A^{*}_{E_{1}}(x,y) = \sum_{(i,j) \in E_{1}} a_{i,j} x^{i-1} y^{j}
 \qquad \text{and} \qquad
 B^{*}_{E_{1}}(x,y) = \sum_{(i,j) \in E_{1}} b_{i,j} x^{i} y^{j-1}. \]
Similarly,
\begin{equation}
\label{eq-ABbarhash}
(\bar{x} \bar{A}(\bar{x},\bar{y}))_{\bar{E}_{1}} = \bar{x} \bar{A}^{*}_{\bar{E}_{1}}(\bar{x},\bar{y})
 \qquad \text{and} \qquad
(\bar{y} \bar{B}(\bar{x},\bar{y}))_{\bar{E}_{1}} = \bar{y} \bar{B}^{*}_{\bar{E}_{1}}(\bar{x},\bar{y}) ,
\end{equation}
with $\bar{A}^{*}_{\bar{E}_{1}}$ and $\bar{B}^{*}_{\bar{E}_{1}}$ defined in an equivalent way.
From Equation \eqref{eq-omegabar1} we have
\begin{equation}
\label{eq-Hc*} {\rm H}_{c}^{*} A^{*}_{E_{1}}(x,y) = \bar{A}^{*}_{\bar{E}_{1}} - c \bar{B}^{*}_{\bar{E}_{1}}
 \qquad \text{and} \qquad
{\rm H}_{c}^{*} B^{*}_{E_{1}}(x,y) = \bar{B}^{*}_{\bar{E}_{1}} .
\end{equation}

Now, invoking our hypothesis, since $e(\bar{E}_{1}) > 1$, all differential binomials $\bar{\omega}_{\bar{i} \, \bar{j} }$, for
 $(\bar{i},\bar{j}) \in \bar{E}_{1}$, have $\bar{y}$ as a factor and, hence,
$\bar{y}$ is a factor of $(\bar{x} \bar{A}(\bar{x},\bar{y} ))_{\bar{E}_{1}}$
 and $\bar{y}^2$ is
a factor of $(\bar{y} \bar{B}(\bar{x},\bar{y}))_{\bar{E}_{1}}$.
Thus, by \eqref{eq-ABbarhash}, $\bar{y}$ is a factor of both $\bar{A}^{*}_{\bar{E}_{1}}(\bar{x},\bar{y})$
and $\bar{B}^{*}_{\bar{E}_{1}}(\bar{x},\bar{y})$ and, by \eqref{eq-Hc*},
also of ${\rm H}_{c}^{*} A^{*}_{E_{1}}(x,y)$ and  ${\rm H}_{c}^{*} B^{*}_{E_{1}}(x,y)$.
This allows us to conclude that $y-cx$ is a common factor of  $A^{*}_{E_{1}}(x,y)$ and  $B^{*}_{E_{1}}(x,y)$ and
hence, by \eqref{eq-ABhash}, also of  $(xA(x,y))_{E_{1}}$ and  $(yB(x,y))_{E_{1}}$.
\end{proof}

If $e(\bar{E}_{1}) = 1$,
since $(\bar x,\bar y)$ are at least semi-convenient,
$\Delta_{\bar\omega}$ reaches the horizontal axis. Thus,
 $\Delta_{\bar{\omega}}$
contains only one edge to the right of $\bar{E}_{1}$, the edge $\bar{E}_{\bar{\kappa}}$ of co-slope $\bar{\kappa} \in \N$, with $\bar{\kappa} >1$.  That is, $\Delta_{\bar{\omega}}$ is pre-primary according to the terminology  introduced in \S\ref{sec-blow-ups}.
Also in this case, something can be said on Newton degeneracy:
\begin{prop}
\label{prop-fact2}
Suppose that  $E_{1}$ is an edge of $\Delta_{\omega}$. If $e(\bar{E}_{1}) = 1$ and
 $\bar{\omega}$ is  Newton $\bar{E}_{\bar{\kappa}}$-degenerate, then $\omega$ is
 Newton $E_{1}$-degenerate.
\end{prop}
\begin{proof}
We have that $\bar{E}_{\bar{\kappa}}$ contains only two particles, say $\Top(\bar{E}_{\bar{\kappa}}) = (i_{0},1)$ and
$\Bot(\bar{E}_{\bar{\kappa}}) = (i_{0} + \bar{\kappa},0)$,  supporting the following   differential monomials:
\[\bar{\omega}_{i_{0}, 1} = a_{i_{0}, 1} \bar{x}^{i_{0} -1} \bar{y} \dd \bar{x} +   b_{i_{0}, 1} \bar{x}^{i_{0}}  \dd \bar{y}
\qquad \text{and} \qquad
\bar{\omega}_{i_{0} + \bar{\kappa},0} = a_{i_{0} + \bar{\kappa},0} \bar{x}^{i_{0} + \bar{\kappa} -1} \dd \bar{x}.
\]
Thus
\[ (\bar{x} \bar{A})_{\bar{E}_{\bar{\kappa}}} = a_{i_{0}, 1} \bar{x}^{i_{0}} \bar{y} + a_{i_{0} + \bar{\kappa},0} \bar{x}^{i_{0} + \bar{\kappa}} = \bar{x}^{i_{0}}( a_{i_{0}, 1}  \bar{y} + a_{i_{0} + \bar{\kappa},0} \bar{x}^{\bar{\kappa}} )
\quad \text{and} \quad
( \bar{y} \bar{B})_{\bar{E}_{\bar{\kappa}}} =  b_{i_{0}, 1} \bar{x}^{i_{0}} \bar{y}.
\]
Hence, $(\bar{x}\bar{A})_{\bar{E}_{\bar{\kappa}}}$ and $(\bar{y}\bar{B})_{\bar{E}_{\bar{\kappa}}}$ have a non-trivial factor if and only if
$b_{i_{0}, 1} = 0$.
We have, in this case, $\bar{\omega}_{i_{0}, 1} = a_{i_{0}, 1} \bar{x}^{i_{0} -1} \bar{y} \dd \bar{x}$, giving that $\bar{y}$
is a factor of all differential monomials supported on $\bar{E}_{1}$.
The argument then finishes, as in the proof of Proposition \ref{prop-fact1}, by undoing the change of coordinates ${\rm H}_{c}$.
 \end{proof}

The following is a consequence of the  calculations above:
 \begin{cor}
\label{cor-E1-dicritical-Hc}
Under the notation above,
$E_{1}$ is dicritical if and only if $\bar{E}_{1}$ is dicritical.
\end{cor}
\begin{proof}
Using equations \eqref{eq-ABhash}, \eqref{eq-ABbarhash} and
\eqref{eq-Hc*}, we have
\begin{eqnarray*}
  {\rm H}_{c}^{*} \left( (xA(x,y))_{E_{1}} + (yB(x,y))_{E_{1}} \right)
  &=&  {\rm H}_{c}^{*} \left( x A^{*}_{E_{1}}(x,y) + y B^{*}_{E_{1}}(x,y) \right) \\
    &=& \bar{x} {\rm H}_{c}^{*} A^{*}_{E_{1}}(x,y)
        + (\bar{y} + c \bar{x})  {\rm H}_{c}^{*} B^{*}_{E_{1}}(x,y)  \\
    &=&  \bar{x} \left( \bar{A}_{\bar{E}_{1}}^{*}(\bar{x},\bar{y})
        - c  \bar{B}^{*}_{\bar{E}_{1}}(\bar{x},\bar{y}) \right)
        + (\bar{y} + c \bar{x})  \bar{B}^{*}_{\bar{E}_{1}}(\bar{x},\bar{y}) \\
    &=&  \bar{x}   \bar{A}_{\bar{E}_{1}}^{*}(\bar{x},\bar{y})
        + \bar{y}   \bar{B}_{\bar{E}_{1}}^{*}(\bar{x},\bar{y}) \\
    &=&  ( \bar{x}   \bar{A}(\bar{x},\bar{y}))_{\bar{E}_{1}}
        +   ( \bar{y}   \bar{B}(\bar{x},\bar{y}))_{\bar{E}_{1}}.
\end{eqnarray*}
Thus, setting $x= \bar{x} = 1$ in the relation above, we have the
following comparison of initial polynomials:
\begin{equation}
\label{eq-dic-eq}
 \Phi_{E_{1}}(\bar{y} + c) = \Phi_{\bar{E}_{1}}(\bar{y}).
\end{equation}
This implies that $\Phi_{E_{1}} \equiv 0$ if and only if
$\Phi_{\bar{E}_{1}} \equiv 0$.
\end{proof}

Let us see what happens when  $E_{1}$ is a single vertex, associated with a separatrix $S$:

\begin{prop}
\label{prop-E1-vertex}
Suppose that the element $E_{1}$ of $\Delta_{\omega}$ is a dicritical  vertex.
Then  $\bar{E}_{1}$ is dicritical and has elevation $e(\bar{E}_{1}) =1$.
\end{prop}
\begin{proof}
Since $E_{1}$ is dicritical, so is $\bar{E}_{1}$  by Corollary \ref{cor-E1-dicritical-Hc}.
Now, $E_{1}$ is   a single vertex $(i,j)$ that,   being dicritical,   corresponds to coefficients
$a_{ij}$ and $b_{ij}$ that  are both non-zero.
From \eqref{eq-principal-polynomial},
we have $\beta_{1}(c) = b_{ij}c^{j-1} \neq 0$ and we finally get
$e(\bar{E}_{1}) =1$ from the conditions given by  \eqref{eq-elevation1*}.
\end{proof}

Let $S$ be an analytic branch at $(\C^{2},0)$.
Suppose that one of its Puiseux series is of the form $y = cx^{k} + c'x^{\ell} + \cdots$, where $c, c' \in \C^{*}$, $1<k<\ell$, with  $k \in \Z$ and
  $\ell \in \Q$, the dots indicating terms of higher order, as usual. A horizontal sequence of $k-1$  blow-ups, say $\pi_{k-1} = \sigma_{1} \circ \cdots, \sigma_{k-1}$, takes $S$
into a branch $\tilde{S} = \pi_{k-1}^{*}S$ having a Puiseux series $\tilde{y} = c \tilde{x}  + c'\tilde{x}^{\ell - k +1} + \cdots$. The horizontalization map $\mathrm{H}_{c}$ then takes $\tilde{S}$ into a branch with horizontal tangent line, having a Puiseux series with first term
$c'\tilde{x}^{\ell - k +1}$. Roughly speaking, through this procedure, the first monomial in the Puiseux series of $S$ is eliminated.
This will be used systematically in the sequel, motivating us to define
the map $\Upsilon_{k,c} =   \mathrm{H}_{c}  \circ \pi_{k-1}$, that we call \emph{initial monomial elimination}.
We have that $\pi_{k-1}(\tilde{x},\tilde{y}) = (\tilde{x},\tilde{x}^{k-1}\tilde{y}) = (x,y)$ and, writing
$\mathrm{H}_{c}(\bar{x},\bar{y}) = (\bar{x},\bar{y} + c\bar{x}) =  (\tilde{x},\tilde{y})$, we get
$\Upsilon_{k,c}(\bar{x},\bar{y}) = (\bar{x},\bar{x}^{k-1}(\bar{y} + c\bar{x})) = (x,y)$.
The next lemma will be used in our future arguments:
\begin{lem}
\label{lem-monomial-elimination}
 Let $\F$ be a germ of foliation at $(\C^{2},0)$, defined by the $1$-form $\omega$ in semi-convenient coordinates $(x,y)$,
admitting a separatrix $S$ with Puiseux series $y = cx^{k}   + \cdots$, where
$c \in \C^{*}$,  $k >1$ and   $k \in \Z$. Let
 $\bar{\F} = \Upsilon_{k,c}^{*}\F$,  defined in the coordinates $(\bar{x},\bar{y})$  by the $1$-form $\bar{\omega}$.
  If $\omega$ is Newton non-degenerate, then
$\Delta_{\bar{\omega}}$ is pre-primary.
\end{lem}
\begin{proof} Recall that $\Upsilon_{k,c} =   \mathrm{H}_{c}  \circ \pi_{k-1}$, the composition of the horizontalization map
$\mathrm{H}_{c}$ and the horizontal sequence of $k-1$ blow-ups $\pi_{k-1}$.
Assume that $\tilde{\F}= \pi_{k-1}^{*} \F$ is defined by the $1$-form $\tilde{\omega}$ and that
$\bar{\F} = \mathrm{H}_{c}^{*} \tilde{\F} = \Upsilon_{k,c}^{*} \F $ is defined by the $1$-form $\bar{\omega}$.
We have that the element $E_{k} \subset \Delta_{\omega}$ is taken by $\pi_{k-1}$ into the element $\tilde{E}_{1}$ of  $\Delta_{\tilde{\omega}}$
of co-slope $\kappa = 1$.
Suppose first that  $E_{k} \subset \Delta_{\omega}$ is an edge.
Then  $\tilde{E}_{1} \subset \Delta_{\tilde{\omega}}$ is also an edge.
Applying Lemma \ref{lem-Ndeg} $k-1$ times, we have that $\tilde{\omega}$ is Newton $\tilde{E}_{1}$-non-degenerate.
Then, using Proposition  \ref{prop-fact1}, we find that the element $\bar{E}_{1} \subset \Delta_{\bar{\omega}}$
of co-slope $\kappa =1$ has elevation $e(\bar{E}_{1}) = 1$, giving that $\Delta_{\bar{\omega}}$ is pre-primary.
On the other hand, if $E_{k}$ is a vertex, since it is associated with the separatrix $S$, it must be dicritical.
Thus, also $\tilde{E}_{1} \subset \Delta_{\tilde{\omega}}$ is a dicritical  vertex
(the multiplier decreases one unit per blow-up, so $\tilde{E}_{1}$ has multiplier  $\tilde{\lambda} = 1$)
associated with the separatrix
$\tilde{S} =  \pi_{k-1}^{*} S$.
Proposition \ref{prop-E1-vertex} gives that
$e(\bar{E}_{1})=1$ and, as before, we conclude  that
$\Delta_{\bar{\omega}}$ is pre-primary.
\end{proof}

\subsection{Ramifications}
\label{sec-ramification}
Suppose that $(x,y)$ are convenient coordinates for the holomorphic foliation $\F$,
defined by the $1$-form $\omega = A(x,y) \dd x + B(x,y) \dd y$.
Let us consider new local analytic coordinates $(\hat{x},\hat{y})$ and the map $\rho(\hat{x},\hat{y}) = (\hat{x}^{n},\hat{y}) = (x,y)$, which is a ramification of order $n > 1$ in
the $x$-variable. Let $S$ be an analytic branch at $(\C^{2},0)$, in the coordinates $(x,y)$, non-tangent to the vertical axis, of order $m$. If $n$ is a multiple
of $m$ then
$\rho^{*}S$ is an analytic curve, in the coordinates $(\hat{x},\hat{y})$, having $m$ smooth branches.
Actually,
given  an equation $f(x,y)$   of $S$, writing its factorization
\[ f(x,y) = \varphi (y - s_{1}(x)) \cdots (y - s_{m}(x)) ,\]
where $\varphi \in \C\{x,y\}$ is a unit  and $s_{1}, \ldots, s_{m} \in \C\{x^{\frac{1}{n}}\}$ are the conjugate Puiseux series of $S$,   we have
that
\[\rho^{*}f(\hat{x},\hat{y}) = (\rho^{*} \varphi) (\hat{y} - s_{1}(\hat{x}^{n})) \cdots (\hat{y} - s_{m}(\hat{x}^{n}))\]
 is an equation for $\rho^{*}S$
whose factors are equations of smooth branches.
We can choose the order of ramification $n$ as a multiple of the order of any separatrix of the foliation $\F$.
Indeed, the number of isolated separatrices is finite and so is the number of different equisingularity classes of dicritical separatrices,
 by the finiteness of the number
of irreducible components of $\D$ in a reduction of singularities (see
Section \ref{sec-holomorphic-foliations}).
Thus, the set of integers that are orders of separatrices of $\F$ is also finite.
We will need the following simple result:

\begin{lem}
\label{lem-remark-factors}
Let $F, G \in \C\{x,y\}$. They have a non-trivial common factor if and only if   $\rho^{*}F, \rho^{*}G \in \C\{\hat{x},\hat{y}\}$ have a non-trivial
common factor.
\end{lem}
\begin{proof}
The necessity in this statement is evident. For the sufficiency,
a non-trivial common factor $\phi \in \C\{\hat{x},\hat{y}\}$ for the pair $\rho^{*}F, \rho^{*}G$ would
engender common factors of the type $\psi_{i}(x^{\frac{1}{n}},y) = \phi(\zeta_{i} x^{\frac{1}{n}},y)$,
where $\zeta_{1}, \ldots, \zeta_{n}$ are the $n$-th roots of unity,  of the pair
$F,G$ in $\C\{x^{\frac{1}{n}},y\}$.
The  product $\psi$ of all different functions appearing in  $\{\psi_{1}, \cdots, \psi_{n}$\}   is a common factor of $F,G$ in $\C\{x^{\frac{1}{n}},y\}$ that is invariant under each conjugation
$\sigma_{\zeta_{i}}: \C\{x^{\frac{1}{n}},y\} \to \C\{x^{\frac{1}{n}},y\}$, defined by $f(x^{\frac{1}{n}},y) \mapsto f(\zeta_{i}x^{\frac{1}{n}},y)$.
Thus,  $\psi \in \C\{x,y\}$ and this is a non-trivial factor of $F,G$ in $\C\{x,y\}$.
\end{proof}

As before, let us suppose that $(x,y)$ are convenient coordinates for $\F$.
Let us consider $\hat{\omega} = \rho^{*} \omega$.
If $\omega = A(x,y) \dd x + B(x,y) \dd y$, then
\begin{equation}
\label{eq-1form-ramification}
\hat{\omega} = n A(\hat{x}^{n},\hat{y}) \hat{x}^{n-1}\dd \hat{x} + B(\hat{x}^{n},\hat{y}) \dd \hat{y}.
\end{equation}
As shown in Lemma \ref{lem-remark-factors} (applied to $F=A$, $G=B$,
which are coprime since $\F$ is a foliation),
 the coefficients of $\hat{\omega}$ are without common factors and, thus, $\hat{\omega}$
defines a foliation $\hat{\F}$, for which the coordinates $(\hat{x},\hat{y})$ are convenient.
The differential  binomial of $\omega$ of bidegree $(i,j)$,
$\omega_{ij} = a_{ij} x^{i-1}y^{j} \dd x +  b_{ij} x^{i}y^{j-1} \dd y $,
is taken to the differential binomial
\[ \hat{\omega}_{ni,j}
= n a_{ij} \hat{x}^{ni-1}\hat{y}^{j}   \dd \hat{x} +  b_{ij} \hat{x}^{ni}\hat{y}^{j-1} \dd \hat{y} , \]
 of bidegree $(ni,j)$.
Thus, the Newton polygon $\Delta_{\hat{\omega}}$ is obtained by a homothety of ratio $n$ in the horizontal direction.
In particular, edges in $\Delta_{\omega}$ and in  $\Delta_{\hat{\omega}}$ are  identified by a one-to-one correspondence:
the edge $E_{\kappa} \subset \Delta_{\omega}$, of co-slope $\kappa$, corresponds to the edge $\hat{E}_{n\kappa} \subset \Delta_{\hat{\omega}}$, of co-slope $n \kappa$.

The following result describes how Newton and Milnor numbers change under ramification:

\begin{prop}
\label{prop-numbers-ramification}
 In the above context, we have:
\begin{enumerate}[label=(\alph*)]
  \item $\nu_{0}(\hat{\F}) = n \nu_{0}(\F) + (n-1) (h(\Ne_{\omega})-1)$
  \item $\mu_{0}(\hat{\F}) = n \mu_{0}(\F) + (n-1) (h(\Ne_{\omega}) -1)$
\end{enumerate}
\end{prop}
\begin{proof}
For (a), using elementary geometry, we have $\cl{A}(\Ne_{\hat{\omega}}) = n  \cl{A}(\Ne_{\omega})$. Since
$\ell(\Ne_{\hat{\omega}}) = n  \ell(\Ne_{\omega})$, we have
\begin{eqnarray*}
  \nu_{0}(\hat{\F}) &=& 2\cl{A}(\Ne_{\hat{\omega}}) -  \ell(\Ne_{\hat{\omega}}) - h(\Ne_{\hat{\omega}}) + 1  \\
    &=& 2 n   \cl{A}(\Ne_{\omega}) -  n\ell(\Ne_{\omega}) - h(\Ne_{\omega}) + 1   \\
     &=&  n (2 \cl{A}(\Ne_{\omega}) -  \ell(\Ne_{\omega}) - h(\Ne_{\omega}) + 1) + (n-1) h(\Ne_{\omega}) - (n-1) \\
     &=&  n  \nu_{0}(\F) + (n-1)( h(\Ne_{\omega}) - 1).
\end{eqnarray*}

For (b), using \eqref{eq-1form-ramification}, we have
\begin{eqnarray*}
  \mu_{0}(\hat{\F}) &=&  (n A(\hat{x}^{n},\hat{y}) \hat{x}^{n-1} , B(\hat{x}^{n},\hat{y}))_{0}  \\
    &=& (A(\hat{x}^{n},\hat{y}), B(\hat{x}^{n},\hat{y}))_{0} + (n-1)( \hat{x}, B(\hat{x}^{n},\hat{y}))_{0}  \\
    &=&  n (A(x,y), B(x,y))_{0} + (n-1) \ma_{0}( B(0,y)) \\
    &=&  n \mu_{0}(\F) +  (n-1)( h(\Ne_{\omega}) - 1).
\end{eqnarray*}
In the third equality, we used Property (4) in \cite[p.146]{camacho1984}, knowing that the ramification map $\rho$
is $n$ to one if $x \neq 0$.
\end{proof}

A consequence of Proposition \ref{prop-numbers-ramification} and of Theorem \ref{teo-fol-nondeg2} is the fact
that, for a fixed system of convenient coordinates $(x,y)$, $\F$ is Newton degenerate if and only if the ramification
$\hat{\F}$ is Newton degenerate.
This can also be seen by the following proposition, which will be essential for several arguments to appear in the coming sections:

\begin{prop}
\label{prop-ramification-edge}
The $1$-form $\omega$ is Newton  $E_{\kappa}$-degenerate if and only if  $\hat{\omega}$ is Newton $\hat{E}_{n\kappa}$-degenerate.
As a consequence, $\omega$ is Newton degenerate if and only if  $\hat{\omega}$ is Newton degenerate.

\end{prop}
\begin{proof} Write $\hat{\omega} = \hat{A}(\hat{x},\hat{y}) \dd \hat{x} + \hat{B}(\hat{x},\hat{y}) \dd \hat{y}$.  We have
\[ (\hat{x} \hat{A})_{\hat{E}_{n\kappa}} = \sum_{(ni,j) \in \hat{E}_{n\kappa}} n a_{i,j} \hat{x}^{ni}\hat{y}^{j}
 = n   \rho^{*}(xA)_{E_{\kappa}} .\]
Similarly,
\[ (\hat{y} \hat{B})_{\hat{E}_{n\kappa}} =  \rho^{*}(yB)_{E_{\kappa}}. \]
The conclusion then follows from Lemma \ref{lem-remark-factors} above.
\end{proof}


\section{Weak separatrices and Newton degeneracy}
\label{sec-well-oriented-saddle-node}
Suppose that the germ of singular holomorphic foliation $\F$ has a reduction of singularities
$\pi:(M,\D) \to (\C^{2},0)$.
As already defined, we say that a branch $S \in \Sep(\F)$  is a \emph{weak separatrix} if the strict transform curve
$\tilde{S} = \pi^{*}S$ is a weak separatrix of a saddle-node singularity of the strict transform foliation $\tilde{\F} = \pi^{*} \F$.
The existence of a weak separatrix
is then equivalent to the existence of a well-oriented (i.e. non-tangent) saddle-node singularity.
In order to
identify $S$ as a weak separatrix, all we need is a reduction of singularities of $\F$ along $S$, i.e. a sequence $\pi$ of blow-ups with centers at base points of the successive strict transforms of  $S$, such that $\pi^{*}\F$ is reduced at the base point of $\pi^{*} S$.
Also noteworthy is the fact that the concept of well-oriented saddle-node, as well as that of weak separatrix, is stable under blow-ups.
Indeed, a blow-up at a   saddle-node singularity transforms its weak separatrix into a weak separatrix of a saddle-node
singularity, whose strong separatrix is contained in the exceptional divisor.
These ideas are implicit in the arguments of this section.

We start by proving that weak separatrices are also  stable under ramifications:

\begin{prop}
\label{prop-ramification-sn}
 Let  $\rho(\hat{x},\hat{y}) = (\hat{x}^{n},\hat{y})$ be a ramification of order $n > 1$. Let $\hat{\F} = \rho^{*} \F$ be the
germ of foliation defined by the $1$-form $\hat{\omega} = \rho^{*} \omega$. If $S$ is a weak separatrix of $\F$, then
all the irreducible components of
$\hat{S} = \rho^{*}S$ are weak separatrices of $\hat{\F}$.
\end{prop}
\begin{proof}
Let $\pi: (M, \D) \to (\C^{2},0)$ and $\hat{\pi}: (\hat{M}, \hat{\D}) \to (\C^{2},0)$ be  reductions of singularities
for $\F$ and $\hat{\F}$, respectively. The ramification $\rho$ can be lifted to  a rational map $\varrho: (\hat{M},\hat{\D}) \to (M,\D)$  such that $\rho \circ \hat{\pi} =   \pi \circ \varrho$. By  performing some additional blow-ups over
$(\hat{M},\hat{\D})$, we can resolve the indeterminacy points of $\varrho$ and still have a reduction of singularities   for
$\hat{\F}$. Thus, since non-tangent saddle-nodes and their weak separatrices are persistent under quadratic blow-ups,
we can assume that $\varrho$ itself is a morphism.
Fix a branch of $\hat{S} = \rho^{*}S$ and, to simplify notation, we keep denoting it by $\hat{S}$.
Let $\hat{p} \in \hat{\D}$ be the base point of $\hat{S}' = \hat{\pi}^{*}\hat{S}$ and set
$p = \varrho(\hat{p}) \in \D$. Then $p$ is the base point of $S' = \pi^{*}S$ and $\varrho(\hat{S}') = S'$.
Since $\pi$ and $\hat{\pi}$ are reductions of singularities, $p$ and $\hat{p}$ are either regular points or simple singularities
for $\F' = \pi^{*} \F$ and  $\hat{\F}'= \hat{\pi}^{*} \hat{\F}$.
They are trace points of $\D$ and $\hat{\D}$, where $S'$ and $\hat{S}'$ are smooth and transverse to
the exceptional divisor. Our assumption is that $\F'$ has a non-tangent saddle-node at $p$, so that $S'$ is its weak separatrix. Thus, if $p$ belongs to the component $D \in \D$, then
$D$ is  invariant by $\F'$ and contains the strong separatrix.
If $S'$ is purely formal, then so is  $\hat{S}'$ and it must be the weak separatrix of    a saddle-node singularity of $\hat{\F}'$
at $\hat{p}$. Otherwise, Dulac's normal form for $\F'$ converges. We can then find a system of local analytic coordinates $(u,v)$ at
$p$, in which $u=0$ and $v=0$ are, respectively, the local equations of $S'$ and $D$, so that
$\F'$ is defined by the meromorphic $1$-form
\[ \eta = \left( \frac{\zeta}{v} + \frac{1}{v^{k+1}}\right) \dd v  -  \frac{\dd u}{u} , \]
where $k \geq 1$ and $\zeta \in \C$.
We also choose at $\hat{p}$ coordinates $(\hat{u},\hat{v})$
in which $\hat{u}=0$ and $\hat{v}=0$ are, respectively, the local equations of $\hat{S}'$ and
$\hat{D}$, the component of
$\hat{\D}$ containing $\hat{p}$.
In these coordinates, write
$\varrho(\hat{u},\hat{v}) = (\hat{u}^{r_{1}} \hat{v} ^{s_{1}} \varrho_{1}, \hat{u} ^{r_{2}} \hat{v} ^{s_{2}} \varrho_{2})$, where
$\varrho_{1}, \varrho_{2}$ are units in $\C\{\hat{u},\hat{v}\}$. Since $\varrho(\hat{D}) \subset D$ and $\varrho(\hat{S}') = S'$, we must have
$r_{1}, s_{2} > 0$. Also, $r_{2} =0$ since, otherwise, $\hat{S}'=\{\hat{u}=0\}$ would be mapped
entirely into $D'=\{v=0\}$, contradicting the transversality of
$S'$ to $D'$. We have
\[ \varrho^{*} \eta  =   -r_{1} \frac{\dd \hat{u}}{\hat{u}} + (\zeta s_2 - s_1) \frac{\dd \hat{v}}{\hat{v}} + \zeta \frac{\dd \varrho_{2}}{\varrho_{2}} - \frac{\dd \varrho_{1}}{\varrho_{1}}
 + \frac{1}{ \hat{v}^{s_{2}k} \varrho_{2}^{k} } \left( s_{2} \frac{\dd \hat{v}}{\hat{v}} + \frac{\dd \varrho_{2}}{\varrho_{2}} \right) .\]
 Cancelling the poles, we see that $\hat{\F}'$ is defined by a holomorphic $1$-form of the kind $ \hat{u} \dd \hat{v} + \hat{v} \theta$,
where $\theta$ is a holomorphic $1$-form of order at least one.
Hence, at $\hat{p}$, the ramified foliation $\hat{\F}'$ is a saddle-node well-oriented with respect to the exceptional divisor $\hat{\D}$ of local equation $\hat{v} = 0$.
\end{proof}

\begin{prop}
\label{prop-non-resonant} Suppose that $S$ is a smooth weak separatrix of $\F$ and that $\Delta_{\omega}$ is primary.
Then $\Delta_{\omega}$ is non-resonant.
\end{prop}
\begin{proof}
Assume, as a contradiction hypothesis, that $\Delta_{\omega}$ is resonant.
  Let $E_{\kappa}$ be the only edge of $\Delta_{\omega}$ and let
$\lambda \in \Q_{>0}$ be the multiplier of $\Top(E_{\kappa})$. Let $y = c x^{k} + \cdots$, where $k \geq 1$ and $c \in \C^{*}$, be the power series of $S$.
We cannot have $k>\kappa$ since, in this
case, the element $E_k$ associated with $S$ would be the vertex $\Bot(E_\kappa)$ which, lying on the
horizontal axis, has multiplier  $\infty$ and cannot be dicritical.
On the other hand, if $k < \kappa$,   $E_{k}$ is a dicritical vertex associated with $S$ and, thus,
has multiplier    $\lambda=k$.
Applying the initial monomial elimination map $\Upsilon_{k,c}$, we take $S$ into a separatrix with horizontal tangent line
and $\Delta_{\omega}$ into a singular primary Newton polygon whose vertex $E_{1}$ is $1$-dicritical.
By Proposition \ref{prop-preprimary-blowup}, another blow-up produces a regular primary Newton polygon,
 which means that the foliation becomes regular and $S$ must correspond  to the leaf at the origin of the new coordinates.
That is, $S$ is a dicritical separatrix, which is not our case.

It remains to consider the case $k =\kappa$, in which the distinguished edge $E_{\kappa}$ is the element
associated with $S$.
First, suppose that the element  $\Top(E_{\kappa})$  has multiplier $\lambda \not \in \N$.
Each horizontal blow-up decreases the multiplier by one unit. Thus, if $\lambda < k - 1$, after a horizontal sequence of at most $k-1$ horizontal blow-ups, we obtain
multiplier in $\Q_{<0}$, corresponding to a simple non-degenerate singularity, which contradicts the hypothesis that $S$ is a weak separatrix.
On the other hand, if $\lambda > k-1$,
we apply the initial monomial elimination map $\Upsilon_{k,c}$. The new Newton polygon is primary and its distinguished edge has as multiplier $\lambda - (k-1)$,
since each blow-up decreases this multiplier in one unit, whereas the horizontalization procedure preserves it.
Thus,  after applying a finite number of initial monomial elimination maps,
 we will eventually obtain a singular primary Newton polygon whose   multiplier is in $\Q_{<0}$.
The singularity is then simple non-degenerate, leading again to a contradiction.
In order to finish the proof,  we have to analyze the situation in which
$E_{\kappa}$ is  the edge associated with $S$ and  $\Top(E_{\kappa})$  has multiplier $\lambda  \in \N$.
Again, if $k' = \lambda \leq \kappa$, a horizontal sequence of  $k'-1$   blow-ups
will produce a primary Newton polygon that is either Poincar\'e-Dulac, if $k' = \kappa$, or $1$-dicritical if $k' < \kappa$.  Both situations are incompatible with the fact of     $S$ being a weak separatrix.
Finally, if $\lambda > \kappa$, we fall in this latter case after the application of a finite number of initial monomial elimination maps,
starting with $\Upsilon_{k,c}$.
\end{proof}

The existence of weak separatrices implies  Newton degeneracy, as shown in the theorem below:
\begin{teoRM}
\label{th-weak} Let $\F$ be a germ of singular holomorphic foliation at $(\C^{2},0)$, defined by the $1$-form $\omega$. Suppose that
$\F$ admits a weak separatrix. Then, for any fixed system of convenient analytic coordinates $(x,y)$,
the Newton polygon $\Delta_{\omega}$ is Newton degenerate.
\end{teoRM}
\begin{proof}
Let $S$ be the weak separatrix we are assuming to exist.
We have convenient coordinates, thus $S$ is neither of the two coordinate axes.
Thus, after possibly switching the coordinates $x$ and $y$, we can suppose that $S$ is not tangent to the $y$-axis.
We can also assume that $S$ is smooth. Indeed, in the non-smooth case, we do a ramification having as order a multiple of
$\ma_{0}(S)$
and the theorem's conclusion follows from the smooth case by Propositions \ref{prop-ramification-edge} and \ref{prop-ramification-sn}.

We write the series of $S$, $y = c x^{k} + \cdots$, where $k\in \N$ and the dots stand for higher order terms.
 We apply the initial monomial elimination map $\Upsilon_{k,c}$, obtaining  new coordinates
$(\bar{x},\bar{y})$ and the $1$-form $\bar{\omega}$ defining $\bar{\F} = \Upsilon_{k,c}^{*} \F$.
We assert that the element $E_{k} \subset \Delta_{\omega}$ cannot be a vertex.  Indeed, if it were a vertex, it would be
dicritical, since it is associated with $S$. Hence, applying Proposition \ref{prop-E1-vertex}
 to $\tilde\omega=\pi_{k-1}^*\omega$
(i.e., just before applying ${\rm H}_c$),
we find that
$e(\bar{E}_{1}) = 1$,  so
that $\Delta_{\bar{\omega}}$ is pre-primary. Since, by the same proposition,  $\bar{E_{1}}$ is dicritical, another
blow-up leads us to a regular primary Newton polygon, so $S$ is a dicritical separatrix, a contradiction.

Thus, we can assume that $E_{k}$ is an edge of $\Delta_{\omega}$.
If  $e(\bar{E_{1}}) > 1$, then $E_{k}$ is Newton-degenerate by Proposition \ref{prop-fact1} and Lemma \ref{lem-Ndeg}.
Otherwise, $e(\bar{E_{1}}) = 1$ and  $\Delta_{\bar{\omega}}$ is pre-primary.
By Proposition \ref{prop-fact2}, we prove the theorem's assertion once we conclude that $\Delta_{\bar{\omega}}$ is Newton degenerate
with respect to its distinguished edge.
However, a new blow-up at the origin of the coordinates $(\bar{x},\bar{y})$ provides new coordinates $(\tilde{x},\tilde{y})$ and
a $1$-form $\tilde{\omega}$ whose Newton polygon $\Delta_{\tilde{\omega}}$ is
singular primary, which is non-resonant by Proposition \ref{prop-non-resonant}. Therefore, it corresponds to
a simple singularity of $\F$ and the
conclusion then follows from Proposition  \ref{prop-primary}.
\end{proof}

As seen in Theorem \ref{teo-fol-nondeg2},  Newton degeneracy can be characterized by a comparison of   Milnor and Newton numbers.
Both invariants, for a foliation $\F$ defined by a $1$-form $\omega$, can be calculated from the expression of  $\omega$,
by taking into account only a finite number of  terms of sufficiently low order.
Thus, as a consequence of Theorem \ref{th-weak}, these invariants can be used to guarantee that the foliation does
not admit purely formal separatrices:
\begin{cor} Suppose that the germ of singular holomorphic foliation $\F$ at $(\C^{2},0)$ is such that $\mu_{0}(\F) = \nu_{0}(\F)$
for some system of convenient analytic coordinates.
Then all separatrices of $\F$ are convergent.
\end{cor}
\begin{proof} Indeed, the equality of Milnor and Newton numbers assures that $\F$ is Newton non-degenerate. By Theorem \ref{th-weak}, $\F$ does not admit weak separatrices.
The conclusion then follows from the fact that all purely formal separatrices of a foliation are weak separatrices.
\end{proof}

We finish this section by giving a pair of concrete examples.  In the first one,   illustrating
Theorem \ref{th-weak},
the foliation has a tangent saddle-node, being thus Newton degenerate:

\begin{example}
\label{ex-newton-deg}
{\rm
Consider the foliation $\F$ defined by
\[ \omega = (-xy + y^2 -x^2y ) \dd x  +x^3\dd y .\]
We have that $0 \in \C^{2}$ is a non-dicritical singularity such that
$\ma_{0}(\F) = 2$. A straightforward calculation shows that $\mu_{0}(\F)=6$.

The  strict transform of $\tilde{\F} = \sigma^{*}\F$
is defined, in the coordinates $(x,t)$ in the first chart of $\sigma$,  by
\[ \tilde{\omega}  = t(t-1) \dd x +  x^2  \dd t . \]
Over the exceptional divisor $D$, the foliation $\tilde{\F}$ has three simple singularities. Two well oriented saddle-node singularities, one
 at $(x,t)=(0,0)$, whose weak separatrix is  the $x$-axis, the other at
$(x,t) = (0,1)$, whose week separatrix has equation $t -1=0$. There is also
a non-degenerate singularity at $q_3$ corresponding to $t = \infty$.

We make the linear change of coordinates $(x = 2\tilde{x} +   \tilde{y},y= \tilde{x} +  2\tilde{y})$ in order
to obtain convenient coordinates.
 Keeping the notation $(x,y)$ for the new coordinate system, we have that $\F$ is defined by the
 $1$-form  $\omega  = A(x,y) \dd x + B(x,y) \dd y $ given by
\[
\omega  = (-12x^2y  -2x^2 -12xy^2  -2xy -3y^3  +4y^2 ) \dd x
+ (12x^3+12x^2y  -x^2 +3xy^2  -xy +  2y^2 ) \dd y.\]
The Newton polygon  $\Delta_{\omega}$ is formed
by a single edge $E$ with vertices $(0,3)$ and $(3,0)$. Thus
the Newton number is $\nu_{0}(\F) = 4$, different from $\mu_{0}(\F)$.
The degeneracy equations are
\begin{eqnarray*}
  (xA)_{E} &=& x(-2x^2-2xy+4y^2) = 2x(x+2y)(y-x) \\
  (yB)_{E} &=&    y(-x^2-xy+2y^2) =  y(x+2y)(y-x) ,
\end{eqnarray*}
having the equations (of the tangent cones)  of the weak separatrices as common factors.
}\end{example}

In our second example, the foliation is Newton non-degenerate and has a single tangent saddle-node in its reduction of singularities, hence no weak separatrix. This shows that
the implication stated    in Theorem \ref{th-weak}, far from being a more general property related to the  existence of saddle-node
singularities,    is a phenomenon specifically associated with
the existence of weak separatrices.

\begin{example}
\label{ex-newton-ndeg}
{\rm
The foliation $\F$ defined by
\[ \omega = (xy + y^2 -x^2y ) \dd x  +(-x^2 + x^3)\dd y  \]
has a non-dicritical singularity at $0 \in \C^{2}$ such that
$\nu_{0}(\F) = 2$. An easy calculation shows that $\mu_{0}(\F) = 4$.
We have that $\tilde{\F}= \sigma^{*}\F$   is
defined in  the coordinates $(x,t)$ of the blow-up chart  by
\[ \tilde{\omega} = t^2 \dd x +  x(x-1)  \dd t . \]
The foliation $\tilde{\F}$ has two simple singularities over the exceptional divisor $D$: one tangent saddle-node
 at $(x,t)=(0,0)$, whose strong separatrix is  the $x$-axis,   and
a non-degenerate singularity corresponding to $t = \infty$.
Making  the linear change of coordinates
$(x,y)  = (2\tilde{x} +   \tilde{y},\tilde{x} +  2\tilde{y})$ in order to obtain convenient coordinates, and
keeping the notation $(x,y)$ for the new coordinate system,
we obtain that $\F$ is defined by the $1$-form
(retaking the foliation of Example \ref{ex-fol-nondeg})
\[\omega = (-12x^2y+ 2x^2-12xy^2+14xy-3y^3+11y^2) \dd x +
 (12x^3+12x^2y-5x^2+3xy^2+xy+4y^2) \dd y .\]
As shown there, $\nu_0(\F)=\mu_0(\F)=4$, so $\F$ is Newton
non-degenerate in these new coordinates.
}\end{example}

In Examples \ref{ex-newton-deg} and \ref{ex-newton-ndeg}, both foliations   have   reductions of singularities formed by a single blow-up ---
hence, by a toric map --- and admit saddle-node singularities. Newton non-degeneracy fails only in the first one, where
weak separatrices are present.


\section{Dicritical families of separatrices}
\label{sec-dicritical-family}

We begin by establishing the following concept: for two distinct power series or fractional power series    $\sigma, \sigma'$,   their \emph{contact} is the integer $\cont(\sigma, \sigma') =  \ma_{0}(\sigma - \sigma' )$. We extend  this definition
by setting $\cont(\sigma, \sigma) = \infty$.

A  germ of   analytic branch $S$ at $(\C^{2},0)$ has  order $n = \ma_{0}(S)$, taken as $\ma_{0}(f)$,
where $f \in \C\{x,y\}$ is one of its irreducible equations.
From the Newton-Puiseux theory (see for instance \cite{casas2000}), if $S$ has non-vertical tangent, then there are
$n$ different fractional power series $s \in \C\{x^{\frac{1}{n}}\}$ such that $f(x,s(x)) = 0$. These are called
\emph{Puiseux series}. Recall that, given $s(x) = \sum_{i \geq 0} a_{i} x^{\frac{i}{n}}$ a fractional power series in $\C\{x^{\frac{1}{n}}\}$,
its conjugation by  $\zeta \in \C^{*}$ is the series 
$\sum_{i \geq 0} a_{i} \zeta^{i} x^{\frac{i}{n}}$.
The $n$ Puiseux series of the branch $S$ are   obtained one from the other through the conjugation by
  $n$-th roots of unity. We define the \emph{Puiseux order} of $S$ as being  the order of any of its Puiseux series, and
denote it by $\oo(S)$. We remark that the description  of this paragraph, with the necessary adaptations, also applies to formal branches.

To a singular analytic  branch $S$ with non-vertical tangent, we can associate a \emph{genus} $g \geq 1$
and $g$ positive rational numbers $\varepsilon_{1}, \ldots, \varepsilon_{g}$, its  \emph{characteristic exponents},
satisfying
\begin{enumerate}
  \item $1 <  \varepsilon_{1} < \cdots < \varepsilon_{g}$;
  \item $\varepsilon_{i} \in \Q \setminus \Z$, for $i= 1,...,g$;
  \item there are $m_{1},\ldots, m_{g}, n_{1},\ldots, n_{g} \in \Z_{>1}$,
  with $m_i,n_i$ relatively prime, such that
$\varepsilon_{i} = m_{i}/n_{1} \cdots n_{i}$ for $i=1,\ldots,g$.
\end{enumerate}

Note that the integer $n= n_{1} \cdots n_{g}$ is $\ma_{0}(S)$. Thus, we will write $\varepsilon_{g} = m /n$, where $m = m_{g}$.
The integers $m,n$ need not  be relatively prime, but they will be so when $g=1$, for in this case $m= m_{1}$ and $n = n_{1}$.
For future reference, we make the convention that $n_{0}=1$. The genus and the set of characteristic exponents are topological invariants
and can be   calculated in an algorithmic way from the exponents of any Puiseux series of $S$, provided $S$ has non-vertical tangent
\cite{casas2000}. In this situation, they are independent of the coordinates.

Given two branches $S,S'$ at $(\C^{2},0)$ with non-vertical tangent, their \emph{contact} is defined
as
\[ \cont(S,S') = \max_{s,s'} \cont(s,s') ,\]
where $s$ and $s'$ run over the conjugated Puiseux series of $S$ and $S'$, respectively.
If $S$ is singular and has horizontal tangent,   we  say that it has \emph{maximal horizontal contact} if  $\oo(S) = \varepsilon_{1}$,
where $\varepsilon_{1}$ is the first characteristic exponent of $S$. This is equivalent to saying that the Puiseux
 series of $S$ have the form $y = c x^{\varepsilon_{1}} + \cdots$, where $c \in \C^{*}$. Explaining the terminology, among
 all singular branches
with non-vertical tangent and having $\varepsilon_{1}$ as the first characteristic exponent
(which, thus, admit  Puiseux expansion in $x$ and actually   have
horizontal tangent direction),
the ones with Puiseux order $\varepsilon_{1}$
 are precisely those having maximal contact with the $x$-axis.

We   make the following definition:

\begin{ddef}
\label{def-dicritical-family-regular}
{\rm
 A \emph{dicritical family of branches} of  \emph{genus} $0$ and \emph{index}
$k \geq 1$, with $k \in \N$,
is a family $\mathrm{F} = \{S_{\alpha}; \, \alpha \in A\}$ of smooth branches with non-vertical tangent
such that:
\begin{enumerate}
  \item If $S_{\alpha}, S_{\alpha'} \in \mathrm{F}$ are distinct branches, then
  \[\cont(S_{\alpha}, S_{\alpha'}) =  k .\]
  \item  If $S$ is a branch satisfying $\cont(S_{\alpha}, S) =  k$ for some  $S_{\alpha} \in \mathrm{F}$, then there exists $S_{\alpha'}  \in \mathrm{F}$
 such that $\cont(S_{\alpha'}, S) > k$.
  \end{enumerate}
}\end{ddef}
Note that item $(2)$ is a condition of maximality, in the sense that if $\mathrm{F}$ and $\mathrm{F}'$ are    families of smooth branches
such that $\mathrm{F} \subset \mathrm{F}'$,
with $\mathrm{F}$ satisfying both conditions of  Definition \ref{def-dicritical-family-regular} and   $\mathrm{F}'$  satisfying (1),  then  $\mathrm{F} = \mathrm{F}'$.
In concrete terms, $\mathrm{F}$ can be described by taking $A = \C$ and then writing the power series of $S_{\alpha}$ as
\[  y = p(x) + \alpha x^{k} + \cdots, \]
where $p(x)$   is a polynomial of degree less than $k$ if non-zero, common to all branches in $\mathrm{F}$,
 the dots
representing terms of higher order.
Seen in another way, if $S \in  \mathrm{F}$  and $\pi = \sigma_{1} \circ \cdots \circ \sigma_{k}$ is a   sequence of $k$ blow-ups with
centers at the base points of the strict transforms of  $S$ --- which are trace points since $S$ is smooth --- then,
for each trace point $q \in D \subset \D$, where $D$ and $\D$ are the exceptional divisors of $\sigma_{k}$ and of $\pi$,
there is a unique branch $S_{\alpha} \in  \mathrm{F}$ such that $q = \pi^{*}S_{\alpha} \cap D$.

We now establish an analogous definition for  singular curves.
We say that an ordered set of positive rational numbers $\E = \{\varepsilon_{1} , \ldots, \varepsilon_{g}\}$, where $g \geq 1$, is a \emph{set of characteristic exponents} if they are characteristic exponents of some singular analytic branch with non-vertical tangent. The \emph{order} of $\E$ is the algebraic multiplicity common to every branch having $\E$ as its set of characteristic exponents.

\begin{ddef}
\label{def-dicritical-family-singular}
{\rm Let $\E = \{\varepsilon_{1} , \ldots,  \varepsilon_{g}\} $  be a set of characteristic exponents
and let $k \in \Z$, with $k \geq 0$. Write $\varepsilon_{g} = m/n$, where $n$ is the order of $\E$.
A \emph{dicritical family    of branches} of  \emph{genus $g$}, \emph{exponents}
$\E$ and \emph{index} $k$ is a family $\mathrm{F} = \{S_{\alpha}; \, \alpha \in A\}$
of singular analytic branches having genus $g$ and characteristic exponents $\E$ such that:
\begin{enumerate}
  \item If $S_{\alpha}, S_{\alpha'} \in \mathrm{F}$ are distinct branches, then
  \[\cont(S_{\alpha}, S_{\alpha'}) =  \frac{m + k}{n} =  \varepsilon_{g} +  \frac{ k}{n}.\]
    \item  If $S$ is a branch, whose first $g$ characteristic exponents are the ones in $\E$,
       satisfying $\cont(S_{\alpha}, S) =  \varepsilon_{g} +  \frac{ k}{n}$ for some  $S_{\alpha} \in \mathrm{F}$, then there exists $S_{\alpha'}  \in \mathrm{F}$
 such that $\cont(S_{\alpha'}, S) > \varepsilon_{g} +  \frac{ k}{n}$.
\end{enumerate}
}\end{ddef}
Putting in other words,
setting $\varepsilon' = \varepsilon_{g} +  \frac{ k}{n}$, for each  $c \in \C$, if $k > 0$, or $c \in \C^{*}$, if $k =0$, there is a unique $S_{\alpha} \in \mathrm{F}$ having
a Puiseux series of the form $y = s(x) + c x^{\varepsilon'} + \cdots$, where $s(x)$ comprises
monomials in $x^{\frac{1}{n'}}$  of order less than $\varepsilon'$,
where $n' = n = n_{1} \cdots n_{g}$, if $k>0$,   $n' =  n_{1} \cdots n_{g-1}$, if $k=0$
(with the convention that
   $n'  =1$ if  $g =1$). Furthermore,  given
Puiseux series $y = s(x) + c x^{\varepsilon'} + \cdots$ of  $S_{\alpha}$
and $y = s'(x) + c' x^{\varepsilon'} + \cdots$ of $S_{\alpha'}$, where $S_{\alpha},S_{\alpha'} \in \mathrm{F}$,
then there is an
$n$-th root of   unity $\zeta$ such that $s'(x)$ is the conjugate of $s(x)$ by $\zeta$. Evidently, from (1), for this $\zeta$
we  must always have
$c' \neq \zeta^{n \varepsilon'} c$ if $S_{\alpha} \neq S_{\alpha'}$. We observe that all branches in $\mathrm{F}$ are equisingular.
Note that, if $\pi$ is a minimal sequence of blow-ups that  desingularizes one --- and hence every --- branch in $\mathrm{F}$, then
all branches in the family $\tilde{\mathrm{F}} = \{ \pi^{*}S_{\alpha}; S_{\alpha} \in \mathrm{F}\}$, when $k>0$, have the same base point $p \in D$, where $D$ is the exceptional divisor of the last blow-up of $\pi$, and thus, at $p$, they form a dicritical family of branches of genus $g=0$ and index $k$. On the other hand,
when $k=0$, different branches in $\tilde{\mathrm{F}}$ have different base points on $D$.

If $\mathrm{F}$ is a  dicritical family    of branches of genus $g>0$  with index $k=0$, as in Definition \ref{def-dicritical-family-singular}, we say that a branch $S_{\mathrm{sr}}$ is a \emph{semi-root} of $\mathrm{F}$  if its order is
  $\ma_{0}(S_{\mathrm{sr}}) =  n_{1} \cdots n_{g-1}$ and
$\cont(S_{\mathrm{sr}},S_{\alpha}) = \varepsilon_{g}$ for every $S_{\alpha} \in \mathrm{F}$.
This means, in particular, that $S_{\mathrm{sr}}$ has genus $g-1$ and has
 $\E' = \{\varepsilon_{1} , \ldots,  \varepsilon_{g-1}\}$  as its set of characteristic exponents.
In the literature, such an $S_{\mathrm{sr}}$ is said to be a $(g-1)$-semi-root  of $S_{\alpha}$, where $S_{\alpha}$
is any of the branches in the family $\mathrm{F}$ (see \cite{popescu1999}).

We call \emph{dicritical family of branches}  a set $\mathrm{F}$ of analytic branches at $(\C^{2},0)$
satisfying either Definition \ref{def-dicritical-family-regular} or Definition \ref{def-dicritical-family-singular}.
This object then comes associated with a genus $g \geq 0$, with an index $k  \geq 0$ --- which is strictly positive if $g=0$ --- and,
when $g>0$, with a set of characteristic exponents $\E$. Dicritical families of branches are of special interest
when their elements are separatrices of a foliation. In this case, we make the following definition:

 \begin{ddef}
 {\rm
\label{def-dicritical-family-separatrices}
Let $\F$ be a germ of holomorphic foliation at $(\C^{2},0)$.
A dicritical family of branches $\mathrm{F}$  is said to be a
\emph{dicritical    family  of separatrices} for $\F$ if $\mathrm{F} \subset \sep(\F)$ and
 the following additional condition is satisfied:
\begin{enumerate}
\setcounter{enumi}{2}
  \item If $\pi: (M,\D) \to (\C^{2},0)$ is the minimal sequence of blow-ups that separates any pair of branches in $\mathrm{F}$, then, setting $p_{\alpha} = \pi^{*}S_{\alpha} \cap \D$, the strict transform foliation $\pi^{*} \F$ is regular at $p_{\alpha}$ for every $\alpha \in A$.
\end{enumerate}
}\end{ddef}
Condition (3) says that,   if $D \subset \D$ is the component introduced by the last blow-up of $\pi$, then
the branches in $\mathrm{F}$ are  dicritical separatrices associated with $D$.
This condition is necessary, since we could have the situation where $D$ is fibered by smooth  $\F$-invariant curves
and contains a singularity, for instance, given by the Poincar\'e-Dulac normal form (see \S\ref{sec-blow-ups}).
If (3) is satisfied, the branches in $\mathrm{F}$ are all the separatrices associated with $D$, except when $\F$
is the radial foliation --- i.e. its reduction of singularities
consists of a single dicritical blow-up ---  for which there is
also  a dicritical separatrix with vertical tangent
(recall that in Definition \ref{def-dicritical-family-regular} all branches have non-vertical tangent).
In particular, this condition says that  if $S_{\alpha} \in \mathrm{F}$, then there is no $S \in \Sep(\F)$ such that
$\cont(S_{\alpha},S)  > \cont(S_{\alpha},S_{\alpha'})$, where $S_{\alpha'}$ is any  branch in $\mathrm{F}$
different from $S_{\alpha}$.
We say that an $\F$-invariant   branch $S_{\mathrm{sr}}$ is a \emph{semi-root} of the dicritical family of separatrices
 $\mathrm{F}$, of genus $g \geq 1$ and index $k=0$, if it is    a semi-root of the underlying dicritical family of branches.

The following technical lemma concerns the use of ramification maps in handling dicritical families of separatrices:
\begin{lem}
\label{lem-proj-dicritical-family}
Let $S \in \sep(\F)$ be a branch having non-vertical tangent,   genus $g=1$, order $\ma_{0}(S) =n$  and characteristic exponent $\varepsilon = m/n$. Let $\rho(\hat{x},\hat{y}) = (\hat{x}^{n},\hat{y})$ be the ramification map of order $n$.
 Suppose that some branch of $\rho^{*} S$ lies in a dicritical
family of separatrices $\hat{\mathrm{F}}$
of $\hat{\F} = \rho^{*} \F$, of genus $0$ and index $m+k$, for some $k \geq 0$.
Then $S$ belongs to a dicritical family of separatrices $\mathrm{F}$ of $\F$, of genus $g=1$ and index $k$.
Besides, if $k=0$, there is a unique branch in
$\sep(\F)$ having the property of being a semi-root of  $\mathrm{F}$.
\end{lem}
\begin{proof}
We first observe that  a branch at $(\C^{2},0)$ is in $\sep(\F)$ if and only if one --- and hence every --- branch of its pre-image
by $\rho$ is in $\sep(\hat{\F})$.
We also know that, if $\hat{S}$ is one of the branches of $\rho^{*}S$, then it has a power series
$\hat{y} = s(\hat{x}^{n})$, for some
Puiseux series $y = s(x)$   of $S$. Since $g=1$, we can write $s(x) = q(x) + c_{0}x^{\varepsilon} + \cdots$, where
$c_{0} \in \C^{*}$, $q \in \C[x]$ is of degree less than $\varepsilon$ and
the dots stand for a  series in $\C\{x^{\frac{1}{n}}\}$ of order greater than $\varepsilon$. Thus, $\hat{S}$ is defined
by a power series of the kind $s(\hat{x}^{n}) = q(\hat{x}^{n}) + c_{0}\hat{x}^{m} + \cdots$.

Suppose that $\hat{S}$ belongs   to the dicritical family of separatrices  $\hat{\mathrm{F}}$ as in the lemma's statement.
The elements of $\hat{\mathrm{F}}$ are the branches with power series of the form
$\hat{s}_{c}(\hat{x}) = \hat{p}(\hat{x}) + c \hat{x}^{m+k} + \cdots$,
where $c$ varies over $\C$ and $\hat{p}(\hat{x})$ is a polynomial of degree less than $m+k$.
In particular, since $\hat{S}$ belongs to $\hat{\mathrm{F}}$, the initial terms of its series
$ s(\hat{x}^{n})$ are given by $\hat{p}(\hat{x})$.
Thus,
\begin{enumerate}[label=(\roman*)]
\item  $\hat{p}(\hat{x}) = q(\hat{x}^{n})$ if $k=0$;
\item   or $\hat{p}(\hat{x}) = q(\hat{x}^{n}) + c_{0}\hat{x}^{m} +  r(\hat{x})$,
where $r \in \C[\hat{x}]$ has terms of order between $m = n \varepsilon$ and $m+k= n\varepsilon' $, if $k>0$.
\end{enumerate}
Denote by $\hat{S}_{c}$ the element of $\hat{\mathrm{F}}$ with power series $\hat{s}_{c}(\hat{x})$.  The curve $S_{c} = \rho(\hat{S}_{c})$
has a Puiseux series of the form
$y = s_{c}(x) = \hat{s}_{c}(x^{\frac{1}{n}}) =   p(x) + c x^{\varepsilon'} + \cdots$,
where $p(x) = \hat{p}(x^{\frac{1}{n}})$ and $\varepsilon' = \frac{m+k}{n} = \varepsilon + \frac{k}{n}$. Thus:
\begin{enumerate}[label=(\roman*)]
  \item if $k=0$, $p(x) = q(x)$ is a polynomial of degree less than $\varepsilon$.
  \item if $k > 0$, we can write $p(x) = q(x) + c_{0} x^{\varepsilon} + r(x^{\frac{1}{n}})$, where $c_{0} \in \C^{*}$,
$q$ and $r$ are polynomials, $q$ of degree less than $\varepsilon$ and $r$ comprising terms of   degrees between $m = n \varepsilon$
and $m+ k =n\varepsilon'$.
\end{enumerate}
Define  $\mathrm{F} = \{S_{c};   c \in \C^{*}\}$, if $k=0$, and $\mathrm{F} = \{S_{c};   c \in \C\}$, if $k \neq 0$.
Taking into account the description of $p(x)$ above, in both cases, the branches in $\mathrm{F}$ have genus $g=1$ and
characteristic exponent $\varepsilon$.
We claim that $\mathrm{F}$, defined in this way, is a dicritical family of separatrices, which evidently
has genus $g=1$ and index $k$. In order to prove this, we begin by  making an assertion:

\begin{assertion} Let $\zeta$ be a non-trivial  $n$-th root of unity.
Then, for every $c \in \C^{*}$,
\begin{enumerate}[label=(\roman*)]
  \item  $S_{c} = S_{c\zeta}$ if  $k=0$;
  \item  $S_{c} \neq S_{c\zeta}$ if  $k \neq 0$.
\end{enumerate}
\end{assertion}

\begin{proof}
Let us  prove (i).
Recalling that   $(m,n)=1$, we first take $r,s \in \Z$ such that $1=rm + sn$.
We have to show that $s_{c}(x)$ and $s_{c \zeta}(x)$
are Puiseux series of the same curve.
If  this    were false, then,  by taking  the conjugation
of $s_{c \zeta}(x)$ by $\zeta^{-r}$, we would find a Puiseux series whose first term with fractional exponent is
$ c \zeta (\zeta^{-r})^{m}x^{\varepsilon} = c \zeta^{sn} x^{\varepsilon} = c  x^{\varepsilon}$.
Hence, we would have two different branches in $\sep(\F)$ having Puiseux series
with initial part $q(x) + c x^{\varepsilon}$, which would be lifted through $\rho$ into two different branches in $\sep(\hat{\F})$
having power series with initial part  $q(\hat{x}^{n}) + c \hat{x}^{m}$. We also know that  one of these branches, $\hat{S}_{c}$, is in $\hat{\mathrm{F}}$. This violates condition (3) in the definition of dicritical family of separatrices (see comment in the paragraph following Definition
\ref{def-dicritical-family-separatrices}).

For the proof of (ii), consider the Puiseux series defining $S_{c}$ and $S_{c\zeta}$,
\[ s_{c}(x) =       q(x) + c_{0} x^{\varepsilon} + r(x^{\frac{1}{n}}) + c x^{\varepsilon'} + \cdots
\quad \text{and} \quad s_{c\zeta}(x) =       q(x) + c_{0} x^{\varepsilon} + r(x^{\frac{1}{n}}) + c\zeta x^{\varepsilon'} + \cdots. \]
Suppose, by contradiction, that $S_{c} = S_{c\zeta}$. If this happens, the Puiseux series $s_{c \zeta}(x)$ must be the conjugate
of $s_{c}(x)$   by some non-trivial $n$-th root of unity $\zeta'$.
We then have  that the conjugate  of $s_{c}(x)$ by $\zeta'$
has  $c_{0}(\zeta')^{m}$ as coefficient of  $x^{\varepsilon}$,
whereas the corresponding coefficient of $s_{c\zeta}(x)$ is $c_{0}$. This gives   a contradiction, since $(m,n)=1$ and, then $(\zeta')^{m} \neq 1$.
\end{proof}

In the case $k=0$, let us set $A = \C^{*}/\sim$, where $``\sim"$    is the equivalence relation
defined in $\C^{*}$ by the product by some $n$-th root of unity. In the case $k \neq 0$, we simply take  $A = \C$.
Next, we define the family of branches
$\mathrm{F} = \{ S_{\alpha}; \alpha \in A\}$, where, when $k=0$, $S_{\alpha}$ is the branch  having
Puiseux series $y = s_{c}(x)$, for some
representative $c \in \C^{*}$ of $\alpha \in A$.
From this construction, it follows easily that, if $\alpha, \alpha' \in A$ are distinct, then
$\cont(S_{\alpha}, S_{\alpha'}) =  \varepsilon'$ and we have
 condition (1) in the definition of dicritical family of branches.
The maximality condition (2) is also evident. Thus, $\mathrm{F}$ is a dicritical family of branches.
Note that, in the case $k \neq 0$, $S_{0}$, the branch corresponding to $c=0$, lies in the same equisingularity class of the curves $S_{\alpha}$, for
$\alpha \in \C^{*}$, and it belongs to $\mathrm{F}$. On the other hand, if $k=0$, even though $\cont(S_{\alpha},S_{0}) = \varepsilon' = \varepsilon$ for
every $\alpha \in A$,
$S_{0}$ does not have $\varepsilon$ as a characteristic exponent and, hence, it cannot belong to the same dicritical family of
branches as the other curves of $\mathrm{F}$. In this case,   we can make the following assertion:

\begin{assertion} For $k=0$, $S_{0}$ is a semi-root of the dicritical family of branches $\mathrm{F}$.
Besides, it is the only separatrix of $\F$ having  this  property.
\end{assertion}
\begin{proof}
We have that $S_{0}$ has a  Puiseux series of the form
$y=s_{0}(x) = p(x) +  z(x^{\frac{1}{n}})$, were $z$ is a power series with terms of
order greater than $n \varepsilon$. In principle, $z(x^{\frac{1}{n}})$ lies in  $\C\{x^{\frac{1}{n}}\}$.
However, we actually have $z(x^{\frac{1}{n}}) \in \C\{x\}$ for, otherwise, $s_{0}(x)$ would be conjugated, by   $n$-th roots of unity, to different
Puiseux series of $S_{0}$, which in turn would be lifted through $\rho$ into different separatrices of $\hat{\F}$ having contact larger than $\varepsilon$
with $\hat{S}_{0}$, the separatrix of $\hat{\mathrm{F}}$ with power series $\hat{s}_{0}(\hat{x})$.
This is impossible, since $\hat{\mathrm{F}}$ satisfies   condition (3) in Definition \ref{def-dicritical-family-separatrices}.
As a conclusion, we get that $S_{0}$ is smooth and, hence, it is a  semi-root of $\mathrm{F}$.
Concerning  uniqueness, it suffices to see that any other separatrix having the property of being
a semi-root of $\mathrm{F}$ is smooth and has contact greater than $\varepsilon$ with $S_{0}$. Its lifting by
$\rho$ would be a smooth separatrix of $\hat{\F}$ with contact greater than $n\varepsilon$ with $\hat{S}_{0}$. This is  impossible,
again by    (3) in Definition \ref{def-dicritical-family-separatrices}.
\end{proof}

We already know that $\mathrm{F}$ is a dicritical family of branches, all these branches in $\sep(\F)$.
In order to  prove that  it is indeed a dicritical family of separatrices, we have to verify, for $\mathrm{F}$, condition (3) in Definition \ref{def-dicritical-family-separatrices}.
For this purpose,
we consider two sequences of blow-ups $\pi: (M, \D) \to (\C^{2},0)$ and $\hat{\pi}: (\hat{M}, \hat{\D}) \to (\C^{2},0)$,
the ones appearing  in condition (3), the first for the dicritical family of branches $\mathrm{F}$,
the second   for the dicritical family of separatrices $\hat{\mathrm{F}}$.
Denote by $D \subset \D$ the component touched by $\pi^{*}S $ and set $p = D \cap \pi^{*} S$.
Similarly,  $\hat{D} \subset \hat{\D}$ is the component touched by $\hat{\pi}^{*} \hat{S}$
and   $\hat{p} = \hat{D} \cap \hat{\pi}^{*} \hat{S}$.
The components $D$ and $\hat{D}$ are   dicritical  for $\F$ and $\hat{\F}$ and, considering condition
(3)  for $\hat{\mathrm{F}}$,   we have that
$\hat{\pi}^{*}\hat{\F}$ is regular at $\hat{p}$.
We want to prove that $\pi^{*} \F$ is regular at $p$.
First, we remark that if $S' \in \sep(\F)$  is such that  $p = D \cap \pi^{*}S'$, then $S' = S$. Indeed,
in this case
there would be a component $\hat{S}'$ of $\rho^{*}S'$
such that $\hat{p} = \hat{D} \cap \hat{\pi}^{*} \hat{S}'$.
 However, since
 $\hat{\pi}^{*}\hat{\F}$ is regular at $\hat{p}$, we must have
$\hat{S}' =  \hat{S}$, implying that $S' = S$.

We then have that  $\pi^{*}S$ is the only separatrix of $\pi^{*}\F$ at $p$, which is also smooth  and transverse to $D$.
Under these conditions,   $\pi^{*}\F$ is singular at $p$
only if there exists, at $p$, a   transcendent leaf of $\pi^{*}\F$. Indeed, if all local leaves
are closed (i.e. contained in analytic curves) then, by \cite{mattei-moussu1980}, $\pi^{*}\F$ has a holomorphic first integral at $p$, so it
is a generalized curve type foliation \cite{camacho1984} and must be regular at $p$ \cite[Lem. 4]{mol2002}.
Let us assume that  $L$ is a transcendent leaf of $\pi^{*}\F$ accumulating at $p$.
We lift  the ramification map  $\rho$  to  a rational map $\varrho: (\hat{M},\hat{\D}) \to (M,\D)$  such that $\rho \circ \hat{\pi} =   \pi \circ \varrho$.
Note that $\hat{p}$ is not a  point of indeterminacy of $\varrho$. Indeed, if this were the case, $\varrho(\hat{p})$ would be a union of irreducible components
of $\D$ and, since $\varrho(\hat{p})$ contains the trace point $p \in D$, it would contain the whole $D$. Hence, the pre-images by $\varrho$ of the separatrices of $\pi^{*} \F$ touching the dicritical component $D$ would all contain $\hat{p}$.
This is absurd, since we know that $\hat{\pi}^{*} \hat{\F}$ is regular at $\hat{p}$.
The image by $\varrho$ of a small neighborhood of $\hat{p}$ in $\hat{M}$, chosen in such a way that
the leaves of $\hat{\pi}^{*} \hat{\F}$ are local regular curves transverse to $\hat{D}$,  is a neighborhood of $p$ in $M$.
Hence,  $\varrho^{-1}(L)$ is a non-empty invariant set for $\hat{\pi}^{*} \hat{\F}$ accumulating at $\hat{p}$.
Thus, locally at $\hat{p}$, $\varrho^{-1}(L) \subset \hat{S}$, which gives $L \subset S$ and, so, $L$  cannot be transcendent.
As a conclusion, $\pi^{*}\F$ is regular at $p$. The same argument works for any branch in $\mathrm{F}$, finishing
the proof that $\mathrm{F}$ is a dicritical family of separatrices.
\end{proof}

Lemma \ref{lem-proj-dicritical-family} could have been formulated in a more general way, involving dicritical
families of higher genus. We opted to present   only the cases we work with in this article.
For instance, with a similar proof, we can state the following:

\begin{lem}
\label{lem-proj-dicritical-family-g2}
Let $S \in \sep(\F)$ be of genus $g=2$,  with non-vertical tangent, having first characteristic exponent   $\varepsilon_{1} = m_{1}/n_{1}$, where $(m_{1},n_{1})=1$. Let $\rho(\hat{x},\hat{y}) = (\hat{x}^{n_{1}},\hat{y})$ be the ramification map of order $n_{1}$.
 Suppose that some branch of $\rho^{*} S$ lies in a dicritical
family of separatrices $\hat{\mathrm{F}}$ of $\hat{\F} = \rho^{*} \F$, of genus $g=1$ and index $k=0$.
Then $\mathrm{F} = \rho(\hat{\mathrm{F}}) =  \{ \rho(\hat{S});  \hat{S} \in \hat{\mathrm{F}} \}$ is a dicritical
family of separatrices of $\F$, of genus $g=2$ and index $k=0$, containing $S$.
\end{lem}

Pre-primary Newton polygons can be used  to identify dicritical families of separatrices.
This is the subject of the next three lemmas.

\begin{lem}
\label{lem-el-dic-smooth1}
 Suppose that the  Newton polygon  $\Delta_{\omega}$ is  pre-primary and that the element
 $E_{1}$ is dicritical.
Then all separatrices with   non-vertical  tangent  are in a dicritical family of genus $g =0$ and index $k=1$.
\end{lem}
\begin{proof}
Since $E_{1}$  is dicritical, so is the first blow-up $\sigma$ at $0 \in \C^{2}$.
Let $\tilde{\omega} = \sigma^{\hash} \omega$ be the $1$-form defining   $\tilde{\F} = \sigma^{*} \F$.
By Proposition \ref{prop-preprimary-blowup},   the Newton polygon $\Delta_{\tilde{\omega}}$  is regular primary.
Hence, $\tilde{\F}$ is regular at the origin of the blow-up coordinates in the first chart. Indeed, by the comment following
the same proposition, $\tilde{\F}$ is transverse to all points of the $\tilde{y}$-axis in the blow-up coordinates,
which correspond to the exceptional divisor $D = \sigma^{-1}(0)$ minus the point
associated with the vertical direction. These local leaves are projected by $\sigma$ onto separatrices which,
following our terminology,  lie in a dicritical family of genus $g=0$ and index $k=1$.
\end{proof}

\begin{lem}
\label{lem-el-dic-smooth}
 Suppose that the  Newton polygon  $\Delta_{\omega}$ is  pre-primary
with      distinguished edge $E_{\kappa}$ and multiplier $\lambda = k \in \N$ with $k>1$.
Assume that there  exists $S \in \sep(\F)$ with  horizontal  tangent.
Then all separatrices with horizontal tangent are in a dicritical family of genus $g =0$ and index $k$.
\end{lem}
\begin{proof} We start by considering the particular case $k \leq \kappa$.
We cannot have  $k = \kappa$, for in this case $\Delta_{\omega}$ would be Poincar\'e-Dulac and,
by Lemma \ref{lem-poincare-dulac},
 $\F$ would not admit
a separatrix with horizontal tangent.
Hence $k < \kappa$.

In this case, we apply a horizontal sequence of $k-1$ blow ups. By Proposition \ref{prop-preprimary-blowup}, the Newton polygon
obtained  is primary and, thus, its element of co-slope $1$ corresponds to $\Top(E_{\kappa})$, having  multiplier
$k - (k-1) =1$. We have then the hypothesis of Lemma
\ref{lem-el-dic-smooth1}.
 The dicritical family of separatrices of this lemma corresponds, by blowing-down, to
the one in the above statement.

Let us complete the proof, supposing that $k > \kappa$.
Under this assumption, $S$ must be an $E_{\kappa}$-separatrix. Indeed, its Puiseux series cannot start with a term in $x^{\ell}$, with $\ell < \kappa$,
for otherwise $\Delta_{\omega}$ would be $\ell$-dicritical. This   is not the case, since $\ell < k$.
Hence, a Puiseux series of $S$
has the form $y = c x^{\kappa} + \cdots $, where $c  \neq 0$ is the unique root of the principal
polynomial $\Phi_{E_{\kappa}}$, since $\Phi_{E_{\kappa}}$ has degree one.
Then, the initial monomial elimination map $\Upsilon_{\kappa,c}$
produces a primary Newton polygon whose edge has co-slope $\kappa'>1$ and whose multiplier
is $k' = k - (\kappa -1)$. If $k' \leq \kappa'$, we are done by the particular case of the first paragraph. Otherwise, we repeat the process,
applying new monomial elimination maps, until we fall in the
particular case.
\end{proof}

The next result is  a version of    Lemma \ref{lem-el-dic-smooth}   for singular separatrices:

\begin{lem}
\label{lem-singular-nondeg}
 Suppose that   the  Newton polygon  $\Delta_{\omega}$ of $\F$ is  pre-primary with distinguished edge $E_{\kappa}$.
If $S \in \sep(\F)$ is singular with
horizontal tangent, then it is
a dicritical separatrix which,  along with all singular separatrices  with horizontal tangent,
lies in a dicritical family of genus $g=1$ and index $k=0$. Besides, there is a unique smooth separatrix with horizontal
tangent, which is a semi-root of this  family.
\end{lem}
\begin{proof} Let  $\varepsilon_{1} > 1$ denote the first characteristic exponent of $S$. Write $\varepsilon_{1} = m_{1}/n_{1}$,
where $(m_{1},n_{1}) = 1$. We first prove the following:

\begin{assertion*} The multiplier of $\Delta_{\omega}$ is $\varepsilon_{1}$.
\end{assertion*}
\begin{proof}
Assume first that $\varepsilon_{1} < \kappa$.
We have that $S$ has a Puiseux series of the form
$y = q(x) + c x^{\varepsilon_{1}} + \cdots$, where $c \in \C^{*}$ and $q(x)$ is a polynomial, of degree less than $\varepsilon_{1}$ if non-zero.
We claim that $q(x) \equiv 0$. Indeed, if $q(x)$ were non-zero, setting  $\ell >1$ to be its order, then, since $\ell < \varepsilon_{1} < \kappa$,
$S$ cannot be an $E_{\kappa}$-separatrix.  Thus $\Delta_{\omega}$ would be $\ell$-dicritical, where $\ell \in \N$ is its multiplier.
By  Lemma \ref{lem-el-dic-smooth},
$S$ would be smooth, which contradicts our hypothesis. Hence, the Puiseux series of $S$ has the form
$y =  c x^{\varepsilon_{1}} + \cdots$. Again, $S$ cannot be an  $E_{\kappa}$-separatrix and we have that $\Delta_{\omega}$ is $\varepsilon_{1}$-dicritical, thus, its multiplier is $\varepsilon_{1}$.
If $\varepsilon_{1} > \kappa$ we observe that   horizontal blow-ups and horizontalization procedures always produce   pre-primary Newton polygons. A horizontal blow-up decreases the multiplier by one unit, whereas a horizontalization procedure preserves the multiplier.
The Assertion's conclusion is then a consequence of the initial assumption,
after applying a convenient number  initial monomial elimination maps, possibly followed by a horizontal sequence of blow-ups.
\end{proof}

Let us proceed with the lemma's proof.
After a  change of coordinates highly tangent to the identity, we can assume that our coordinates are convenient.
Considering that $\varepsilon_{1} = m_{1}/n_{1}$, we   apply a ramification map $\rho$ with   order of ramification $n_{1}$.
The Newton polygon
$\Delta_{\hat{\omega}}$, where $\hat{\omega} = \rho^{*} \omega$ defines $\hat{\F} = \rho^{*} \F$, is not, in principle, pre-primary:
its last edge has height one and elevation zero (co-slope $n_{1} \kappa$ and   multiplier $m_{1} = n_{1} \varepsilon_{1} \in \N$), but it may contain other edges of co-slope greater than one. However,
 after applying
a horizontal sequence of $n_{1}-1$ blow-ups, we get a foliation $\hat{\F}'$  defined by a $1$-form $\hat{\omega}'$ having,
in coordinates $(\hat{x}',\hat{y}')$,  a pre-primary Newton polygon $\Delta_{\hat{\omega}'}$ with
multiplier $m_{1} - n_{1}  + 1 > 1$.
The above Puiseux series of  $S$ generates $\hat{S}' \in \sep(\hat{\F}')$ of Puiseux order $m_{1} - n_{1}  + 1$.
In particular, $\hat{S}'$ has horizontal tangent with respect  to the coordinates $(\hat{x}',\hat{y}')$.
Thus, by Lemma \ref{lem-el-dic-smooth}, all separatrices of $\hat{\F}'$ with horizontal tangent lie in a dicritical family of
genus $g=0$ and index $m_{1} - n_{1}  + 1 \in \N$.
Thus, the corresponding  separatrices of $\hat{\F}$  lie in a dicritical family of
genus $g=0$ and index $m_{1} \in \N$.
The foliation $\hat{\F}$  may have other separatrices,  associated with the elements of $\Delta_{\hat{\omega}}$ whose co-slope  less than $n_{1} \kappa$. These separatrices, however, project through $\rho$ onto separatrices of $\F$ having non-horizontal tangent lines
and can be disregarded.
We have, in particular, that the branches of $\rho^{*}S$ are smooth and, so $S$ has genus $g=1$.
The proof is then finished by the application of Lemma \ref{lem-proj-dicritical-family}.
\end{proof}

Concerning pre-primary Newton polygons and isolated separatrices, we can state the following:

\begin{lem}
\label{lem-isolated-pre-primary}
Suppose that $\Delta_{\omega}$ is pre-primary. Then $\F$ has at most one isolated separatrix with horizontal tangent.
\end{lem}
\begin{proof}
Suppose, by contradiction, that there are two different isolated separatrices of $\F$, say $S_{1}$ and $S_{2}$, with horizontal tangent.
By Lemma  \ref{lem-singular-nondeg}, both $S_{1}$ and $S_{2}$ are smooth.
Also, denoting  by $E_{\kappa}$  the distinguished edge of $\Delta_{\omega}$,  both separatrices
must be $E_{\kappa}$-separatrices.
Indeed, if one of these separatrices, say $S_{1}$, were not associated with $E_{\kappa}$, its associated element would be
 $E_{k}$, with $k \in \N$ and $k>1$ and,  by    Lemma \ref{lem-el-dic-smooth},
$S_{1}$ would be dicritical, giving a contradiction.
We then write their
power series in the form  $y = c_{1} x^{\kappa} + \cdots$ and $y = c_{2} x^{\kappa} + \cdots$, where $c_{1},c_{2} \in \C^{*}$
are roots of the principal polynomial $\Phi_{E_{\kappa}}$.
Since $\Delta_{\omega}$ is pre-primary, $\Phi_{E_{\kappa}}$ has degree one and, thus, $c_{1} = c_{2} = c$.
Next,  we apply the initial monomial elimination $\Upsilon_{\kappa,c}$,
finding a foliation with primary (and hence pre-primary)
Newton polygon, having two isolated separatrices, $\tilde{S}_{1} = \Upsilon_{\kappa,c}^{*}S_{1}$ and $\tilde{S}_{2} = \Upsilon_{\kappa,c}^{*}S_{2}$,
with horizontal tangent. The above argument
gives that, in the new coordinates $(\bar{x},\bar{y})$, the initial terms of   their power series coincide.
Hence,   the second term of the power series of $S_{1}$ and $S_{2}$ in the coordinates $(x,y)$ also coincide.
Proceeding by induction,   after the application of $n >1 $ initial monomial elimination maps, we find   a foliation
with primary Newton polygon having two isolated separatrices with horizontal tangent (the transforms of $S_{1}$ and $S_{2}$)
 whose power series have coinciding initial terms in
the new coordinates. Then,   the $n+1$ first terms of the power series of $S_{1}$ and $S_{2}$ are the same.
This inductive process allows us to conclude that the separatrices $S_{1}$ and $S_{2}$ have the same power series and, hence, coincide.
\end{proof}


\section{Newton sets of separatrices}
\label{sec-newton-germs}

In this section, we characterize the set of separatrices of a Newton non-degenerate foliation $\F$.
This is done through a series of propositions and corollaries, which are later systematized in Definition
\ref{def-Newtontype-separatrices} and Theorem \ref{teo-newton-curve1}.
Let us fix  a   system of convenient  coordinates $(x,y) \in \C^{2}$
 denote by $\sep_{c}(\F)$ the family
of separatrices of $\F$ whose tangent line has slope   $c \in \C$.
Under this notation, the branches in $\sep_{0}(\F)$ are separatrices with horizontal tangent line.
We will denote by $\sep_{\infty}(\F)$ the family of separatrices having vertical tangent line.
The first two propositions  examine separatrices with oblique tangent line.
In the case of  smooth separatrices, we can state the following:

\begin{prop}
\label{prop-nondeg-smooth}
Suppose that $\Delta_{\omega}$ is Newton non-degenerate.
Let $S \in \sep_{c}(\F)$ be a smooth separatrix, where $c \in \C^{*}$ is the slope of its tangent line.
Suppose that all separatrices in $\sep_{c}(\F)$ are smooth.
Then either $S$ is dicritical and the separatrices in $\sep_{c}(\F)$ form a dicritical family of genus $g=0$,
or $S$ is isolated and there are no elements in $\sep_{c}(\F)$ other than $S$.
\end{prop}
\begin{proof}
Apply the horizontalization  map $\mathrm{H}_{c}$ and set $\bar{\omega} = \mathrm{H}_{c}^{*} \omega$.
This map takes $\Delta_{\omega}$ into the Newton polygon $\Delta_{\bar{\omega}}$,
which is pre-primary, by Proposition \ref{prop-fact1}, if $E_{1}$ is an edge, or by
Proposition \ref{prop-E1-vertex}, if $E_{1}$ is a vertex.
The family $\sep_{c}(\F)$ is taken into the family  $\sep_{0}(\bar{\F})$ of separatrices of $\bar{\F} = \mathrm{H}_{c}^{*} \F$ having horizontal tangent.  In particular
 $\bar{S} = \mathrm{H}_{c}^{*}S$  has the same classification of $S$ as dicritical or isolated.
If $\bar{S}$ is a dicritical separatrix, $\Delta_{\bar{\omega}}$ is resonant with  multiplier $\bar{\lambda} \in \N$.
Indeed, if $\bar{\lambda} \neq 1$, recalling that $\bar{S}$ is smooth,
we apply a finite number of initial monomial elimination maps of the form $\Upsilon_{k,c}$, possibly followed by
a sequence of horizontal blow-ups,   up to
the moment that precedes a dicritical blow-up. At this point we know that the Newton polygon is $E_{1}$-dicritical and, hence,
the multiplicity of the distinguished edge is $\lambda = 1$. The conclusion on $\bar{\lambda}$ then follows by going
back to $\Delta_{\bar{\omega}}$.
We can then apply  Lemmas \ref{lem-el-dic-smooth1} and  \ref{lem-el-dic-smooth}
to conclude that the branches in $\sep_{0}(\bar{\F})$ are in a dicritical family of genus $g=0$.
This gives the result in the dicritical case and also  implies  that, if $\bar{S}$ is an isolated separatrix, then $\sep_{0}(\bar{\F})$ does not contain dicritical separatrices.    We
can then apply  Lemma \ref{lem-isolated-pre-primary} in order to obtain the conclusion in the non-dicritical case.
\end{proof}

We have the following statement concerning singular separatrices with oblique tangents:
\begin{prop}
\label{prop-nondeg-oblique}
Suppose that $\Delta_{\omega}$ is Newton non-degenerate. Let $S \in \sep_{c}(\F)$ be a singular separatrix, where $c \in \C^{*}$ is the slope of its tangent line. Then $S$ lies in a dicritical family of separatrices of genus $g=1$ and index $k=0$. Besides, the separatrices in $\sep_{c}(\F)$ are the
other elements of this dicritical family along with a single smooth separatrix, which is a semi-root of this family.
\end{prop}
\begin{proof}
The horizontalization of coordinates by the map $\mathrm{H}_{c}$ takes $\Delta_{\omega}$ into a pre-primary Newton polygon $\Delta_{\bar{\omega}}$.
Again, this follows from  Proposition   \ref{prop-fact1}, if $E_{1}$ is an edge, or from
Proposition  \ref{prop-E1-vertex}, if  $E_{1}$ is a vertex.
We then apply  Lemma \ref{lem-singular-nondeg}  to the separatrix of horizontal tangent
$\bar{S} = \mathrm{H}_{c}^{*}S$ of $\bar{\F} = \mathrm{H}_{c}^{*} \F$
in order to obtain the proposition's conclusion.
\end{proof}

The above propositions have the  following corollaries:
\begin{cor}
\label{cor-nondeg-smooth}
Suppose that $\omega$ is Newton non-degenerate and
that there are two smooth separatrices $S, S' \in \sep(\F)$ with non-vertical tangent lines
such that
$$\cont(S,S') > \min\{\mathrm{o}(S),\mathrm{o}(S')\}.$$
Then $S$ and $S'$ belong to the same dicritical family of separatrices.
\end{cor}
\begin{proof} The hypothesis says that $S$ and $S'$ are
defined by power  series of the form $y =  c  x^{\ell} + \cdots$, where $\ell = \min\{\mathrm{o}(S),\mathrm{o}(S')\} \in \N$ and $c \in \C^{*}$ is the same
for both curves.
If $\ell=1$ then, by Proposition \ref{prop-nondeg-oblique}, since  $\sep_{c}(\F)$ already contains two smooth separatrices,   all
its separatrices are smooth.
Thus, Proposition \ref{prop-nondeg-smooth} gives the result. If $\ell>1$, we apply a horizontal sequence of $\ell-1$ blow-ups, using Lemma \ref{lem-Ndeg} in order to assure that non-degeneracy is preserved, getting two separatrices with the same oblique tangent line. The result then follows from the particular case $\ell = 1$.
\end{proof}

\begin{cor}
\label{cor-nondeg-oblique1}
Suppose that $\omega$ is Newton non-degenerate. Let  $S  \in \sep(\F)$ be a singular  separatrix with non-vertical
tangent line without the property of
maximal horizontal contact. Then $S$ lies in a dicritical family of separatrices of genus $g=1$ and index $k=0$.
Besides, the separatrices $S'  \in \sep(\F)$ satisfying $\cont(S,S') >  \mathrm{o}(S)$ are the ones
in this dicritical family along with the only smooth separatrix  which is a semi-root of this family.
\end{cor}
\begin{proof} We have that $S$ has Puiseux series of the form $y =  c  x^{\ell} + \cdots$, where $c \in \C^{*}$ and $\ell = \mathrm{o}(S) \in \N$ is such that
$\ell < \varepsilon_{1}$, where $\varepsilon_{1}$ is the first characteristic exponent of $S$.
Besides, the hypothesis on the contact says that the Puiseux series of the separatrices $S'$ all have the same initial term $c x^{\ell}$.
The corollary then follows from  Proposition \ref{prop-nondeg-oblique}, preceded by the application of a horizontal sequence of $\ell -1$ blow-ups if $\ell >1$.
\end{proof}

The following proposition deals with separatrices of genus $g \geq 2$:

\begin{prop}
\label{prop-nondeg-g2}
Suppose that $\omega$ is Newton non-degenerate. Let  $S \in \sep(\F)$ be a singular separatrix with non-vertical tangent such that $g(S) \geq 2$. Then $S$ has horizontal tangent and has maximal horizontal  contact, $g(S) = 2$ and $S$ belongs to a dicritical family of separatrices of genus $g=2$ and index $k=0$. Also, there is
a unique  separatrix of genus $g=1$ that is a semi-root of this family.
Furthermore, if $S' \in \sep(\F)$
is such that $\cont(S,S') > \mathrm{o}(S)$, then either $g(S') = 2$ and $S'$ belongs to the   same dicritical family of separatrices of $S$ or $g(S') = 1$ and $S'$ is the unique   semi-root of this dicritical family.
\end{prop}
\begin{proof}
From Proposition \ref{prop-nondeg-oblique}, we already know that $S$ must have horizontal tangent.
Let $\varepsilon_{1} = m_{1}/n_{1}$ be its first characteristic exponent, written as a reduced fraction, and
$y = p(x) + c x^{\varepsilon_{1}} + \cdots$ be one of its Puiseux series, where $ c \in \C^{*}$,  $p(x)$ is a polynomial  of degree
 less than $\varepsilon_1$, if
non-zero,  the dots comprising terms of order greater then $\varepsilon_1$.
  If $p(x) \not\equiv 0$, then, after a horizontal sequence of $\ell-1$ blow-ups,
where $\ell = \ma_{0}(p(x))$, we reach the situation in which the strict transform of $S$ is a separatrix   with oblique tangent of genus $g \geq 2$.
Since, by Lemma \ref{lem-Ndeg}, the Newton polygon at this step is non-degenerate, we have a contradiction with Proposition \ref{prop-nondeg-oblique}.
Thus, we must have $p(x) \equiv 0$, which means that $S$ has maximal horizontal contact.

Now, supposing that the initial
coordinates $(x,y)$ are convenient, which may be achieved through a change of coordinates  highly tangent to the identity, we   take the ramification map $\rho(\hat{x},\hat{y}) = (\hat{x}^{n_1},\hat{y})$.
The separatrix $S$ is lifted in  $n_1$ separatrices of $\hat{\F} = \rho^{*} \F$. Fix one of these separatrices, say $\hat{S}$.
Its genus is  $g(\hat{S}) = g(S) - 1$ and
its Puiseux series are of the form $\hat{y} =  c \hat{x}^{m_{1}} + \cdots$.
By Corollary  \ref{cor-nondeg-oblique1}, we have that $g(\hat{S}) = 1$ and that $\hat{S}$ lies in a dicritical family of separatrices of genus $g=1$ and index $k=0$. Besides, for this
family, there exists a unique semi-root  $\hat{S}_{0}$, which is  smooth.
After Lemma \ref{lem-proj-dicritical-family-g2},
we can conclude that $g(S) = 2$ and that $S$ lies in a dicritical family of separatrices of genus $g=2$ and index $k=0$.  The  unique semi-root of this dicritical family of separatrices is $S_{0} = \rho(\hat{S}_{0})$, which has genus $g=1$.

Finally, suppose that inequality in the proposition's statement is verified. Then, following the steps of the previous paragraph,
 $\rho^{*} S'$ contains a branch $\hat{S}'$ whose
Puiseux series also have  $c \hat{x}^{m_{1}}$ as initial term.
Again, by  Corollary  \ref{cor-nondeg-oblique1}, either $\hat{S}'$  belongs to the same dicritical family of separatrices as $\hat{S}$ or it
is the  semi-root of this family.
 Descending this information to $S'$, we get the proposition's conclusion.
\end{proof}

The next result concerns singular isolated separatrices:
\begin{prop}
\label{prop-nondeg-sing-isolated}
 Suppose that $\omega$ is Newton non-degenerate. Let  $S \in \sep(\F)$ be a singular isolated separatrix with non-vertical tangent.    Then  $g(S) = 1$, $S$ has horizontal tangent and maximal horizontal contact.
\end{prop}
\begin{proof} Since $\Delta_{\omega}$ is Newton non-degenerate,   Proposition \ref{prop-nondeg-g2} gives that   $g(S) \leq 2$.
However, by the same proposition, we cannot  have $g(S)=2$  since, in this case, $S$ would be a dicritical separatrix (isolated and dicritical
separatrices are mutually exclusive   by definition).
Thus $g(S) = 1$. Also, since $S$ is isolated, it follows from Proposition  \ref{prop-nondeg-oblique} that $S$ must have horizontal tangent.
On the other hand, if we had $\ell = \mathrm{o}(S) < \varepsilon$, where  $\varepsilon$ denotes the characteristic exponent of $S$, then,
after $\ell-1$ blow-ups, we would find a foliation
with non-degenerate Newton polygon for which the strict transform of $S$ is a singular isolated separatrix with
oblique tangent. This would contradict Proposition  \ref{prop-nondeg-oblique}. Thus $\mathrm{o}(S) = \varepsilon$ and $S$ has maximal horizontal contact, as wished.
\end{proof}

 Newton non-degeneracy imposes a condition of minimal contact between  separatrices. This is true, for instance, for
 pairs of isolated separatrices:

 \begin{prop}
\label{prop-isolated-contact}
 Suppose that $\omega$ is Newton non-degenerate and that    $S_{1}, S_{2} \in \sep(\F)$ are  distinct    isolated separatrices with non-vertical tangent.
Then  we have $\cont(S_{1},S_{2}) = \min \{ \oo(S_{1}) , \oo(S_{2}) \}$.

\end{prop}
\begin{proof}
Let us first assume that both $S_{1}$ and $S_{2}$ are smooth. Then, the equality in the statement fails to be true only
if $\ell = \oo(S_{1}) = \oo(S_{2})$ and the power series of both curves have the same initial term $c  x^{\ell}$, where $c \in \C^{*}$.
If $\ell=1$, this is not possible by Propositions \ref{prop-nondeg-smooth} and
\ref{prop-nondeg-oblique}. On the other hand, if $\ell >1$,
 a horizontal sequence of $\ell-1$ blow-ups   takes $S_{1}$ and $S_{2}$ into a pair of smooth isolated separatrices
of a foliation having the  Newton non-degeneracy property, both sharing the same oblique tangent line of slope $c$. We have thus fallen in the case $\ell=1$, which has just been eliminated.

Now, let us assume that one of the separatrices, say $S_{1}$, is singular.
Then, by Proposition \ref{prop-nondeg-sing-isolated}, $S_{1}$ has genus $g=1$, it has horizontal tangent and
maximal horizontal contact.  That is,   $\oo(S_{1}) = \varepsilon$, where $\varepsilon$ is its only characteristic exponent.
Thus, the condition on the contact holds unless $S_{2}$ is also singular, having $\varepsilon$ as its only characteristic exponent
and, besides, $S_{1}$ and  $S_{2}$ have Puiseux series with the same initial term $c  x^{\varepsilon}$, where $c \in \C^{*}$.
Write $\varepsilon = m/n$, with $m  \in \N$ and $n = \mathrm{m}_{0}(S_{1}) = \mathrm{m}_{0}(S_{2})$.
We  assume that the initial coordinates $(x,y)$ are convenient and
consider the ramification map $\rho(\hat{x},\hat{y}) = (\hat{x}^{n},\hat{y})$.
Take a pair of   branches   $\hat{S}_{1}$ and  $\hat{S}_{2}$ of, respectively,  $\rho^{*}S_{1}$ and  $\rho^{*}S_{2}$,
having power series of the form $\hat{y} = c   \hat{x}^{m} + \cdots$.
By Lemma \ref{lem-proj-dicritical-family}, neither of them belongs to a dicritical family of separatrices,
since $S_{1}$ and $S_{2}$ are both isolated. Thus  $\hat{S}_{1}$ and  $\hat{S}_{2}$ are  isolated separatrices of $\hat{\F} = \rho^{*}\F$
such that $\cont(\hat{S}_{1},\hat{S}_{2}) > m$.
Taking into account that $\Delta_{\hat{\omega}}$ is Newton non-degenerate by Proposition \ref{prop-ramification-edge},
this is impossible
by the smooth case  above.
This concludes the proof of the proposition.
\end{proof}

Our description of separatrices of a Newton non-degenerate foliation is completed with the following result:

\begin{prop}
\label{prop-flat-g1}
Suppose that $\omega$ is Newton non-degenerate. Let  $S \in \sep(\F)$ be a dicritical separatrix of genus $g=1$,  having
 maximal horizontal contact. Then $S$ lies in a dicritical family of separatrices of genus $g=1$ and index $k$.
Furthermore, if $\tilde{S} \in \sep(\F)$, with $\tilde{S} \neq S$,
is such that $\cont(S,\tilde{S}) >  \mathrm{o}(S)$, then  $k>0$ and $\tilde{S}$ is singular, with  $g(\tilde{S}) = 1$,  belonging to the   same dicritical family of separatrices as $S$.
\end{prop}
\begin{proof}
Denote by $\varepsilon$ the characteristic exponent of $S$.
 Let $D \subset \D$ be the dicritical component associated with $S$, where
$\D$ is the exceptional divisor of the reduction of singularities of $\F$. Let us first
assume that if $S'$ is any other dicritical separatrix associated with $D$, then
$\cont(S,S') = \varepsilon$.
Writing $\varepsilon = m/n$, where $m,n$ are relatively prime, we do the ramification
$\rho(\hat{x},\hat{y}) = (\hat{x}^{n},\hat{y})$, producing the foliation $\hat{\F} = \rho^{*} \F$
defined by the $1$-form $\hat{\omega} = \rho^{*} \omega$.
If $y = cx^{\varepsilon }+ \cdots$ is a Puiseux series of $S$, then $\rho^{*}S$ has a smooth branch $\hat{S}$
given by the power series $\hat{y} = c\hat{x}^{m}+ \cdots$, which is a separatrix for $\hat{\F}$.
Since $\Delta_{\omega}$ is Newton non-degenerate,  $\Delta_{\hat{\omega}}$ is Newton non-degenerate, by
Proposition \ref{prop-ramification-edge}. Besides,  the
element $\hat{E}_{m} \subset \Delta_{\hat{\omega}}$ of co-slope $m$ is dicritical. Indeed, taking other dicritical separatrices of
$\F$ associated with the component $D$, we  find that $\hat{\F}$ has  separatrices of Puiseux parametrization
$\hat{y} = \hat{c}\hat{x}^{m}+ \cdots$ for infinitely many values of $\hat{c} \in \C^{*}$.
Also, by Proposition \ref{prop-dicritical-edge}, $\hat{E}_{m}$ is a vertex.
Since $m>1$, we apply  a horizontal sequence of $m-1$ blow-ups,   obtaining a foliation $\hat{\F}'$, whose Newton polygon
$\Delta_{\hat{\omega}'}$ has the element $\hat{E}_{1}'$, corresponding to $\hat{E}_{m}$, as a dicritical  vertex of co-slope one.
The branch $\hat{S}$, in particular, is taken to the branch $\hat{S}'$, which has oblique tangent of slope $c$.

Let us first consider the case   $e(\hat{E}_{1}') = 1$. In this situation, $\Delta_{\hat{\omega}'}$ is pre-primary and  Lemmas
\ref{lem-el-dic-smooth1} gives that the  separatrices of $\hat{\F}'$ with
non-vertical tangent form a dicritical family of genus $g=0$ and index $k=1$. Hence, the separatrices
of $\hat{\F}$ associated with $\hat{E}_{m}$ form a dicritical family of genus $0$ and index $m$.  The conclusion
that $S$ belongs to a dicritical family of separatrices of genus $g=1$ and index $k=0$ then comes from
Lemma \ref{lem-proj-dicritical-family}.

On the other hand, if $e(\hat{E}_{1}') > 1$,
 the horizontalization map $\mathrm{H}_{c}$, where $c$ is the slope of the tangent line of $\hat{S}'$,  takes
this branch into one with horizontal tangent and
$\Delta_{\hat{\omega}'}$ into a pre-primary Newton polygon, which has a dicritical edge of co-slope $\kappa =1$.
This follows from  Proposition \ref{prop-E1-vertex}.
 We are then in the case of a pre-primary Newton polygon and, again, we find that
 $S$ belongs to a dicritical family of separatrices
of genus $g=1$ and index $k=0$.

Let us now suppose that $\cont(S,S') = \varepsilon' > \varepsilon$
for some --- and hence for every --- dicritical separatrix $S'$  associated with the component $D$, with $S' \neq S$.
Hence, these separatrices associated with  $D$   have
Puiseux series of the form  $y = cx^{\varepsilon}+ ... + c' x^{\varepsilon'} + \cdots$, where $\varepsilon' = \varepsilon + \frac{k}{n}$ for some $k \geq 1$, being the coefficient $c \in \C^{*}$  fixed while
$c' \in \C$   depends on the dicritical separatrix. In particular, for different dicritical separatrices
$S,S'$, their corresponding  coefficients $c'$ are different.
We proceed as in the previous paragraph:
through the ramification $\rho(\hat{x},\hat{y}) = (\hat{x}^{n},\hat{y})$, the branch $S$, which has a Puiseux series as above,
 is lifted by $\rho$
 in $n$ smooth branches, one of them, say $\hat{S}$, having power series
$\hat{y} = c\hat{x}^{m}+ \cdots + c' \hat{x}^{m + k} + \cdots$.
Starting with $\Upsilon_{m,c}$, we apply a finite number of initial monomial elimination maps   in order to get rid of the terms of order lower than
$m+k$. The same argument carried out above gives that
$\hat{S}$ belongs to a dicritical family of separatrices of genus $g=0$ and index $m+k$.
By Lemma \ref{lem-proj-dicritical-family} again, we have that $S$ lies in a dicritical family of separatrices
of genus $g=1$ and index $k$.

Finally, as for the last part of the proposition's statement, if $\tilde{S} \in \sep(\F)$ is such that $\cont(S,\tilde{S}) >  \mathrm{o}(S) = \varepsilon$,
then, in all operations above, the strict transforms of $S$ and $\tilde{S}$ will follow together until they both become separatrices
with non-vertical tangents associated with a pre-primary Newton polygon. By  Lemma \ref{lem-el-dic-smooth},
they are both in the same dicritical family of separatrices of genus $g=0$ and index $m + k$.
Hence, by Lemma  \ref{lem-proj-dicritical-family}, either $\tilde{S}$ is a dicritical separatrix   in the same dicritical family of  $S$,
implying that $\tilde{S}$ is singular of genus $g=1$,
or $\tilde{S}$  is the only  semi-root of this dicritical family, giving that $\tilde{S}$ is smooth and $k=0$.
However, the semi-root case does not happen, since we would have   $\cont(S,\tilde{S}) =  \mathrm{o}(S)$.
We then conclude that $k>0$ and that the separatrices $S$ and $\tilde{S}$ belong to the same dicritical families.
\end{proof}

Our objective in this section is to describe the set of separatrices of a Newton non-degenerate foliation.  With this goal in mind, we establish
Definitions \ref{def-Newtontype-horizontal},  \ref{def-Newtontype-oblique} and \ref{def-Newtontype-separatrices}
below. In them we assume
that a system of convenient analytic $(x,y)$ is fixed.

\begin{ddef}
\label{def-Newtontype-horizontal}
{\rm
 We say that $\sep_{0}(\F)$, the set of separatrices of $\F$ having horizontal
tangent line, is of  \emph{Newton type}  if, when non-empty, is
composed of:
\begin{enumerate}
\item finitely many smooth isolated separatrices;
\item  finitely many isolated separatrices of genus $g=1$ with maximal horizontal contact;
\item dicritical families of separatrices of genus $g=0$;
\item dicritical families of separatrices of genus $g=1$ and index $k=0$  without maximal  horizontal contact
(having a unique associated semi-root);
\item dicritical families of separatrices of genus $g=1$   having maximal horizontal contact;
\item dicritical families of separatrices of genus $g=2$ and index $k=0$    having maximal horizontal contact
(having a unique associated semi-root).
\end{enumerate}
Besides,
  the equality $\cont(S,S') = \min\{\mathrm{o}(S),\mathrm{o}(S')\}$ is true
  for any pair of distinct branches $S,S' \in   \sep_{0}(\F)$, unless they  belong to a  dicritical family
  of  separatrices of  items  (3) and (5) or,
up to switching their roles, $S$ belongs to one of the dicritical families of singular separatrices of items (4) and (6) and $S'$ is
either another element of this family or its only semi-root.

We say that $\sep_{\infty}(\F)$ if of \emph{Newton type} if, after switching coordinates
$(x,y) \mapsto (y,x)$, the set of separatrices with horizontal tangent line is of Newton type.
}\end{ddef}

\begin{ddef}
\label{def-Newtontype-oblique}
{\rm
Let $c \in \C^{*}$. We say that $\sep_{c}(\F)$, the set of separatrices of $\F$ having
tangent line of slope $c$, is of  \emph{Newton type}  if, when non-empty, it is composed of one of the following
\begin{enumerate}
\item a single smooth isolated  separatrix;
\item  a dicritical family of separatrices of genus $g=0$;
\item a dicritical family of separatrices of genus $g=1$ and index $k=0$, along with a single
associated smooth semi-root.
\end{enumerate}
}\end{ddef}

In view of definitions \ref{def-Newtontype-horizontal} and  \ref{def-Newtontype-oblique}, we have:

\begin{ddef}
\label{def-Newtontype-separatrices}
{\rm We say that $\sep(\F)$ is a Newton set of separatrices if there are convenient analytic coordinates
$(x,y)$ such that $\sep_{c}(\F)$ is of Newton
type for every   $c \in \C \cup \{\infty\}$.
}\end{ddef}

Let $C$ be a germ of analytic curve at $(\C^{2},0)$. A system of analytic coordinates $(x,y)$ is convenient
for $C$ if the coordinate axes are not branches of $C$. For $c \in \C \cup \{\infty\}$, denote by $C_{c}$ the union
of branches of $C$ whose tangent line has slope $c$ ($c= \infty$ corresponding to the vertical tangent).
We present here an adaptation of the main definition of  \cite{garciabarroso2007,garciabarroso2009}:
\begin{ddef}
\label{def-Newtoncurve}
{\rm In the above context, we say that $C$ is a \emph{Newton curve} if there is a system of convenient analytic coordinates $(x,y)$
such that:
\begin{itemize}
\item the branches in $C_{0}$ are either smooth or singular of genus $g=1$ with maximal horizontal contact and, besides,
 for any pair of branches $S,S' \subset C_{0}$, the equality
$\cont(S,S') = \min\{\mathrm{o}(S),\mathrm{o}(S')\}$ is true;
\item  $C_{\infty}$ satisfies the conditions of $C_{0}$, after the change of coordinates
$(x,y) \mapsto (y,x)$;
\item  $C_{c}$, for $c \in \C^{*}$, when non-empty, is formed by a single smooth branch.
\end{itemize}
}\end{ddef}

Observe that, for a germ of holomorphic foliation $\F$,
if  $\sep_{c}(\F)$ is of Newton type and contains only non-dicritical separatrices, then we only
have   cases (1) and (2) in Definition
\ref{def-Newtontype-horizontal}, which automatically satisfy the contact condition without exception,
 together with  case (1)  in Definition
\ref{def-Newtontype-oblique}. Therefore, if $\F$ is a  non-dicritical foliation
and has a Newton set of separatrices (Definition \ref{def-Newtontype-separatrices}), then $C=\Sep(\F)$ is
a Newton curve (Definition \ref{def-Newtoncurve}).

The main result in \cite{garciabarroso2007} (Theorem 3.4) can be stated in the following way:

\begin{teo}
\label{teo-newton-curve} A germ of analytic curve at $(\C^{2},0)$ is a Newton curve  if and only
if it is Newton non-degenerate.
\end{teo}

As already mentioned, we intend to
characterize  Newton non-degenerate foliations in terms of their sets of separatrices. In one
 direction, we have:

\begin{teoRM}
\label{teo-newton-curve1}
Suppose that $\F$ is a holomorphic foliation at $(\C^{2},0)$,   Newton non-degenerate
with respect to a system of convenient analytic coordinates $(x,y)$. Then $\sep(\F)$ is a Newton
set of separatrices.
\end{teoRM}
\begin{proof}
The set of separatrices $\sep_{c}(\F)$ for $c \in \C^{*}$ is described by  Propositions  \ref{prop-nondeg-smooth}  and \ref{prop-nondeg-oblique}.
Let us then consider $\sep_{0}(\F)$, the set of separatrices with horizontal tangent.
By Proposition \ref{prop-nondeg-g2},  $\sep_{0}(\F)$ contains no branch of genus $g > 2$. Also, the branches of genus $g=2$ have maximal horizontal
contact and are in dicritical families
of index $k =0$ ((6) in Definition \ref{def-Newtontype-horizontal}).  Each of these dicritical families has a single  semi-root of genus $g=1$.
In particular, we get that there are no isolated separatrices of genus  $g \geq 2$.
By Corollary \ref{cor-nondeg-oblique1}, separatrices of genus $g=1$ without maximal horizontal  contact are in dicritical families of genus $g=1$ and index $k=0$, and each one of these families has a single associated semi-root
((4) in Definition \ref{def-Newtontype-horizontal}).
As a consequence, we obtain that isolated separatrices of genus $g=1$ have maximal horizontal contact and we obtain  items (1) and (2) in Definition \ref{def-Newtontype-horizontal}.
 On the other hand, by Proposition \ref{prop-flat-g1},  dicritical separatrices of genus $g=1$ with maximal horizontal contact
are in dicritical families of separatrices of genus $g=1$ and any index $k$, yielding  (5).
Evidently, aside from the cases listed above, we can still have dicritical families of separatrices of genus $g=0$, completing
the list of possibilities with item (3).

In order to finish the proof that $\sep_{0}(\F)$ is of Newton type, we have to check the condition on the
contact of pairs of separatrices. To this end, let $S,S' \in \sep_{0}(\F)$ be such that
$\cont(S,S') > \min\{\mathrm{o}(S),\mathrm{o}(S')\}$. Remark that this condition implies, in particular, that $\mathrm{o}(S) = \mathrm{o}(S')$
and that there are  Puiseux series of  $S$ and $S'$  with the same initial term.
If $S$ is singular and has maximal horizontal contact, then this condition implies that also $S'$ is singular and has maximal horizontal contact.
If one of the separatrices has genus $g=2$, say $g(S) = 2$, then,  by Proposition \ref{prop-nondeg-g2}, $S$ belongs to
a dicritical family of separatrices $\mathrm{F}$ and,   either $g(S') = 2$ and $S' \in \mathrm{F}$ or
$g(S') = 1$ and $S'$ is the only   semi-root of $\mathrm{F}$.
On the other hand, if $g(S) = g(S') = 1$, then, by Proposition \ref{prop-isolated-contact}, at least one
of the separatrices  is dicritical.
 Proposition \ref{prop-flat-g1} then gives that both separatrices lie in the same dicritical family
of separatrices of index $k>0$.
If $S$ is singular and does not have maximal horizontal contact then, by Corollary  \ref{cor-nondeg-oblique1}, it lies necessarily in a dicritical family of separatrices of genus $g=1$ and index $k=0$.
The same corollary implies that $S'$ is either another separatrix in this dicritical family or its unique    semi-root.
Finally, if both $S$ and $S'$ are smooth,
Corollary \ref{cor-nondeg-smooth} says that they both belong to the same dicritical family of separatrices
of genus $g=0$.
\end{proof}

We finish this section with a concrete example of a Newton non-degenerate foliation having a dicritical family
of separatrices of genus $g=2$. Note that this is a  generalized curve type foliation --- since it is  logarithmic --- whose
reduction of singularities is not of toric type. Its construction was done by means of the technique   of logarithmic models of foliations \cite{mol2021}.
\begin{example}
\label{ex-dicfamily-g2}
{\rm Let $\F$ be the germ of holomorphic foliation defined by the $1$-form
\begin{eqnarray*}
  \omega &=& (y^{2} - x^{3})(y-x^{2}) \left( 6 \, \frac{\dd (y^{2} - x^{3})}{y^{2} - x^{3}} - 13 \, \frac{\dd ( y-x^{2})}{y-x^{2}} \right) \medskip \\
   &=& (-8 x^4 - 18 x^2y + 26xy^2) \dd x  + (13 x^3  -12x^2y - y^2) \dd y .
\end{eqnarray*}
$\F$ has a dicritical family of separatrices of genus $g=2$ and index $k=0$, having exponents
$\E = \{\frac{3}{2},\frac{7}{4}\}$, along with two isolated separatrices, $S_{1}$ of equation $y^{2} - x^{3}=0$,
which is a semi-root of this dicritical family, and $S_{2}$ of equation $y-x^{2}=0$.
Note that, since the coordinate axes are not $\F$-invariant, the coordinates are convenient.
The Newton polygon $\Delta_{\omega}$ has two edges: $E_{\frac{3}{2}}$ and $E_{2}$, of co-slopes $\kappa = \frac{3}{2}$
and $\kappa = 2$  (see Figure~\ref{fig-example-g2}) .
The reduction of singularities of $\F$ is achieved by the sequence of 5 blow ups that desingularizes
any singular branch of genus $g=2$ and characteristic exponents $\E = \{\frac{3}{2},\frac{7}{4}\}$, having
maximal horizontal contact.

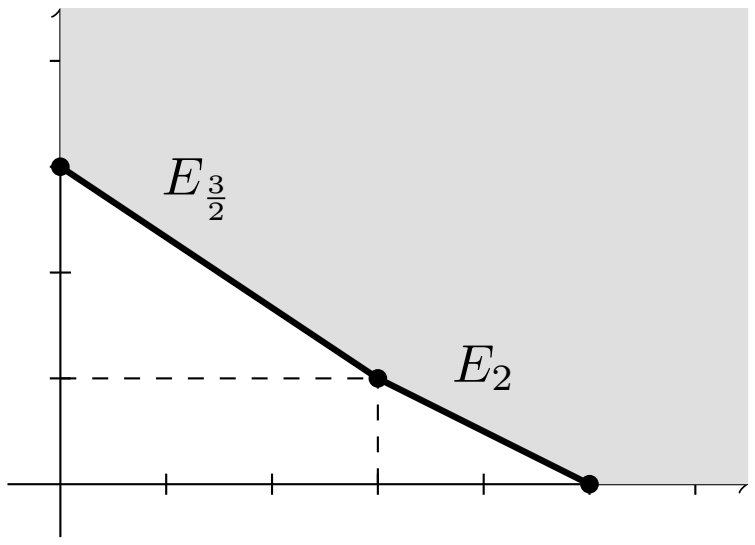
\begin{figure}[h]
\centering
\begin{tikzpicture}[scale=0.7]
  \draw[->] (-0.5,0) -- (6.5,0);
  \draw[->] (0,-0.5) -- (0,4.5);
  \foreach \x in {1,...,6} \draw (\x,-0.1) -- (\x,0.1);
  \foreach \y in {1,...,4} \draw (-0.1,\y) -- (0.1,\y);

  \fill[gray!25] (0,3) -- (0,4.5) -- (6.5,4.5) -- (6.5,0) -- (5,0)
    -- (3,1) -- (0,3) -- cycle;

  \draw[dashed] (0,1) -- (3,1);
  \draw[dashed] (3,0) -- (3,1);

  \draw[very thick] (0,3) -- (3,1) -- (5,0);

  \fill (0,3) circle (2.5pt);
  \fill (3,1) circle (2.5pt);
  \fill (5,0) circle (2.5pt);

  \node[above] at (1.3,2.3) {\small $E_{\frac{3}{2}}$};
  \node[above] at (4,0.7) {\small $E_{2}$};
\end{tikzpicture}
\caption{Newton polygon of the foliation in Example \ref{ex-dicfamily-g2}}
\label{fig-example-g2}
\end{figure}

The degeneracy equation are, for  $E_{\frac{3}{2}}$,
\[ (xA)_{E_{\frac{3}{2}}} = -18x^3y  \qquad \text{and}   \qquad  (yB)_{E_{\frac{3}{2}}} = - y^3+13x^3y \]
and, for $E_{2}$,
\[ (xA)_{E_{2}} = -8x^5-18x^3y \qquad \text{and}   \qquad  (yB)_{E_{2}} = 13 x^3y. \]
In both case, there are no non-trivial common factors. Hence, $\Delta_{\omega}$ is Newton non-degenerate.

}\end{example}


\section{Generalized curve type foliations}
\label{sec-generalizecurve-nondeg}

The   objective of this section is to establish a result
converse to Theorem \ref{teo-newton-curve1}.
We have the following question:
given a germ of foliation $\F$  and a   system of convenient coordinates $(x,y)$  with respect to which  $\sep(\F)$ is a Newton set of
separatrices,  can we conclude that $\F$ is Newton non-degenerate?
Evidently, taking into account Theorem \ref{th-weak}, the answer to this question is negative if $\sep(\F)$ contains a weak separatrix, which is equivalent to asking  the reduction of singularities of $\F$ to admit a well-oriented saddle-node singularity. Thus, for our purposes, some condition must be imposed on the simple models of $\F$. We will see below in Theorem \ref{teo-newton-curve2} that the answer
to our question will be positive if $\F$  is asked to be  of generalized curve type.

We need the following definition
 given in \cite{mariana2025}:

\begin{ddef}
\label{def-effective}
{\rm
 A germ of holomorphic foliation $\F$ at $(\C^2,0)$ is called \emph{effective} either if it is non-dicritical or
if it is dicritical admitting, in the exceptional divisor of its reduction of singularities, dicritical components of
valence at most two.
}\end{ddef}
The name \emph{effective}  for these foliations derives from the fact that they admit effective balanced divisors of separatrices, formed by all isolated separatrices and by one separatrix associated with each dicritical component of valence one.
Many arguments involving non-dicritical foliations are based on the comparison with  a ``reference    foliation'' defined by the level sets of a reduced equation of its curve of   separatrices. This kind of resource  is  ineffective in the dicritical case, considering that the number of separatrices is infinite.
When the foliation is dicritical and effective, however, an effective balanced divisor  of separatrices can play the role of a curve of reference.
Illustrating this fact, we cite the following theorem, which is one of the main results of \cite{mariana2025}:

\begin{teo}
\label{teo-equal-polygons} Let $\F$ be a germ of effective holomorphic foliation at $(\C^{2},0)$ and $\B$ be an  effective balanced divisor of separatrices. Then, for any   $\B$-compatible system of coordinates, $\F$ is a second type foliation if and only if
the Newton polygons of $\F$ and  $\B$ coincide.
\end{teo}

A system of analytic coordinates $(x,y)$ is \emph{$\B$-compatible} if, given
$\pi: (M,\D) \to (\C^{2},0)$  a reduction of singularities of $\F$,
if $p \in \D$ is the base point of the strict transform of one of the coordinate axes,
then either $p$ lies in a
non-dicritical component of $\D$ or it is the base point of some branch of $\pi^{*}\B$.
The existence of a Newton set of separatrices implies that the foliation is effective:

\begin{prop}
\label{prop-B-compatible}
Suppose that the germ of holomorphic foliation $\F$ at $(\C^{2},0)$ has a Newton set of separatrices, for some system of convenient analytic coordinates $(x,y)$.
Then $\F$ is an  effective foliation and any effective balanced divisor of separatrices $\B$ is a Newton curve.
Moreover, $\B$ can be chosen in such a way that the coordinates
$(x,y)$ are  $\B$-compatible.
\end{prop}
\begin{proof} Let $\pi: (M,\D) \to (\C^{2},0)$ be the reduction of singularities of $\F$.
Since $\sep(\F)$ is a Newton set of separatrices, a dicritical component  $D \subset \D$ comes from dicritical families
of separatrices of genus $g = 0$ and $g=1$ of any index $k$,  or $g=2$ and index $k=0$.
For $g=0$ and $g=1$ with $k>0$, we have $\val(D) = 1$. For the other cases, $g=1$ or $2$ with $k=0$,
 we have $\val(D) = 2$.
We thus have that $\F$ is an effective foliation.

If $\B$ is an effective balanced divisor of separatrices for $\F$, then the curve formed by isolated separatrices of $\sep(\F)$
 is a Newton curve. This is explained by  the same argument given  right below Definition \ref{def-Newtoncurve}.
The construction of $\B$ is then completed by the addition of one separatrix   to each dicritical component of valence one.
If one of these separatrices has  slope   $c \in \C^{*}$, then, according to Definition \ref{def-Newtontype-oblique}, it has genus $g=0$ and it is the only element of $\sep_{c}(\F)$ in $\B$.
Thus, the contact condition of Definition
\ref{def-Newtoncurve} is true for the pair formed by this separatrix and any other branch of $\B$.
On the other hand, for a pair of branches in $\B$ of  slope $c=0$,  the contact condition  is also valid,  since
we do not fall in the exceptional cases of Definition \ref{def-Newtontype-horizontal}.
Indeed, two different dicritical separatrices of $\B$ belong to different  dicritical families and, besides, since the dicritical
components of $\D$ associated with  elements of $\B$ have valence one, the corresponding dicritical families of separatrices
have genus $g=0$ or genus $g=1$ and index $k>0$ and, thus,  are devoid of associated semi-roots.

Concerning  $\B$-compatibility, remark that the
hypothesis of $\sep_{0}(\F)$ being of Newton type implies that if
  the base point of the strict transform of  the $x$-axis lies in a dicritical component $D \subset \D$, then
$\val(D) = 1$. This comes from the fact that the separatrices cutting $D$ form a dicritical family of genus $g=0$,
which always has index $k\geq1$, hence
$\val(D)=1$ as noted above.
Hence, as the separatrix in $\B$ associated with $D$, we can choose the one whose strict transform contains $p$.
The same remark applies to $\sep_{\infty}(\F)$. In this way,
 the balanced divisor of separatrices $\B$ can be chosen in such a way that the coordinates $(x,y)$
are $\B$-compatible.
\end{proof}

Foliations of generalized curve type minimize Milnor numbers, once the equisingularity type of its set of separatrices is fixed \cite{camacho1984}. When the foliation $\F$ is non-dicritical,
denoting by $C:f(x,y)=0$ its curve of separatrices,
the statement above is expressed in the following way: $\mu_{0}(C):=\mu_{0}(f) \leq \mu_{0}(\F)$, the two numbers being
equal if and only if $\F$ is of generalized curve type.
The proposition below  is a partial extension of this result to dicritical foliations.

\begin{prop}
\label{prop-milnor-effective} Let $\F$ be a germ of holomorphic foliation at $(\C^{2},0)$. Suppose that $\F$ is effective and that $\B$  is an effective balanced divisor
of separatrices for $\F$.
Then $\mu_{0}(\B) \leq  \mu_{0}(\F)$ and the equality holds if and only if $\F$ is of generalized curve type.
\end{prop}
\begin{proof}
Let $\sigma: (M,D) \to (\C^{2},0)$ be a blow-up with center $0 \in \C^{2}$ and denote  the strict transform foliation by $\tilde{\F} = \sigma^{*} \F$.
Let $p_{1},\ldots,p_{r} \in D$ be set of points formed by all  singular points of $\tilde{\F}$ with respect to $D$ ---
that is, the points where $\tilde{\F}$ and $D$ fail to be transverse ---  together with
the base points of  $\sigma^{*}\B$ on $D$. When $\sigma$ is non-dicritical, these points are the singularities of $\tilde{\F}$ over $D$.
When $\sigma$ is dicritical, since $\F$ is an effective foliation, we have that $r=2$.
In this case,   each of these points can be either a singularity $\tilde{\F}$ with respect to $D$   or,
when the strict transform of  $D$ in the reduction of singularities of $\F$ has valence one, the point  corresponding to
the branch of $\B$ associated with $D$. Order these points in such a way that the irreducible components of $\pi^{*}\B$ intersect
$D$ precisely at $p_{1},\ldots, p_{r'}$, where $r' \leq r$.
For $i=1,\ldots,r$,  denote by $\tilde{\F_{i}}$ the germ of   $\sigma^{*} \F$ at $p_{i}$.
Also, for $i=1,\ldots,r'$, denote by $\tilde{\B}_{i}$ the union of all branches of  $\sigma^{*} \B$ containing $p_{i}$.
Setting
\begin{equation}
\label{eq-balanced-divisor-blowup}
\B(\tilde{\F_{i}}) =
 \begin{cases}  \tilde{\B}_{i} \cup D_{i}   \ \ \text{if $\sigma$ is non-dicritical} \\
\tilde{\B}_{i}  \ \ \text{if $\sigma$ is dicritical},
\end{cases}
\end{equation}
where $D_{i}$ denotes the germ of $D$ at $p_{i}$, then  $\B(\tilde{\F_{i}})$ is an effective balanced divisor of separatrices
for $\tilde{\F_{i}}$ at $p_{i}$ \cite[Lem. 3.11]{genzmer2018}.

Let $m = \ma_{0}(\F)$. The following well-known formulas allow the calculation of the Milnor number of $\F$ from the blow-up data (see \cite[Prop. 4.13]{cano-cerveau-deserti2013}):

\begin{equation}
\label{eq-milnor-blowup}
{\displaystyle
 \mu_{0}(\F) =
 \begin{cases}  m^{2} - (m+1) + \sum_{i=1}^{r } \mu_{p_{i}}(\tilde{\F}_{i})  \ \ \text{if $\sigma$ is non-dicritical} \medskip \\
m^2 + m -1 + \sum_{i=1}^{r } \mu_{p_{i}}(\tilde{\F}_{i}) \ \ \text{if $\sigma$ is dicritical}.
\end{cases}
}\end{equation}

By \cite{camacho1984}, we have that $\mu_{0}(\B) = \mu_{0}(\F_{\B})$, where $\F_{\B}$ is the foliation defined by the level curves of some reduced analytic equation of
 $\B$. This foliation has order $\ma_{0}(\F_{\B}) = m_{\B} -1 $, where $m_{\B} = \ma_{0}(\B)$. Recall from Section
 \ref{sec-holomorphic-foliations} that $m \geq m_{\B} - 1$
 and that the equality holds if and only if $\F$ is a foliation of second type.
At $p_{i}$, for $i=1,\ldots,r'$, the strict transform foliation $\tilde{\F}_{\B} = \sigma^{*} \F_{\B}$  is defined by the level curves of a reduced equation of
$\tilde{\B}_{i} \cup D_{i}$ and, hence, $ \mu_{p_{i}}(\tilde{\F}_{\B}) =\mu_{p_{i}}(\tilde{\B}_{i} \cup D_{i})$.
The non-dicritical part of \eqref{eq-milnor-blowup} applied to $\F_{\B}$   gives:
 \begin{equation}
\label{eq-mu-balanced-blowup11}
 \mu_{0}(\B)   = (m_{\B} -1)^{2} - ((m_{\B} -1) + 1) + \sum_{i=1}^{r' } \mu_{p_{i}}(\tilde{\B}_{i} \cup D_{i}) .
\end{equation}
On the other hand, by \cite[Prop. 6.4.4]{casas2000},
\[ \mu_{p_{i}}(\tilde{\B}_{i} \cup D_{i})  =
 \mu_{p_{i}}(\tilde{\B}_{i})   +
2(\tilde{\B}_{i}, D_{i})_{p_{i}} -1.\]
Summing this over all points $p_{i}$, we have
\begin{equation}
\label{eq-mu-balanced-blowup2} \sum_{i=1}^{r' } \mu_{p_{i}}(\tilde{\B}_{i} \cup D_{i})  =  \sum_{i=1}^{r' }  \mu_{p_{i}}(\tilde{\B}_{i}) +  2(\sigma^{*} \B, D) - r'
=  \sum_{i=1}^{r' }  \mu_{p_{i}}(\tilde{\B}_{i}) +  2 m_{\B}  - r' .
\end{equation}
Putting \eqref{eq-mu-balanced-blowup2} in \eqref{eq-mu-balanced-blowup11}, we get
\begin{equation}
\label{eq-mu-balanced-blowup}
 \mu_{0}(\B) =   m_{\B} ^{2} -  m_{\B} + 1   - r'  + \sum_{i=1}^{r'}  \mu_{p_{i}}(\tilde{\B}_{i}).
\end{equation}

The proof of the proposition then goes by induction in the   number $n$ of blow-ups of a minimal reduction of singularities of $\F$.
If $n=0$ and $\F$ is a regular foliation then $\B$ is formed by the local leaf of $\F$ through the origin.
Hence, $\mu_{0}(\F) = \mu_{0}(\B) = 0$ and there is nothing to prove. On the other hand, if $n=0$ and $\F$ has a simple singularity at $0 \in \C^2$, then $\B$ is formed by two transverse separatrices and, so,  $\mu_{0}(\B) = 1$. We evidently have
 $\mu_{0}(\F) \geq 1$  and the equality holds if and only if $\F$ is a non-degenerate singularity, giving the result in this case.

Now, let us suppose that the minimal reduction of singularities of $\F$ has $n \geq 1$ blow-ups and that the proposition's statement
is true for foliations having reductions of singularities composed by less than $n$ blow-ups.
We apply the first blow-up $\sigma: (M,D) \to (\C^{2},0)$ at $0 \in \C^{2}$. The analysis then splits into two cases:

\par \noindent 1. \underline{$\sigma$ is non-dicritical.}

In this case,
applying  \eqref{eq-milnor-blowup}, \eqref{eq-mu-balanced-blowup11}
and the first line of  \eqref{eq-balanced-divisor-blowup},  we find
\[ \mu_{0}(\F) =  m^{2} - (m+1) + \sum_{i=1}^{r } \mu_{p_{i}}(\tilde{\F}_{i})
\geq
(m_{\B}-1)^{2} - ((m_{\B}-1) + 1) + \sum_{i=1}^{r'} \mu_{p_{i}}(\tilde{\B}_{i} \cup D_{i}) =  \mu_{0}(\B) \]
The inequality comes from the
induction hypothesis and from the fact that $m \geq m_{\B} - 1$.
If, in this  expression, the equality is true, then $r' =r$ and  $\mu_{p_{i}}(\tilde{\F}_{i}) = \mu_{p_{i}}(\tilde{\B}_{i} \cup D_{i})$ for every $i=1,\ldots,r$. Then,
the induction hypothesis again gives that  each $\tilde{\F}_{i}$ is a generalized curve foliation and so is $\F$. Conversely, if $\F$ is a generalized curve foliation, then  so is each  $\tilde{\F}_{i}$ and, also, there must be a branch of $\pi^{*} \B$ containing $p_{i}$ for $i=1,\ldots,r$
\cite[Lem. 4]{mol2002}.
Thus   $r' = r$, $\mu_{p_{i}}(\tilde{\F}_{i}) = \mu_{p_{i}}(\tilde{\B}_{i} \cup D_{i})$ for every $i=1,\ldots,r$ by induction and, also, $m = m_{\B} - 1$, since
$\F$ is in particular a foliation of second type. As a consequence, equality holds in the expression above and  $\mu_{0}(\F) = \mu_{0}(\B)$
as wished.

\par \noindent 2. \underline{$\sigma$ is dicritical.}
By hypothesis, $\F$ is  an effective foliation, thus, by definition, the valence of the strict transform of $D$ in the reduction of singularities of $\F$
is one or two. In the case of valence one, a dicritical separatrix of $D$ is attached to the balanced divisor of separatrices. In both cases, $r=2$. In order to simplify the exposition, we also assume that  $r' = 2$. Thus, formula \eqref{eq-mu-balanced-blowup} becomes
\begin{equation}
\label{eq-mu-balanced-blowup-r2}
 \mu_{0}(\B) = \mu_{0}(\F_{\B}) =  m_{\B} ^{2} -  m_{\B} -1 + \sum_{i=1}^{2 }  \mu_{p_{i}}(\tilde{\B}_{i}).
\end{equation}
Now, the second line of  \eqref{eq-balanced-divisor-blowup} gives that $\tilde{\B}_{i}$ is an effective balanced divisor of separatrices for $\tilde{\F}_{i}$ at each $p_{i}$. As an induction hypothesis, $\mu_{p_{i}}(\tilde{\B}_{i}) \leq \mu_{p_{i}}(\tilde{\F}_{i})$ for
each $i=1,2$, with equality holding if and only if $\tilde{\F}_{i}$ is a germ of generalized curve foliation at $p_{i}$.
Next, we apply formula \eqref{eq-milnor-blowup}, obtaining,  from the
induction hypothesis, from   $m \geq m_{\B}  - 1$ and also from \eqref{eq-mu-balanced-blowup-r2}, the   inequality in the proposition's statement:
\[ \mu_{0}(\F) = m^2 + m -1 + \sum_{i=1}^{2 } \mu_{p_{i}}(\tilde{\F}_{i}) \geq
(m_{\B}-1)^{2} + (m_{\B}-1) - 1 + \sum_{i=1}^{2 } \mu_{p_{i}}(\tilde{\B}_{i}) \]
\begin{equation}
\label{eq-inequality-milnor} = m_{\B}^{2} - m_{\B} -1 + \sum_{i=1}^{2 } \mu_{p_{i}}(\tilde{\B}_{i}) = \mu_{0}(\B) .
\end{equation}
If, in this expression, the equality is true, then  $\mu_{p_{i}}(\tilde{\F}_{i}) = \mu_{p_{i}}(\tilde{\B}_{i})$ for   $i=1,2$ and, by
the induction hypothesis,    each $\tilde{\F}_{i}$ is a generalized curve foliation and so is $\F$. Conversely, if $\F$ is a generalized curve foliation, then so is each  $\tilde{\F}_{i}$  and the equality involving their Milnor numbers is true and, also, $m= m_{\B}  - 1$, since
$\F$ is in particular a foliation of second type. Thus, equality holds in \eqref{eq-inequality-milnor} and  $\mu_{0}(\F) = \mu_{0}(\B)$
as wished.
\end{proof}

We can then prove that, for foliations of generalized curve type, a result converse to Theorem \ref{teo-newton-curve1} is true:

\begin{teoRM}
\label{teo-newton-curve2}
Let $\F$ be a holomorphic foliation at $(\C^{2},0)$ such that $\sep(\F)$ is a Newton set of separatrices with respect to some
system of convenient analytic coordinates.
If $\F$ is of generalized curve type then it is Newton non-degenerate.
\end{teoRM}
\begin{proof} Since $\sep(\F)$ is a Newton set of separatrices, by Proposition \ref{prop-B-compatible}, $\F$ is an effective foliation
 having an effective balanced divisor of separatrices $\B$ which
can be chosen making the coordinates in question $\B$-compatible.
By   Theorem \ref{teo-equal-polygons}, being $\F$ an effective foliation of generalized curve type,
the Newton polygons of $\F$ and  $\B$ coincide and, so,  their Newton numbers are equal:
$\nu_{0}(\F) = \nu_{0}(\B)$.
On the other hand, also by  Proposition \ref{prop-B-compatible}, $\B$ is a Newton curve (Definition \ref{def-Newtoncurve}) and
therefore,
 by Theorem \ref{teo-newton-curve},   it is Newton non-degenerate as a curve.
Corollary \ref{cor-mu-nu-curves} then implies that $\mu_{0}(\B) = \nu_{0}(\B)$.
Using again the fact that $\F$ is an effective foliation of generalized curve type,
Proposition \ref{prop-milnor-effective} gives  $\mu_{0}(\F) = \mu_{0}(\B)$.
This chain of equalities  gives   $\mu_{0}(\F) = \nu_{0}(\F)$ and it follows from Theorem \ref{teo-fol-nondeg2} that the foliation $\F$ is Newton non-degenerate.
\end{proof}

 \begin{example}
\label{ex-newton-set-not-enough}
{\rm
We retake Example \ref{ex-newton-deg} to show that the hypothesis
of $\F$ being of generalized curve type cannot be dropped in
Theorem \ref{teo-newton-curve2}. Recall that, in that example,
$0\in\C^2$ is a non-dicritical singularity, and the reduction of
singularities has two well-oriented saddle-node singularities ---
hence $\F$ is not of generalized curve type --- together with a
non-degenerate singularity, and that $\Delta_\omega$ is Newton
degenerate.

Searching for invariant lines $y=cx$ of $\F$, the condition
$A(x,cx)+cB(x,cx)\equiv0$ factors as
\[
x^{2}(c-1)(c+2)(2c+1) \equiv 0,
\]
giving exactly three invariant lines: $y=x$ and $y=-\tfrac12 x$
(the tangent cones of the two weak separatrices already
identified) and a third one, $y=-2x$, coming from the
non-degenerate singularity at $t=\infty$. Thus
$\Sep(\F)=\{y=x,\,y=-\tfrac12x,\,y=-2x\}$ consists of three smooth
isolated separatrices with distinct oblique tangent lines, and no
other direction contains a separatrix.

Consequently, $\sep(\F)$ is trivially a Newton set of separatrices
(Definition \ref{def-Newtontype-separatrices}), even though $\F$
is Newton degenerate. This shows that having a Newton set of
separatrices does not by itself force Newton non-degeneracy: the
generalized curve type hypothesis in Theorem
\ref{teo-newton-curve2} is essential.
}\end{example}

We finish by proving that, for foliations of generalized curve type, the property of Newton non-degeneracy
is a topological invariant. We say that two germs of holomorphic foliations $\F, \hat{\F}$ at   $(\C^{2},0)$ are \emph{topologically
equivalent} if there is a germ of homeomorphism $h : (\C^{2},0) \to (\C^{2},0)$ that, when realized in sufficiently small
neighborhoods of $0 \in \C^2$, sends leaves of $\F$ into leaves of $\hat{\F}$.
We have:

\begin{teoRM}
\label{teo-topological}
Let $\F,\hat{\F}$ be topologically equivalent holomorphic foliations at $(\C^{2},0)$.
Suppose that $\F$ is of generalized curve type and Newton non-degenerate with respect
to convenient analytic coordinates $(x,y)$. Then $\hat{\F}$ is of generalized curve type and there are
analytic coordinates with respect to which $\hat{\F}$ is Newton non-degenerate.
\end{teoRM}
\begin{proof} It is well known that the family of generalized curve type foliations is stable
under topological equivalences \cite{camacho1984}, giving at once that $\hat{\F}$ is of generalized curve type.
Since $\F$ is Newton non-degenerate with respect to the coordinates $(x,y)$, by Theorem
\ref{teo-newton-curve1}, $\sep(\F)$ is a Newton set of separatrices. We have:

\begin{assertion*}
 $\sep(\hat{\F})$ is a Newton set of separatrices with respect to some system
of analytic coordinates $(\hat{x},\hat{y})$.
\end{assertion*}
\noindent  \emph{Proof {\rm (of the Assertion)}. }
As consequence of  Zariski's equisingularity theorem
\cite{casas2000}, all properties involved in Definitions \ref{def-Newtontype-horizontal} and \ref{def-Newtontype-oblique}
are preserved by topological equivalences. Indeed, denoting by $h$ the germ of homeomorphism conjugating $\F$ and $\hat{\F}$,
we have:
\begin{itemize}
  \item   $S \in \Sep(\F)$ and $\hat{S} = h(S) \in \sep(\hat{\F})$ are equisingular (same genus, same
   characteristic exponents).
  \item If $S, S' \in \Sep(\F)$ and  $\hat{S} = h(S), \hat{S}' = h(S') \in \sep(\hat{\F})$, then
  $\cont(S,S') = \cont(\hat{S},\hat{S}')$.
  \item Separatrices of $\F$ with the same tangent line will be taken in separatrices of $\hat{\F}$ with
  the same tangent line.
  \item A dicritical family of separatrices of $\F$ is taken into a dicritical family of separatrices of $\hat{\F}$  with the same genus $g$ and the same index $k$.
\end{itemize}
It remains to establish how to choose analytic coordinates  $(\hat{x},\hat{y})$ and, for this purpose,
the condition of maximal horizontal contact is essential.
First we take, among all  elements of   $\Sep_{0}(\F)$, a separatrix $S$ whose contact  $\varepsilon \in \Q_{>0}$ with the $x$-axis is maximal
and set $\hat{S} = h(S) \in \Sep(\hat{\F})$. Next,  take any smooth branch $\hat{S}_{0}$, non-invariant
by $\hat{\F}$, such that  $\cont(\hat{S}_{0}, \hat{S}) = \varepsilon$. If   $\Sep_{0}(\F) = \emptyset$, we can take any smooth branch tangent to the $x$-axis as $\hat{S}_{0}$.
We do the same construction for $\Sep_{\infty}(\F)$, after interchanging the roles of $x$ and $y$, and take
another smooth branch denoted $\hat{S}_{\infty}$. The coordinates $(\hat{x}, \hat{y})$ are then chosen in such a way that $\hat{S}_{0}$ and $\hat{S}_{\infty}$
are, respectively, the horizontal and the vertical axes. This completes the proof of the assertion.
\qed

We thus have that $\hat{\F}$ is a generalized curve type foliation having a Newton set of separatrices with respect to the analytic
coordinates $(\hat{x},\hat{y})$. By Theorem \ref{teo-newton-curve2}, for these coordinates, $\hat{\F}$ is Newton non-degenerate.
\end{proof}

\begin{example}
\label{ex-topological-invariance}
{\rm
We illustrate Theorem \ref{teo-topological}. Consider
\[ \F:\ \omega = (-x^{2}+2xy^{2}-y)\,\dd x + (-2x^{2}y+x+y^{2})\,\dd y \]
\noindent and
\[ \hat{\F}:\ \hat{\omega} = (-4X^{2}Y-4XY+4Y^{3}+2Y)\,\dd X
+ (4X^{3}-2X^{2}-4XY^{2}-2X-2Y^{2})\,\dd Y. \]
The linear map $h(X,Y)=(X+Y,X-Y)$ satisfies $h^{*}\omega=\hat{\omega}$,
so $\F,\hat{\F}$ are analytically, hence topologically,
equivalent.
For $\F$: $A(x,0)=-x^{2}$ and $B(0,y)=y^{2}$ are non-zero, so
$(x,y)$ are convenient. The Newton polygon $\Delta_{\omega}$ has vertices $(0,3)$,
$(1,1)$, $(3,0)$, and its degeneracy equations,
\[ (xA)_{E_{1/2}}=-xy,\ \ (yB)_{E_{1/2}}=y(x+y^{2}),\qquad
   (xA)_{E_{2}}=-x(y+x^{2}),\ \ (yB)_{E_{2}}=xy, \]
share only the trivial factors $y$, $x$: $\F$ is Newton
non-degenerate.

For $\hat{\F}$: $\hat{A}(X,0)\equiv0$, so $Y=0$ is invariant and
$(X,Y)$ are \emph{not} convenient, so non-degeneracy cannot even
be tested there. Theorem \ref{teo-topological} guarantees
convenient coordinates for $\hat{\F}$; we exhibit them via the
construction in its proof.
The functions $f_{1}=y-x^{2}$ and $f_{2}=x-y^{2}$ define
separatrices $S_{1}:y=x^2$ (horizontal tangent, $\oo(S_1)=2$) and
$S_{2}:x=y^2$ (vertical tangent, $\oo(S_2)=2$) of $\F$, as one
checks from $v(f_i)\equiv0 \bmod f_i$ for the dual field
$v=B\partial_x-A\partial_y$.
Set $\hat{S}=h(S_1)$: its Puiseux series is $Y=X-4X^{2}+\cdots$, no longer
horizontal. Looking for $\hat{S}_0$ non-invariant with
$\cont(\hat{S}_0,\hat{S})=\oo(S_1)=2$: the line $Y=X$ works
($\hat{S}-\hat{S}_0=-4X^2+\cdots$, order $2$; non-invariance gives
remainder $8X^2\neq0$). Analogously, $\hat{S}'=h(S_2)$ gives
$X=-Y+4Y^2+\cdots$, and $\hat{S}_\infty:Y=-X$ satisfies
$\cont(\hat{S}_\infty,\hat{S}')=2$ and is non-invariant
(remainder $-8X^2$).
Since $Y=X$ and $Y=-X$ are transverse, they give the new
convenient coordinates $(\hat{x},\hat{y})$ for $\hat{\F}$
promised by the theorem.
}\end{example}

\end{document}